\RequirePackage{fix-cm}

\documentclass[leqno, fleqn, centertags, 12pt]{article}

\usepackage{latexsym}
\usepackage{amsmath}
\usepackage{amssymb}
\usepackage{verbatim}

\usepackage[T1]{fontenc}

\usepackage[upright, widespace]{fourier}

\usepackage{bm}

\usepackage{fullpage}

\usepackage{tikz-cd}

\tikzcdset{row sep/special/.initial=12ex}
\tikzcdset{row sep/spe/.initial=10.5ex}
\tikzcdset{row sep/speh/.initial=13ex}
\tikzcdset{row sep/speh/.initial=13ex}
\tikzcdset{row sep/sma/.initial=1.2em}

\tikzcdset{column sep/spec/.initial=5em}
\tikzcdset{column sep/spech/.initial=6em}
\tikzcdset{column sep/speci/.initial=6em}
\tikzcdset{column sep/specc/.initial=8em}
\tikzcdset{column sep/speccc/.initial=10em}
\tikzcdset{column sep/sboom/.initial=5.6em}

\tikzcdset{row sep/red/.initial=3.79em}

\tikzcdset{row sep/tri/.initial=5em}
\tikzcdset{column sep/tri/.initial=2.4em}

\tikzcdset{row sep/boom/.initial=5em}
\tikzcdset{column sep/boom/.initial=5.1em}

\tikzcdset{row sep/sboom/.initial=5.6em}
\tikzcdset{column sep/sboom/.initial=5.712em}

\tikzcdset{row sep/qboom/.initial=6.4em}
\tikzcdset{column sep/qboom/.initial=6.528em}
\tikzcdset{row sep/hqboom/.initial=1.65em}

\tikzcdset{row sep/boomm/.initial=5.6em}
\tikzcdset{column sep/boomm/.initial=4.25em}

\tikzcdset{row sep/bboom/.initial=8em}
\tikzcdset{column sep/bboom/.initial=8.16em}

\tikzcdset{row sep/booms/.initial=4.25em}
\tikzcdset{column sep/booms/.initial=4.335em}

\tikzcdset{row sep/boomss/.initial=4.75em}
\tikzcdset{column sep/boomss/.initial=4.845em}

\tikzcdset{row sep/free/.initial=12ex}
\tikzcdset{column sep/freee/.initial=4em}

\usepackage{comment}

\newskip\nineskipamount \nineskipamount=9pt plus 0pt minus 0pt
\newskip\zeroskipamount \zeroskipamount=0pt plus 0pt minus 0pt
\usepackage[noorphans,vskip=\nineskipamount]{quoting}

\makeatletter
\renewcommand{\@makefntext}[1]{\vspace*{0.5ex}\parindent=0em
\hspace*{-0.4em}
\hbox to 0.4em{\hss\@makefnmark}\hspace*{0.4em}{#1}
}
\makeatother

\newcounter{mysectionnumber}
\newcommand{\mysection}[2]{\setcounter{footnote}{0}
\setcounter{myparnum}{0}
\refstepcounter{mysectionnumber}
\vspace{27pt}{\Large {\themysectionnumber.} {#1}}\label{#2}\vspace*{15pt}}

\newcommand{\myuppar}[1]{\vspace{\medskipamount}\textbf{#1}\hspace*{0.5em}}

\newcommand{\myit}[1]{\textbf{\textit{#1}}\hspace{0.0em}}

\newcounter{myparnum}[mysectionnumber]
\renewcommand{\themyparnum}{\arabic{mysectionnumber}.\arabic{myparnum}}
\newcommand{\mypar}[2]{\refstepcounter{myparnum}{\vspace{\medskipamount}\textbf{{\themyparnum. #1}\label{#2}}\hspace{0.5em}}}

\newcounter{mylemmanum}[myparnum]
\newcounter{myappendnumber}
\newcounter{myaparnum}[myappendnumber]
\newcommand{\myappend}[2]{\setcounter{footnote}{0}
\setcounter{myaparnum}{0}
\setcounter{myparnum}{0}
\refstepcounter{myappendnumber}
\vspace{27pt}{\Large A\dff. 
{#1}}\label{#2}\vspace*{15pt}}

\newcommand{\myapar}[2]{\refstepcounter{myaparnum}{\vspace{\medskipamount}\textbf{{\themyaparnum. #1
}\label{#2}}\hspace{0.5em}}}
\renewcommand{\themyaparnum}{A\halfff\fff.\fff
\arabic{myaparnum}}

\newcounter{myapparnum}[mysectionnumber]
\newcommand{\proof}{\vspace{\medskipamount}{\textbf{{\emph{Proof}.}}\hspace*{1em}}}

\newcommand{\eproof}{ $\blacksquare$}

\newcommand{\dis}{\displaystyle}

\def\sss{\hspace{0.05em}\ }
\def\dss{\hspace{0.1em}\ }
\def\trs{\hspace{0.15em}\ }
\def\qss{\hspace{0.2em}\ }
\def\pss{\hspace{0.3em}\ }
\def\oss{\hspace{0.4em}\ }

\def\halfff{\hspace*{0.025em}}
\def\fff{\hspace*{0.05em}}
\def\dff{\hspace*{0.1em}}
\def\trf{\hspace*{0.15em}}
\def\qff{\hspace*{0.2em}}
\def\pff{\hspace*{0.3em}}
\def\off{\hspace*{0.4em}}

\newcommand{\nsp}{\hspace*{-0.1em}}
\newcommand{\nnsp}{\hspace*{-0.15em}}

\newcommand{\dnsp}{\hspace*{-0.2em}}

\renewcommand{\leq}{\leqslant}
\renewcommand{\geq}{\geqslant}

\newcommand{\zzz}{\mathbf{Z}}

\newcommand{\rrr}{\mathbf{R}}

\newcommand{\image}{\operatorname{Im}\trf}
\newcommand{\kernel}{\operatorname{Ker}\trf}
\newcommand{\coker}{\operatorname{coker}\trf}
\newcommand{\pr}{\operatorname{p{\fff}r}}

\newcommand{\ssk}{\operatorname{S{\fff}k}}
\newcommand{\id}{\operatorname{id}}

\newcommand{\gr}{\operatorname{G{\fff}r}}

\newcommand{\ob}{\operatorname{Ob}}

\newcommand{\ske}{\operatorname{Sk}}
\newcommand{\map}{\operatorname{Map}}
\newcommand{\fin}{\operatorname{fin}}

\newcommand{\num}[1]{|\qff #1 \qff|}
\newcommand{\bnum}[1]{\bigl|\qff #1 \qff\bigr|}
\newcommand{\bbnum}[1]{\llbracket\qff #1 \qff\rrbracket}

\newcommand{\norm}[1]{\|\qff #1 \qff\|}

\newcommand{\ffin}{^{\dff \mathrm{fin}}}

\newcommand{\ttoo}{\hspace*{0.2em}\longrightarrow\hspace*{0.2em}}

\begin{document}

\setlength{\baselineskip}{12pt plus 0pt minus 0pt}
\setlength{\parskip}{12pt plus 0pt minus 0pt}
\setlength{\abovedisplayskip}{12pt plus 0pt minus 0pt}
\setlength{\belowdisplayskip}{12pt plus 0pt minus 0pt}

\newskip\smallskipamount \smallskipamount=3pt plus 0pt minus 0pt
\newskip\medskipamount   \medskipamount  =6pt plus 0pt minus 0pt
\newskip\bigskipamount   \bigskipamount =12pt plus 0pt minus 0pt

\author{Nikolai\qss V.\qss Ivanov}
\title{Topological\qss categories\qss related\qss to\qss
Fredholm\qss operators\fff:\oss
II.\oss The\qss analytical\qss index}
\date{}

\footnotetext{\hspace*{-0.65em}\copyright\oss 
Nikolai\qss V.\qss Ivanov,\oss 2021.\oss}

\maketitle

\renewcommand{\baselinestretch}{1}
\selectfont

\vspace*{6ex}
\vspace*{6ex}

\myit{\hspace*{0em}\large Contents}  \vspace*{1ex} \vspace*{\bigskipamount}\\ 
\hbox to 0.8\textwidth{\myit{\phantom{1}1.}\hspace*{0.5em} Introduction \hfil 2}\hspace*{0.5em} 
\vspace*{0.25ex}\\
\hbox to 0.8\textwidth{\myit{\phantom{1}2.}\hspace*{0.5em} Analytical\dss index of\dss self-adjoint\trs Fredholm\dss families \hfil 8}\hspace*{0.5em} \vspace*{0.25ex}\\
\hbox to 0.8\textwidth{\myit{\phantom{1}3.}\hspace*{0.5em} Analytical\dss index of\trs Fredholm\dss families \hfil 17}\hspace*{0.5em} \vspace*{0.25ex}\\
\hbox to 0.8\textwidth{\myit{\phantom{1}4.}\hspace*{0.5em} Hilbert\dss bundles and strictly\trs Fredholm\dss families \hfil 21}\hspace*{0.5em} \vspace*{0.25ex}\\
\hbox to 0.8\textwidth{\myit{\phantom{1}5.}\hspace*{0.5em} Analytical\dss index of\dss strictly\dss Fredholm\dss families \hfil 27}\hspace*{0.5em} \vspace*{0.25ex}\\
\hbox to 0.8\textwidth{\myit{\phantom{1}6.}\hspace*{0.5em} Discrete-spectrum\dss families \hfil 31}  \hspace*{0.5em} \vspace*{0.25ex}\\
\hbox to 0.8\textwidth{\myit{\phantom{1}7.}\hspace*{0.5em} The classical\dss index of\trs Fredholm\dss families \hfil 37}  \hspace*{0.5em} \vspace*{0.25ex}\\
\hbox to 0.8\textwidth{\myit{\phantom{1}8.}\hspace*{0.5em} The classical\dss index of\dss self-adjoint\trs 
Fredholm\dss families \hfil 42}  \hspace*{0.5em} \vspace*{1.5ex}\\
\hbox to 0.8\textwidth{\myit{Appendix.}\hspace*{0.5em} Polarizations and\sss trivializations of\dss Hilbert\dss bundles\hfil 48}\hspace*{0.5em}
\vspace*{1.5ex}\\
\hbox to 0.8\textwidth{\myit{References}  \hfil 52}\hspace*{0.5em}  \hspace*{0.5em}  \vspace*{0.25ex}

\renewcommand{\baselinestretch}{1}
\selectfont

\vspace*{6ex}

\newpage
\mysection{Introduction}{introduction}

\vspace*{\medskipamount}
\textbf{A naive approach\sss to\sss the analytical\dss index of\trs families 
of\dss self-adjoint\qss Fredholm\qss operators.}
The classical\sss definition of\dss the analytical\dss index of\dss families 
of\dss self-adjoint\trs Fredholm\dss operators,\oss
due\sss to\dss Atiyah\dss and\dss Singer\qss \cite{as},\oss \cite{as5},\oss
is\dss indirect.\oss
Let $H$ be a separable infinite dimensional\dss Hilbert\sss space,\oss
and\sss let\sss $A_{\dff x}\dff,\pff x\qff \in\qff X$\sss 
be a family of\dss self-adjoint\trs Fredholm\dss operators\sss $H\qff \ttoo\qff H$\nnsp.\oss
Atiyah\dss and\dss Singer\dss define\sss the analytical\dss index of\dss
$A_{\dff x}\dff,\pff x\qff \in\qff X$\sss 
as\sss the analytical\dss index of\dss a related\sss family of\dss general\qss
(not\sss assumed\sss to be self-adjoint\fff)\qss
Fredholm\dss operators\sss $H\qff \ttoo\qff H$\sss parameterized\sss by\sss the suspension\sss
$\Sigma\dff X$\sss of\sss $X$\sss and defined\sss by a simple formula.\oss

The index of\dss the\sss family $A_{\dff x}\dff,\pff x\qff \in\qff X$ 
ought\sss to be\sss the family\sss
of\dss kernels\sss $\kernel\dff A_{\dff x}\dff,\pff x\qff \in\qff X$\nnsp,\oss
or,\oss rather,\oss some equivalence class of\dss this family of\dss kernels.\oss
This\sss idea\sss immediately encounters a couple of\dss difficulties.\oss
Namely,\oss the kernels usually do not\sss depend continuously on\sss $x$\nnsp,\oss
and\sss it\dss is\dss unclear\sss what\sss kind of\dss an object\sss the family\sss
$\kernel\dss A_{\dff x}\dff,\pff x\qff \in\qff X$\sss is.\oss
Moreover,\oss if\dss the\sss kernels happen\sss to continuously depend on $x$
and\sss hence\sss to define a vector bundle on\sss $X$\nnsp,\oss
then\sss the family can\sss be deformed\sss to a family\sss with zero kernels
and\dss its\sss index ought\sss to be zero.\oss
In\sss fact,\sss the analytical\dss index of\sss
$A_{\dff x}\dff,\pff x\qff \in\qff X$\sss 
can\sss be\sss thought\sss as a measure of\dss 
the failure of\dss the family of\dss kernels\sss
$\kernel\dff A_{\dff x}\dff,\pff x\qff \in\qff X$\sss
to be a vector\sss bundle.\oss

Nevertheless,\oss there\dss is\dss a definition of\dss the analytical\dss index\sss
of\dss families of\dss self-adjoint\trs Fredholm\dss operators
making\sss this idea\sss rigorous,\oss
and\sss the present\sss paper\dss is\dss devoted\sss to such a definition.\oss
It\dss is\dss based on\sss ideas of\trs Segal\qss \cite{s4}\qss
as developed\sss by\sss the author\qss \cite{i2}.\oss

\myuppar{Enhanced operators.}
The starting\sss point\dss is\dss the notion of\qss
\emph{enhanced\sss operators}\pss introduced\sss in\qss \cite{i2}.\oss
Recall\sss that\sss if\dss $A\dff \colon\dff H\qff \ttoo\qff H$\sss
is\dss a self-adjoint\trs Fredholm\dss operator,\oss
then\sss there exists\sss $\varepsilon\qff >\qff 0$\sss such\sss that\sss
$\varepsilon\fff,\qff -\qff \varepsilon$\sss 
do not\sss belong\sss to\sss
the spectrum\sss 
$\sigma\trf(\trf A\trf)$\sss of\sss $A$\nsp,\oss
the interval\sss 
$[\trf -\qff \varepsilon\fff,\qff \varepsilon\trf]$\sss
is\dss disjoint\sss from\sss the essential\sss spectrum of\sss $A$\nnsp,\oss
and\sss hence\sss the image of\dss the spectral\dss projection\sss
$P_{\dff [\dff -\dff \varepsilon\fff,\qff \varepsilon\dff]}\trf
(\trf A\trf)$\sss
is\dss finitely dimensional.\oss
An\qss \emph{enhanced\dss operator}\pss is\dss defined as a pair\sss
$(\trf A\fff,\qff \varepsilon\trf)$\sss
having\sss these properties.\oss
By a basic property of\dss self-adjoint\trs Fredholm\dss operators,\oss
if\dss $(\trf A\fff,\qff \varepsilon\trf)$\sss is\dss an enhanced operator\sss
and\sss $A'\dff \colon\dff H\qff \ttoo\qff H$\sss is\dss a
self-adjoint\trs Fredholm\dss operator sufficiently close\sss to\sss $A$\sss
in\sss the norm\sss topology,\oss then\sss $(\trf A'\fff,\qff \varepsilon\trf)$\sss
is\dss also an enhanced operator.\oss
Moreover,\oss the spectral\sss projection\sss
$P_{\dff [\dff -\dff \varepsilon\fff,\qff \varepsilon\dff]}\trf
(\trf A'\trf)$ and\sss its image
continuously depend on\sss $A'$ for\sss $A'$ close\sss to\sss $A$\nnsp.\oss

\myuppar{Toward\sss a rigorous form of\dss the naive approach.}
As above,\oss let\sss
$A_{\dff x}\dff,\pff x\qff \in\qff X$\sss
be a family\sss of\dss self-adjoint\trs Fredholm\dss operators in\sss $H$\nnsp.\oss
We would\sss like\sss to\sss turn\sss the family\sss
of\dss kernels\sss $\kernel\dff A_{\dff x}\dff,\pff x\qff \in\qff X$\sss 
into something\sss manageable.\oss
Let\sss $z\qff \in\qff X$\nnsp.\oss
Since\sss $A_{\dff z}$\sss is\dss a self-adjoint\trs Fredholm\dss operator,\oss
there exists\sss 
$\varepsilon\off =\off \varepsilon_{\dff z}\qff >\qff 0$\sss 
such\sss that\sss
$(\trf A\fff,\qff \varepsilon\trf)$\sss
is\dss an enhanced operator.\oss 
Clearly,\oss the kernel\sss 
$\kernel\dff A_{\dff z}$\sss
is\dss contained\sss in\sss this image of\dss the spectral\sss projection\sss
$P_{\dff [\dff -\dff \varepsilon\fff,\qff \varepsilon\dff]}\trf
(\trf A_{\dff z}\trf)$\nnsp,\oss 
and\dss if\dss $\varepsilon$\sss
is\dss sufficiently small,\pss 
$\kernel\dff A_{\dff z}$\sss
is\dss equal\sss to\sss this image.\oss
Even\sss if\dss this\dss is\dss not\sss the case,\oss
one may\sss consider\sss the image\sss 
$\image\dff 
P_{\dff [\dff -\dff \varepsilon\fff,\qff \varepsilon\dff]}\trf
(\trf A_{\dff z}\trf)$\sss
as an approximation\sss to\sss the kernel\sss $\kernel\dff A_{\dff z}$\nsp.\oss
If\dss the family\sss 
$A_{\dff x}\dff,\pff x\qff \in\qff X$\sss
is\dss moderately continuous,\oss 
there\dss is\dss a neighborhood\sss $U_{\fff z}$\sss of\sss $z$\sss
such\sss that\sss $(\trf A_{\dff y}\fff,\qff \varepsilon\trf)$\sss
is\dss an enhanced operator for every\sss 
$y\qff \in\qff U_{\fff z}$ and\sss
the spectral\sss projections\sss 
$P_{\dff [\dff -\dff \varepsilon\fff,\qff \varepsilon\dff]}\trf
(\trf A_{\dff y}\trf)$\sss
continuously depends on\sss $y$\sss
for\sss 
$y\qff \in\qff U_{\fff z}$\nsp.\oss
The family of\dss images\sss 
$\image\dff P_{\dff [\dff -\dff \varepsilon\fff,\qff \varepsilon\dff]}\trf
(\trf A_{\dff y}\trf)$\nnsp,\qss $y\qff \in\qff U_{\fff z}$\sss
can\sss be considered as an approximation\sss to\sss the family of\dss kernels\sss
$\kernel\dff A_{\dff y}$\nsp,\qss $y\qff \in\qff U_{\fff z}$\nsp.\oss
The continuous dependence of\dss
$P_{\dff [\dff -\dff \varepsilon\fff,\qff \varepsilon\dff]}\trf
(\trf A_{\dff y}\trf)$\sss 
on\sss $y\qff \in\qff U_{\fff z}$\sss implies\sss that\sss the family\sss of\dss images\sss
$\image\dff 
P_{\dff [\dff -\dff \varepsilon\fff,\qff \varepsilon\dff]}\trf
(\trf A_{\dff y}\trf)$\nsp,\qss $y\qff \in\qff U_{\fff z}$\sss
defines a finitely dimensional\sss vector bundle over\sss $U_{\fff z}$\nsp.\oss
We will\sss define\sss the analytical\dss index as\sss the homotopy class of\dss a continuous map\sss
build\sss from\sss these vector bundles,\oss
or,\oss rather,\oss directly from\sss the images
$\image\dff 
P_{\dff [\dff -\dff \varepsilon\fff,\qff \varepsilon\dff]}\trf
(\trf A_{\dff y}\trf)$
and\sss the restrictions of\sss $A_{\dff y}$\sss 
on\sss these images.

The main\sss problem\dss is\dss in\sss turning\sss a collection of\dss
locally defined continuous maps\sss into a globally\sss defined continuous map.\oss
The usual\dss tool\sss of\dss partitions of\dss unity does not\sss apply,\oss
at\sss least\sss not\sss directly,\oss
because\sss the\sss target\sss of\dss our maps\dss is\dss not\sss a vector\sss space,\oss
but\sss a\dss Grassmannian,\oss the space of\dss all\dss finitely dimensional\sss subspaces of\sss $H$\nnsp.\oss

Let\sss us start\sss with a point\sss $x\qff \in\qff U_{\fff z}$\sss
and\sss move it\sss away\sss from\sss $z$\sss along a path in\sss $X$\nnsp.\oss
After\sss we move sufficiently\dss far,\oss 
at\sss some point\sss $y\qff \in\qff X$ \sss 
we will\sss be forced\sss to replace\sss
$\varepsilon\off =\off \varepsilon_{\dff z}$\sss
by another\sss number\sss $\varepsilon'\off =\off \varepsilon_{\dff z'}$.\oss
This will\dss force us\sss to\sss jump\sss from\sss\vspace{1.5pt}\vspace{-0.87pt}
\[
\quad
V
\off =\off
\image\dff 
P_{\dff [\dff -\dff \varepsilon\fff,\qff \varepsilon\dff]}\trf
(\trf A_{\dff y}\trf)
\quad\qff
\mbox{to}\quad\qff
V\fff'
\off =\off
\image\dff 
P_{\dff [\dff -\dff \varepsilon'\fff,\qff \varepsilon'\dff]}\trf
(\trf A_{\dff y}\trf)
\pff.
\]

\vspace{-12pt}\vspace{1.5pt}\vspace{-0.87pt}
The key\sss idea\dss is\dss to connect\sss $V$\sss with\sss $V\fff'$\sss
by a segment\sss of\pss \emph{formal\trs convex\sss combinations}\vspace{0pt}
\[
\quad
t\qff
V
\off\dff +\off\dff
(\trf 1\qff -\qff t\trf)\qff
V\fff'
\pff,
\]

\vspace{-12pt}\vspace{0pt}
where\sss $0\qff \leq\qff t\qff \leq\qff 1$\nnsp.\oss
This expression\dss is\dss formal\sss in\sss the sense\sss that\sss
the multiplication and\sss the addition are purely symbolic and do not\sss
correspond\sss to any operations on subspaces of\sss $H$\nnsp.\oss 
In\sss fact\sss one needs\sss to be more precise.\oss 
If,\oss say,\pss $\varepsilon'\qff >\qff \varepsilon$\nnsp,\oss 
then\sss the orthogonal\sss complement\sss $V\fff'\dff \ominus\dff V$\sss
splits into\sss two direct\sss summands corresponding\sss the negative and\sss
positive parts of\dss the spectrum of\sss $A_{\dff y}$\nsp,\oss
and\sss one should\sss introduce a\sss formal\sss connecting segment\sss for
each\sss possible splitting.\oss
If,\pss say,\qss $y$\sss belongs\sss to\sss three
neighborhoods\sss $U_{\fff z}$,\qss $U_{\fff z'}$,\qss $U_{\fff z''}$\sss
and\sss $V$\dnsp,\qss $V\fff'$\dnsp\dnsp,\qss $V\fff''$\sss
are\sss the corresponding\sss vector\sss subspaces,\oss
then one needs\sss to use\sss formal\dss convex combinations\vspace{0pt}
\[
\quad
t_{\trf 0}\qff
V
\off\dff +\off\dff
t_{\dff 1}\qff
V\fff'
\off\dff +\off\dff
t_{\trf 2}\qff
V\fff''
\pff,
\]

\vspace{-12pt}\vspace{0pt}
where\sss 
$t_{\trf 0}\dff,\qff t_{\dff 1}\dff,\qff t_{\trf 2}\qff \geq\qff 0$\dss
and\sss
$t_{\trf 0}\qff +\qff t_{\dff 1}\qff +\qff t_{\trf 2}\off =\off 1$\nnsp,\oss
and,\oss again,\oss to\sss take into account\sss splittings.\oss
In\sss general,\oss one needs\sss to use simplices having such subspaces as vertices.\oss

\myuppar{Topological\sss categories and\sss functors.}
The proper\sss language for\sss making\sss this outline rigorous\dss
is\dss provided\sss by\sss the notions of\dss 
topological\sss categories and\sss their classifying spaces
introduced\sss by\trs Segal\qss \cite{s1}.\oss 
The reader\sss may\sss find a review of\dss definitions and\sss basic properties
in\sss the first\sss sections of\qss \cite{i2}.\oss
We will\sss now recall\sss some key constructions from\qss \cite{s1}\qss and\qss \cite{i2}.\oss

We will\sss identify each\sss topological\sss space\sss $X$\sss with\sss
the\sss topological\sss category\sss having\sss $X$\sss as\sss the space of\dss objects
and only\sss identity\sss morphisms.\oss
Every\sss covering\sss 
$U_{\dff a}\dff,\pff a\qff \in\qff \Sigma$\sss 
of\dss a space\sss $X$\sss defines a\sss topological\sss
category\sss $X_{\dff U}$\sss together\sss with a functor\sss
$\pr\dff \colon\dff
X_{\dff U}\qff \ttoo\qff X$\nnsp.\oss
This construction\dss is\dss due\sss to\dss Segal\qss \cite{s1}\qss 
and\dss is\dss recalled\sss in\dss Section\qss \ref{analytic-index-section}.\oss
If\dss the covering\dss is\dss open and\sss $X$\sss is\dss paracompact,\oss
then\sss the geometric realization\sss
$\num{\pr}\dff \colon\dff
\num{X_{\dff U}}\qff \ttoo\qff X$\sss
is\dss a homotopy equivalence.\oss
In\dss Section\qss \ref{analytic-index-section}\qss we also
recall\sss definitions of\dss a\sss topological\sss category\sss $\hat{\mathcal{S}}$\sss
having as objects finitely dimensional\sss subspaces of\sss $H$\nnsp,\oss
and of\qss Quillen--Segal\qss \cite{s4}\qss topological\sss category\sss $Q$\nnsp,\oss
a version of\qss Quillen's\dss $Q$\dnsp-construction,\oss
having as objects finitely dimensional\dss Hilbert\sss spaces.\oss
There\dss is\dss a canonical\sss homotopy equivalence\sss
$\num{\hat{\mathcal{S}}}
\qff \ttoo\qff \num{Q}$\nnsp.\oss\vspace{0.5pt}

\myuppar{Index\sss functors,\pss index\sss map,\pss and\sss the analytical\dss index.}
Let\sss us\sss return\sss to\sss the family
$A_{\dff x}\dff,\qff x\qff \in\qff X$\sss
and\sss choose a subset\sss $\Sigma\qff \subset\qff X$\sss
such\sss that\sss the neighborhoods\sss
$U_{\dff z}\dff,\pff z\qff \in\qff \Sigma$\sss
form a covering of\sss $X$\nnsp.\oss
Let\sss $\bm{\varepsilon}$\sss be\sss the family
of\dss the corresponding\sss parameters\sss $\varepsilon_{\dff z}$\nsp.\oss
The category\sss $\hat{\mathcal{S}}$\sss is\dss defined\sss in such a way\sss
that\sss the covering\sss $U$\sss together with\sss the family\sss 
$\bm{\varepsilon}$\sss define a continuous functor\vspace{1.5pt}
\[
\quad
\mathbb{A}_{\qff U,\qff \bm{\varepsilon}}\qff \colon\qff
X_{\dff U}\qff \ttoo\qff \hat{\mathcal{S}}
\pff
\]

\vspace{-12pt}\vspace{1.5pt}
which may\sss be called\dss the\qss \emph{index\dss functor},\oss
and\sss hence a continuous map\sss\vspace{1.5pt}
\[
\quad
\num{\mathbb{A}_{\qff U,\qff \bm{\varepsilon}}}\qff \colon\qff
\num{X_{\dff U}}\qff \ttoo\qff \num{\hat{\mathcal{S}}}
\pff.
\]

\vspace{-12pt}\vspace{1.5pt}
By composing\sss it\sss with a homotopy inverse\sss
$s\dff \colon\dff
X\qff \ttoo\qff \num{X_{\dff U}}$\sss
of\dss the map\sss $\num{\pr}$\sss
we get\sss a map\sss
$X\qff \ttoo\qff \num{\hat{\mathcal{S}}}$\sss
well\sss defined up\sss to homotopy.\oss
One can also combine\sss it\sss with\sss
the canonical\sss homotopy equivalence\sss
$\num{\hat{\mathcal{S}}}
\qff \ttoo\qff \num{Q}$\sss
and\sss get\sss a map\sss
$X\qff \ttoo\qff \num{Q}$\sss
well\sss defined up\sss to homotopy.\oss
The maps\sss
$X\qff \ttoo\qff \num{\hat{\mathcal{S}}}$\sss
and\sss
$X\qff \ttoo\qff \num{Q}$\sss
both deserve\sss to be called\dss the\qss \emph{index\dss maps}\pss
of\dss the family\sss $A_{\dff x}\dff,\pff x\qff \in\qff X$\nnsp.\oss
The\qss \emph{analytical\dss index}\pss of\sss
$A_{\dff x}\dff,\pff x\qff \in\qff X$\sss
is\dss defined as\sss the homotopy class of\dss either of\dss these\sss two map.\oss
Of\dss course,\oss these homotopy classes carry\sss the same information.\oss

If\dss the operators\sss $A_{\dff x}$\sss
are bounded and\sss the family\sss 
$A_{\dff x}\dff,\pff x\qff \in\qff X$\sss
is\dss continuous in\sss the norm\sss topology,\oss
then\sss the analytical\dss index carries\sss the same information\sss
as\sss the homotopy class of\dss the map\sss
$x\off \longmapsto\off A_{\dff x}$\sss
from\sss $X$\sss to\sss the space\sss $\hat{\mathcal{F}}$\sss
of\dss bounded self-adjoint\trs Fredholm\dss operators\sss
$H\qff \ttoo\qff H$\nnsp.\oss

But\sss there are definite advantages of\dss replacing\sss
$\hat{\mathcal{F}}$\sss by\sss $\num{\hat{\mathcal{S}}}$\sss
or\sss $\num{Q}$\sss as\sss the\sss target\sss of\dss the index\sss map
even\sss for\sss norm-continuous families of\dss bounded operators.\oss
Both\sss index\sss maps\sss 
$X\trf \ttoo\trf \num{\hat{\mathcal{S}}}$\sss
and\sss
$X\qff \ttoo\qff \num{Q}$\sss
realize\sss the idea of\dss index\sss as 
an equivalence class of\dss the family of\dss kernels,\oss
while\sss the map\sss $X\qff \ttoo\qff \hat{\mathcal{F}}$\sss
is\dss just\sss another name of\dss the family.\oss
The space\sss $\num{\hat{\mathcal{S}}}$\sss is\dss defined\sss
in\sss terms of\dss finitely dimensional\sss subspaces of\sss $H$\sss
and\dss is\dss more accessible\sss than\sss the space\sss $\hat{\mathcal{F}}$\sss
of\dss operators,\oss
while\sss the definition of\sss $\num{Q}$\sss involves\sss only\sss 
finitely dimensional\dss linear algebra.\oss\vspace{0.5pt}

\myuppar{Fredholm\dss families and\dss Hilbert\dss bundles.}
These aspects of\dss our definition\sss lead\sss to\sss its\sss main advantages.\oss
First,\oss the index\sss maps\sss
$X\qff \ttoo\qff \num{\hat{\mathcal{S}}}$\sss
and\sss
$X\qff \ttoo\qff \num{Q}$\sss
are defined\sss even\sss when\sss the family\sss
$A_{\dff x}\dff,\pff x\qff \in\qff X$\sss
is\dss not\sss continuous in\sss the norm\sss topology.\oss
In\dss Section\qss \ref{analytic-index-section}\qss
we will\sss define\qss \emph{Fredholm\dss families}\pss
essentially as families for\sss which\sss the construction
of\dss the index\sss functors\sss
$\mathbb{A}_{\qff U,\qff \bm{\varepsilon}}\qff \colon\qff
X_{\dff U}\qff \ttoo\qff \hat{\mathcal{S}}$\sss
makes sense.\oss
Such\sss families don't\dss need\sss to be continuous.\oss 
Second,\oss the construction of\dss these index\sss functors 
and\sss index\sss maps\sss easily\sss extends\sss to 
families of\dss operators acting\sss in\sss the fibers
of\dss a\sss Hilbert\sss bundle.\oss
This\dss is\dss the\sss topic of\dss the second\sss part\sss of\trs 
Section\qss \ref{analytic-index-section}.

Let\sss $\mathbb{H}$\sss be a\dss Hilbert\dss bundle with\sss the base $X$\nnsp,\oss
thought\sss as a family\sss of\trs Hilbert\sss spaces\sss 
$H_{\dff x}\dff,\qff x\qff \in\qff X$\nnsp.\oss
It\dss is\dss worth\sss to stress\sss that\sss in\sss general\dss there\dss is\dss
no analogue of\dss the space\sss $\hat{\mathcal{F}}$\sss suitable for\sss working\sss
with families of\dss operators acting\sss in\sss the fibers.\oss
More precisely,\oss there\dss is\dss no reasonable\sss bundle over\sss $X$\sss having\sss
as\sss the fiber over $x\qff \in\qff X$ the space of\dss self-adjoint\trs
Fredholm\dss operators\sss $H_{\dff x}\qff \ttoo\qff H_{\dff x}$\sss
with\sss the norm\sss topology.\oss
The reason\dss is\dss that\sss the\sss transition\sss maps\sss of\dss Hilbert\dss
bundles\sss are only\sss rarely\sss continuous in\sss the norm\sss topology.\oss

At\sss the same\sss time one can always define an analogue\sss 
$\hat{\mathcal{S}}\dff(\trf \mathbb{H}\trf)$\sss
of\dss the category\sss $\hat{\mathcal{S}}$\dnsp.\oss
This\dss is\dss a\sss topological\sss category\sss having as objects\sss
finitely\sss dimensional\sss subspaces of\dss fibers\sss 
$H_{\dff x}$\sss of\dss $\mathbb{H}$\nnsp.\oss
There\dss is\dss also analogue\sss $Q\dff(\trf \mathbb{H}\trf)$\sss of\sss $Q$\nnsp.\oss
It\dss is\dss especially simple and\sss has as objects\sss pairs\sss
$(\trf x\fff,\qff V\trf)$\nnsp,\oss
where\sss $x\qff \in\qff X$\sss and\sss $V$\sss is\dss an object\sss of\sss $Q$\nnsp.\oss
Both definitions are recalled\sss in\dss Section\qss \ref{analytic-index-section}.\oss
By\sss the definition\sss the categories\sss $Q\dff(\trf \mathbb{H}\trf)$\sss depend only on\sss $X$\nnsp.\oss
Also,\qss 
$\num{Q\dff(\trf \mathbb{H}\trf)}
\off =\off
X\dff \times\dff \num{Q}$\nnsp.\oss
The categories $\hat{\mathcal{S}}\dff(\trf \mathbb{H}\trf)$
depend only on\sss $X$\sss up\sss to an\sss isomorphism.\oss
In\sss fact,\oss by a well\sss known\sss theorem of\qss Dixmier\dss and\dss Douady\qss \cite{dd}\qss
every\dss Hilbert\dss bundle over a paracompact\sss space\dss
is\dss isomorphic\sss to\sss the\sss trivial\sss bundle\sss
$X\dff \times\dff H\qff \ttoo\qff X$\nnsp.\oss
Moreover,\oss a\sss trivialization of\sss $\mathbb{H}$\sss allows\sss to identify\sss
$\num{\hat{\mathcal{S}}\dff(\trf \mathbb{H}\trf)}$\sss
with\sss
$X\dff \times\dff \num{\hat{\mathcal{S}}}$\nnsp,\oss
and such an\sss identification\dss is\dss unique up\sss to homotopy.\oss

With\sss the obvious changes\sss the definition of\pss \emph{Fredholm\dss families}\pss
applies\sss to families of\dss self-adjoint\sss operators\sss
$A_{\dff x}\dff \colon\dff H_{\dff x}\qff \ttoo\qff H_{\dff x}$\nsp,\qss
$x\qff \in\qff X$\nnsp.\oss
For such\sss families one can define\qss 
\emph{index\dss functors}\vspace{0.5pt}
\[
\quad
\mathbb{A}_{\qff U,\qff \bm{\varepsilon}}\qff \colon\qff
X_{\dff U}\qff \ttoo\qff \hat{\mathcal{S}}\dff(\trf \mathbb{H}\trf)\dff,\quad
X_{\dff U}\qff \ttoo\qff Q\dff(\trf \mathbb{H}\trf)
\pff,
\]

\vspace{-12pt}\vspace{0.25pt}
and\qss \emph{index\dss maps}\dss\vspace{0.25pt}
\[
\quad
X\qff \ttoo\qff \num{\hat{\mathcal{S}}\dff(\trf \mathbb{H}\trf)}
\off =\off
X\dff \times\dff \num{\hat{\mathcal{S}}}\dff,\qquad
X\qff \ttoo\qff \num{Q\dff(\trf \mathbb{H}\trf)}
\off =\off
X\dff \times\dff \num{Q}
\pff.
\]

\vspace{-12pt}\vspace{0.5pt}
The composition\sss of\dss the index\sss maps\sss
with\sss the projections\sss to\sss $X$\sss
are equal\sss to\sss the identity\sss map and 
carry\sss no information.\oss
This suggests\sss to define another\qss \emph{index\dss maps}\qss 
$X\qff \ttoo\qff \num{\hat{\mathcal{S}}}$\sss
and\dss
$X\qff \ttoo\qff \num{Q}$\dss
as\sss the compositions of\dss the original\sss ones
with\sss the projections\sss to\sss the second\sss factors.\oss
The\qss \emph{analytical\dss index}\pss is\dss defined as\sss the homotopy 
class of\dss either of\dss these\sss index\sss maps.\oss
Both\sss versions carry\sss the same\sss information,\oss
as also\sss the homotopy class of\sss
$X\qff \ttoo\qff \num{\hat{\mathcal{S}}\dff(\trf \mathbb{H}\trf)}$\nnsp.\oss\vspace{-0.125pt}

\myuppar{Strictly\trs Fredholm\dss families.}
In order\sss to go further\sss than\sss the definitions and\sss their correctness,\oss
we need\sss to\sss restrict\sss a\sss little\sss the class of\dss considered\sss families
and\sss require\sss that\sss spectral\sss projections corresponding\sss to some\sss
half-lines continuously depend on parameter.\oss

Suppose\sss first\sss that\sss $\mathbb{H}$\sss is\dss the\sss trivial\sss bundle
with\sss the fiber\sss $H$\nnsp.\oss
We will\sss say\sss that\sss the family\sss
$A_{\dff x}\dff \colon\dff H\qff \ttoo\qff H$\nsp,\qss
$x\qff \in\qff X$\sss
is\qss \emph{strictly\trs Fredholm}\pss
if\dss it\dss is\trs Fredholm\dss and\sss
for every\sss $z\qff \in\qff X$\sss there exist\sss $\varepsilon\qff >\qff 0$\sss
and\dss a neighborhood\sss $U_{\fff z}$\sss of\sss $z$\sss
such\sss that\sss $(\trf A_{\dff y}\fff,\qff \varepsilon\trf)$\sss
is\dss an enhanced operator for every\sss 
$y\qff \in\qff U_{\fff z}$ and\sss
the spectral\sss projections\sss 
$P_{\dff [\trf \varepsilon\fff,\qff \infty\dff)}\trf
(\trf A_{\dff y}\trf)$\sss
continuously depends on\sss $y$\sss
for\sss 
$y\qff \in\qff U_{\fff z}$\nsp.\oss

In\sss general,\pss 
$A_{\dff x}\dff \colon\dff H_{\dff x}\qff \ttoo\qff H_{\dff x}$\nsp,\qss
$x\qff \in\qff X$\sss
is\dss said\sss to be a\qss \emph{strictly\trs Fredholm\dss family}\pss
if\dss for every\sss $z\qff \in\qff X$\sss there exists\sss 
a neighborhood\sss $U_{\fff z}$\sss of\sss $z$\sss
and a\sss trivialization of\dss $\mathbb{H}$\sss
over\sss $U_{\fff z}$\sss turning\sss the restriction\sss
$A_{\dff x}\dff \colon\dff H_{\dff x}\qff \ttoo\qff H_{\dff x}$\nsp,\qss
$x\qff \in\qff U_{\fff z}$\sss
into a strictly\trs Fredholm\dss family\sss
in\sss the above sense.

\myuppar{Families of\dss elliptic operators.}
A continuous\sss family\sss of\dss elliptic pseudo-differential\sss operators
of\dss order\sss $0$\sss
acting on\sss the fibers of\dss locally\sss trivial\sss bundle\sss 
$\mathbb{M}$\sss over\sss $X$\sss with closed\sss manifolds as fibers 
defines a family of\dss bounded operators\sss acting in\sss 
the fibers of\dss a\dss Hilbert\dss bundle\sss $\mathbb{H}$\nnsp.\oss
Classical\sss results of\trs Seeley\qss \cite{see}\qss
(namely,\oss the integral\sss presentation\sss in\sss the proof\dss of\trs
Theorem\qss 5\qss in\qss \cite{see})\qss 
imply\sss that\sss the\sss latter\sss family\dss is\dss strictly\trs Fredholm.\oss

\myuppar{Strictly\sss adapted\dss trivializations.}
A\sss trivialization of\dss $\mathbb{H}$\sss
is\dss said\sss to be\pss \emph{strictly\sss adapt\-ed}\pss to the family\sss
$A_{\dff x}\dff \colon\dff H_{\dff x}\qff \ttoo\qff H_{\dff x}$\nsp,\qss
$x\qff \in\qff X$\trs
if\dss its restrictions\sss to sufficiently small\sss neighborhoods\sss
$U_{\fff z}$\sss turn\sss the restrictions\sss
$A_{\dff x}\dff \colon\dff H_{\dff x}\qff \ttoo\qff H_{\dff x}$\nsp,\qss
$x\qff \in\qff U_{\fff z}$\sss
into strictly\trs Fredholm\dss families,\oss
i.e.\qss turn\sss the families\sss
$P_{\dff [\trf \varepsilon\fff,\qff \infty\dff)}\trf
(\trf A_{\dff y}\trf)$\nsp,\sss
$y\qff \in\qff U_{\fff z}$\sss
with\sss appropriate\sss $\varepsilon$\sss
into norm-continuous ones.\oss
Strictly adapted\sss trivializations
exist\sss for every\sss strictly\trs Fredholm\dss family\sss
with\sss paracompact\sss $X$\nnsp.\oss
This\dss is\dss the main\sss result\sss of\trs Section\qss \ref{strictly-fredholm-families}.\oss
It\dss is\dss proved\sss first\sss for\sss triangulable spaces\sss $X$\nnsp,\oss
when one can\sss use\sss an induction\sss by skeletons,\oss
and\sss then\sss for paracompact\sss $X$\nnsp.\oss
See\dss Theorems\qss \ref{adapted}\qss and\qss \ref{adapted-paracompact}.\oss
The proofs depend\sss on\sss contractibility\sss theorems of\trs Kuiper\qss \cite{ku}\qss
and of\qss Dixmier--Douady\qss \cite{dd}\qss
in\dss Atiyah--Segal\qss \cite{ase}\qss form.\oss
For\sss the rest\sss of\dss the introduction\sss 
we assume\sss that\sss $X$\sss is\dss paracompact.

\myuppar{Polarizations.}
Recall\sss that\sss a\qss \emph{polarization}\pss of\dss a\sss Hilbert\sss space\sss $H$\sss
is\dss a presentation of\sss $H$\sss as an orthogonal\sss direct\sss sum\sss
$H\off =\off H_{\dff -}\dff \oplus\dff H_{\dff +}$\sss
of\sss two closed\sss infinitely dimensional\sss subspaces.\oss
A polarization\sss leads\sss to\sss the\qss
\emph{restricted\dss Grassmannian}\qss $\gr$\nnsp,\oss
the space of\dss subspaces\sss $K\qff \subset\qff H$\qss
\emph{commensurable}\pss with $H_{\dff -}$\nsp,\oss i.e.\qss
such\sss that\dss 
$K\dff \cap\dff H_{\dff -}$\sss
has finite codimension\sss in\sss both\sss $K$\sss and\dss
$H_{\dff -}$\nsp.\oss
The\sss topology of\sss $\gr$\sss
is\dss defined\sss by\sss the norm\sss topology of\dss orthogonal\dss projections\sss
$H\qff \ttoo\qff K$\nnsp.\oss
A strictly\trs Fredholm\dss family\sss leads\sss to 
a\qss \emph{polarization}\pss of\dss
$\mathbb{H}$\nnsp,\pss a structure defined\sss in\sss the\dss Appendix.\oss
In\sss the main\sss part\sss of\dss the paper\sss
these structures are not\sss even mentioned,\oss 
but\sss in\sss the\dss Appendix\dss Theorems\qss \ref{adapted}\qss and\qss \ref{adapted-paracompact}\qss 
and\dss their proofs are\sss translated\sss into\sss the\sss language of\dss polarizations.\oss

In\dss Section\qss \ref{analytic-index-strict}\qss
we recall\sss from\qss \cite{i2}\qss the definition of\dss a\sss topological\sss category\sss
$\mathcal{P}{\nsp}\hat{\mathcal{S}}$\dnsp,\oss
a version of\dss the category\sss $\hat{\mathcal{S}}$\sss
involving\sss polarization.\oss
There\dss is\dss a\sss forgetting\sss functor\sss
$\mathcal{P}{\nsp}\hat{\mathcal{S}}
\qff \ttoo\qff 
\hat{\mathcal{S}}$\sss
inducing a homotopy equivalence\sss
$\num{\mathcal{P}{\nsp}\hat{\mathcal{S}}}
\qff \ttoo\qff 
\num{\hat{\mathcal{S}}}$\nnsp.\oss
If\dss we use a strictly adapted\sss trivialization\sss to identify\sss
$\num{\hat{\mathcal{S}}\dff(\trf \mathbb{H}\trf)}$\sss
with\sss
$X\dff \times\dff \num{\hat{\mathcal{S}}}$\nnsp,\oss
then\sss the\sss index\sss map\sss
$X\qff \ttoo\qff \num{\hat{\mathcal{S}}}$\sss
naturally\sss lifts\sss to a map\sss 
$X\qff \ttoo\qff \num{\mathcal{P}{\nsp}\hat{\mathcal{S}}}$\nnsp,\oss
which we call\sss the\qss \emph{polarized\dss index\dss map}.\oss
See\dss Section\qss \ref{analytic-index-strict}.\oss

\myuppar{Grassmannian\dss bundles.}
Let\sss 
$A_{\dff x}\dff \colon\dff H_{\dff x}\qff \ttoo\qff H_{\dff x}$\nsp,\qss
$x\qff \in\qff X$\sss
be a strictly\trs Fredholm\dss family.\oss
We will\sss denote it\sss also by\sss $\mathbb{A}$\nnsp.\oss
If\sss $z\qff \in\qff X$\sss and\sss
$(\trf A_{\dff z}\fff,\qff \varepsilon\trf)$\sss 
is\dss an enhanced operator,\oss then\vspace{0pt}
\[
\quad
H_{\dff z}
\off =\off
\image\dff
P_{\dff (\dff -\qff \infty\fff,\qff \varepsilon\trf]}\trf
(\trf A_{\dff z}\trf)
\qff \oplus\qss
\image\dff
P_{\dff [\trf \varepsilon\fff,\qff \infty\dff)}\trf
(\trf A_{\dff z}\trf)
\]

\vspace{-12pt}\vspace{0pt}
is\dss a polarization of\sss $H_{\dff z}$\nsp.\oss
This polarization\sss leads\sss to\sss a\sss
restricted\dss Grassmannian,\oss
which\sss will\sss denote by\sss $\gr\trf(\trf z\trf)$\nnsp.\oss
Since\sss the family\sss $\mathbb{A}$\dss
is\dss strictly\trs Fredholm,\oss
the family\sss of\dss restricted\dss Grassmannians\sss 
$\gr\trf(\trf z\trf)\dff,\qff
z\qff \in\qff X$\sss
forms a\sss locally\sss trivial\sss bundle\vspace{0pt}
\[
\quad
\bm{\pi}\trf(\trf \mathbb{A}\trf)\dff \colon\dff
\gr\trf(\trf \mathbb{A}\trf)
\qff \ttoo\qff 
X
\pff.
\]

\vspace{-12pt}\vspace{0pt}
In\trs Theorem\qss \ref{induced-bundle}\qss we prove\sss
that\sss this bundle\dss is\dss equal\sss
to\sss the bundle induced\sss by\sss the polarized\sss index\sss map\sss
$X\qff \ttoo\qff \num{\mathcal{P}{\nsp}\hat{\mathcal{S}}}$
from\sss the standard\sss bundle\sss
$\bm{\pi}\dff \colon\dff
\mathbf{G}
\qff \ttoo\qff 
\num{\mathcal{P}{\nsp}\hat{\mathcal{S}}}$\dss
constructed\sss in\qss \cite{i2}.\oss
Motivated\sss by\trs
Melrose\dss and\dss Piazza\qss \cite{mp},\oss
we introduce another name,\oss \emph{weak\sss spectral\dss sections}\pss
of\dss the family\sss $\mathbb{A}$\nnsp,\oss
for continuous sections of\dss the bundle\sss $\bm{\pi}\trf(\trf \mathbb{A}\trf)$\nnsp.\oss
In an agreement\sss with\sss the ideas of\qss Melrose\dss and\dss Piazza\qss \cite{mp},\oss
a weak\sss spectral\sss section exists\dss if\trs and\dss only\trs if\trs
the analytical\dss index\sss vanishes,\oss
i.e.\qss the index\sss map\dss is\dss homotopic\sss to a constant\sss map.\oss
See\dss Theorem\qss \ref{weak-spectral-sections}.\oss
The proof\dss of\dss the\qss ``only\trs if''\qss part\sss depends on\sss
one of\dss the main\sss results of\pss \cite{i2},\oss
the contractibility of\dss the\sss total\sss space\sss $\mathbf{G}$\sss
of\dss the bundle\sss $\bm{\pi}$\nnsp.\oss
See\qss \cite{i2},\oss Theorem\qss E\qss or\trs Theorem\qss 13.6.\oss\vspace{-0.8pt}

\myuppar{Discrete-spectrum\dss families.}
Let\sss us\sss say\sss that\sss
$A_{\dff x}\dff \colon\dff H_{\dff x}\qff \ttoo\qff H_{\dff x}$\nsp,\qss
$x\qff \in\qff X$\sss
is\dss a\qss \emph{discrete-spectrum\dss family}\pss
if\dss for every\sss $\lambda\qff \in\qff \rrr$\sss the family\sss of\dss operators\sss 
$A_{\dff x}\qff -\pff \lambda$\nsp,\qss $x\qff \in\qff X$\sss
is\dss a\dss Fredholm\dss family.\oss
As usual,\oss here\sss $\lambda$\sss stands for\sss $\lambda$\sss times\sss the identity operator.\oss
The operators $A_{\dff x}$ of\dss such a family\sss have discrete spectrum,\oss
i.e.\qss their essential\sss spectrum\dss is\dss empty.\oss
By\sss this reason\sss these operators cannot\sss be bounded.\oss
We assume\sss that\sss they are closed densely defined operators.\oss\vspace{-0.8pt}

Following\trs Melrose\dss and\dss Piazza\qss \cite{mp},\oss
we say\sss that\sss a weak\sss spectral\sss section\sss
$\mathbb{S}\dff \colon\dff
X\qff \ttoo\qff \gr\trf(\trf \mathbb{A}\trf)$\sss
is\dss a\qss \emph{spectral\sss section}\pss
if\trs there exists a continuous function\dss
$r\dff \colon\dff X\qff \ttoo\qff \rrr_{\trf >\dff 0}$\sss 
such\sss that\vspace{1.5pt}
\[
\quad
\image
P_{\qff [\trf r\dff(\dff x\trf)\fff,\pff \infty\qff)}\qff(\dff A_{\dff x}\dff)
\off \subset\off\dff
\mathbb{S}\dff(\trf x\trf)
\off \subset\off\dff
\image
P_{\qff [\qff -\qff r\dff(\dff x\trf)\fff,\pff \infty\qff)}\qff(\dff A_{\dff x}\dff)
\pff
\]

\vspace{-12pt}\vspace{1.5pt}
for every\sss $x\qff \in\qff X$\nnsp.\oss
This\sss notion seems\sss to be useful\sss only\sss for
discrete-spectrum strictly\trs Fredholm\dss families.\oss
In\sss this case we prove\sss that\sss every weak\sss spectral\sss section\dss
can\sss be deformed\sss into a spectral\sss sections,\oss
and\sss hence spectral\sss sections exist\dss if\trs and\dss only\trs if\dss
the analytical\dss index\sss vanishes.\oss
See\dss Theorem\qss \ref{two-types-sections}\qss and\dss
Corollary\qss \ref{index-and-mp-sections}.\oss
These results generalize and\sss clarify\dss Proposition\qss 1\qss of\qss
Melrose\dss and\dss Piazza\trs \cite{mp}.\oss
See\sss the discussion at\sss the end of\trs Section\qss \ref{discrete-spectrum}.\oss\vspace{-0.8pt}

\myuppar{Families of\dss general\trs Fredholm\dss operators.}
For\sss families of\dss general\qss (not\sss assumed\sss to be self-adjoint\fff)\qss
Fredholm\dss operators one can develop a\sss theory\sss parallel\sss 
to\sss the one outlined above.\oss
See\dss Section\qss \ref{analytic-index-fredholm-non-sa}.\oss
In\sss the present\sss paper we\sss limit\sss the discussion\sss
by\sss the results needed\sss to prove\sss the equivalence of\dss
our definition of\dss the analytical\dss index\sss with\dss
Atiyah--Singer\dss one when\sss the\sss latter\sss or\dss a natural\sss
generalization of\trs it\sss applies.\oss\vspace{-0.8pt}\vspace{-0.125pt}

\myuppar{Comparing\dss two definitions of\dss analytical\dss index.}
Since\trs Atiyah--Singer\dss definition of\dss
the analytical\dss index of\dss families of\dss self-adjoint\sss operators
depends on\sss the definition\sss in\sss the non-self-adjoint\sss case,\oss
we have\sss to deal\sss with\sss this case first.\oss
The most\sss classical\sss definition\sss applies only\sss to compact\sss $X$\nnsp.\oss
It\sss was extended\sss to paracompact $X$ by\dss Segal\qss \cite{sfc}.\pss
We prove\sss that\sss two definitions agree\sss for compact\sss $X$
in\trs Theorem\qss \ref{fredholm-agree},\pss and\sss for paracompact\sss $X$
in\trs Theorem\qss \ref{fredholm-agree-paracompact}.\vspace{-0.8pt}

In\dss Section\qss \ref{classical-analytic-index-sa}\qss we deal\sss with\sss
families of\dss self-adjoint\sss operators.\oss
First,\oss we show\sss that\sss one can extend\sss the approach of\qss 
Atiyah\dss and\dss Singer\dss to\dss Fredholm\dss families as defined above.\oss
See\dss Theorem\qss \ref{fr-to-fr}.\oss
Then we prove\sss that\sss two definitions agree for\sss strictly\trs Fredholm\dss families.\oss
See\dss Theorem\qss \ref{sa-index-agree}.\oss
Along\sss the way\sss we prove\sss that\sss every\sss strictly\trs Fredholm\dss family\sss
can\sss be replaced\sss by a much better\sss family\sss without\sss affecting\sss its index.\oss 
See\dss Theorem\qss \ref{index-replacement}.\oss
This construction\dss is\dss similar\sss to a spectral\sss deformation\sss used\sss in\qss \cite{as}\qss
and\sss may\sss be of\dss independent\sss interest.\oss

\newpage
\mysection{Analytical\qss index\qss of\qss self-adjoint\qss Fredholm\qss families}{analytic-index-section}

\myuppar{Self-adjoint\dss Fredholm\sss operators.}
Let\sss us\sss fix a separable infinitely dimensional\dss Hilbert\sss space $H$\nnsp.\oss
By\sss $\hat{\mathcal{F}}$\sss we denote\sss the space of\dss self-adjoint\dss
Fredholm\dss operators\sss $H\qff \ttoo\qff H$\nnsp.\oss
Actually,\oss there are\sss at\sss least\sss two versions of\dss the space\sss $\hat{\mathcal{F}}$\dnsp.\oss
In\sss the first\sss version\sss the operators are assumed\sss to be bounded and\sss
$\hat{\mathcal{F}}$\sss is\dss equipped\sss with\sss the norm\sss topology.\oss
In\sss the second\sss version\sss the operators are allowed\sss to be unbounded
closed densely defined operators.\oss
In\sss this version $\hat{\mathcal{F}}$\sss is\dss equipped\sss with\dss Riesz\dss topology.\oss
Our discussion applies equally\sss well\sss to both versions.\oss
Most\sss of\dss our arguments work also for\sss
the\sss topology of\dss convergence in\sss the uniform resolvent\sss sense.\oss

\myuppar{Families of\dss self-adjoint\dss Fredholm\sss operators.}
Let\sss $X$\sss be a\sss compactly\sss generated\sss topological\sss space.\oss 
A family\sss of\dss self-adjoint\dss Fredholm\sss operators parameterized\sss by $X$\sss
is\dss a continuous map\sss
$\mathbb{A}\dff \colon\dff X\qff \ttoo\qff \hat{\mathcal{F}}$\nsp\dnsp.\oss
It\dss is\dss convenient\sss to\sss denote\sss $\mathbb{A}\trf(\trf x\trf)$\sss
by\sss $A_{\dff x}$ and\sss use for\sss $\mathbb{A}$\sss the notation\sss
$A_{\dff x}\fff,\qff x\qff \in\qff X$\nnsp,\oss
which\sss reflects\sss the idea of\dss 
a\qss ``family''\qss of\dss operators better.\oss

\myuppar{The analytical\dss index of\dss families.}
The\qss \emph{analytical\dss index}\pss of\dss a family\sss
$A_{\dff x}\fff,\qff x\qff \in\qff X$\dss should\sss be an element\sss of\dss
the group\sss $K^{\dff 1}\dff(\trf X\trf)$\nnsp.\oss
The group\sss $K^{\dff 1}\dff(\trf X\trf)$\sss can\sss be defined\sss
as\sss the group of\dss homotopy classes of\dss maps from $X$\sss to a classifying space
for\sss the functor\sss
$Y\off \longmapsto\off K^{\dff 1}\dff(\trf Y\trf)$\nnsp,\oss
initially defined only for compact\sss $Y$\dnsp.\oss 
The group structure results from a canonical\sss structure of\dss an $H$\dnsp-space on\sss
the classifying space.\oss
In\sss the case of\dss bounded operators\sss the space $\hat{\mathcal{F}}$\sss
is\dss such a classifying space by\sss the results of\trs 
Atiyah\dss and\trs Singer\qss \cite{as}.\oss
This allows\sss to define\sss the analytical\dss index of\dss a family\sss
$A_{\dff x}\fff,\qff x\qff \in\qff X$\dss simply as\sss the homotopy class
of\dss the corresponding\sss map\sss
$\mathbb{A}\dff \colon\dff X\qff \ttoo\qff \hat{\mathcal{F}}$\nsp\dnsp,\oss
and\sss this\dss is\dss the most\sss widely,\oss 
if\dss not\sss exclusively,\oss used definition,\oss
even\sss for compact\sss $X$\nnsp.\oss

This definition seems\sss to be somewhat\sss tautological,\oss
in particular\dss if\dss compared\sss with\sss the standard definition
of\dss the analytical\dss index of\dss a\sss family\sss of\qss Fredholm\sss operators.\oss
At\sss the same\sss time\sss this definition\dss is\dss not\sss quite suitable
for generalizations\sss to families of\dss unbounded operators and\sss to families
of\dss operators in\dss Hilbert\sss spaces varying with\sss $x\qff \in\qss X$\nnsp,\oss
i.e.\qss to operators in\dss Hilbert\dss bundles over $X$\nnsp.\oss
Such\sss generalizations\sss encounter some subtle continuity\sss problems.\oss

\myuppar{Topological\sss categories\sss related\dss to
self-adjoint\dss Fredholm\sss operators.}
Let\sss us\sss recall\sss some definitions from\qss \cite{i2}.\oss
To begin with,\oss the space $\hat{\mathcal{F}}$ can\sss be considered 
as a\sss topological\sss category\sss
having $\hat{\mathcal{F}}$ 
as\sss the space of\dss objects\qss and only\sss identity\sss morphisms.\oss
Then\sss the classifying space $\num{\hat{\mathcal{F}}}$\sss
is\dss equal\sss to $\hat{\mathcal{F}}$ considered as a\sss topological\sss space.\oss

An\qss \emph{enhanced\sss self-adjoint\trs Fredholm\dss operator},\oss
or\sss simply\sss an\qss \emph{enhanced operator}\pss is\dss defined as
a pair\sss $(\trf A\dff,\qff \varepsilon\trf)$\sss consisting of\dss an operator\sss
$A\qff \in\qff \hat{\mathcal{F}}$ and a real\sss number\sss $\varepsilon\qff >\qff 0$\sss
such\sss that\sss $-\qff \varepsilon\fff,\qff \varepsilon$\sss do not\sss belong\sss to\sss
the spectrum\sss $\sigma\dff(\dff A\dff)$ of\sss $A$\sss and\sss the spectral\sss projection\sss
$P_{\dff [\dff -\dff \varepsilon\fff,\qff \varepsilon\dff]}\trf(\trf A\trf)$\sss
has\sss finite rank.\oss
The set\sss $\hat{\mathcal{E}}$ of\dss enhanced operators\dss is\dss
equipped\sss with\sss the\sss topology defined\sss by\sss
the\sss topology of\sss $\hat{\mathcal{F}}$\sss
and\sss the discrete\sss topology on\sss the space\sss $\rrr_{\qff >\dff 0}$ of\dss
parameters $\varepsilon$\nnsp.\oss
The space\sss $\hat{\mathcal{E}}$\sss 
is\dss ordered\sss by\sss the relation\sss $\leq$\nnsp,\oss where
$(\trf A\dff,\qff \varepsilon\trf)
\off \leq\off
(\trf A'\dff,\qff \varepsilon'\trf)$\trs
if\qss
$A\off =\off A'$\sss
and\dss
$\varepsilon\qff \leq\qff \varepsilon'$\nnsp.\oss
This order defines a structure of\dss a\sss topological\sss
category on $\hat{\mathcal{E}}$\sss having\sss $\hat{\mathcal{E}}$\sss as\sss the space of\dss objects
and a single morphism\sss
$(\trf A\dff,\qff \varepsilon\trf)
\qff \ttoo\qff
(\trf A'\dff,\qff \varepsilon'\trf)$\dss
if\trs $A\off =\off A'$\sss and\sss
$\varepsilon\qff \leq\qff \varepsilon'$\nnsp.\oss
The obvious forgetting\sss functor\sss
$\hat{\varphi}\dff \colon\dff \hat{\mathcal{E}}\qff \ttoo\qff \hat{\mathcal{F}}$\sss
induces a map of\dss classifying spaces
$\num{\hat{\varphi}}\dff \colon\dff 
\num{\hat{\mathcal{E}}}
\qff \ttoo\qff 
\num{\hat{\mathcal{F}}}
\off =\off
\hat{\mathcal{F}}$\dnsp.\oss

The\sss topological\sss category\sss $\hat{\mathcal{S}}$\sss 
has\sss finitely\sss dimensional\sss subspaces
of\sss $H$\sss as objects.\oss
The\sss topology on\sss the space of\dss such subspaces\dss is\dss the obvious one.\oss
Morphisms\sss $V\qff \ttoo\qff V\fff'$\sss exist\sss only\dss if\trs
$V\qff \subset\qff V\fff'$\nnsp,\oss and\sss in\sss this case morphisms correspond\sss to
orthogonal\sss decompositions\dss\vspace{1.5pt}
\[
\quad
V\fff'
\off =\off
U_{\dff -}\qff \oplus\qff
V\qff \oplus\qff U_{\dff +}
\pff.
\]

\vspace{-12pt}\vspace{1.5pt}
The subspaces\sss $U_{\dff -}$ and\sss $U_{\dff +}$\sss are called\sss the\qss
\emph{negative}\qss and\qss \emph{positive parts}\oss of\dss the morphism\sss in question.\oss
The composition\dss is\dss defined\sss by\sss taking\sss the sum of\dss
negative parts and\sss the sum of\dss positive parts\sss to get,\oss
respectively,\oss the negative and\sss
the positive parts of\dss the composition.\oss 
The\sss topology on\sss the set\sss of\dss
morphisms\dss is\dss the obvious one.\oss
There\dss is\dss a\sss forgetting\sss functor\sss
$\hat{\pi}\dff \colon\dff \hat{\mathcal{E}}\qff \ttoo\qff \hat{\mathcal{S}}$\dnsp,\oss
defined as follows.\oss
The functor\sss $\hat{\pi}$\sss takes
an enhanced operator\sss $(\trf A\dff,\qff \varepsilon\trf)$\sss to\vspace{1.5pt}
\[
\quad
\hat{\pi}\trf(\trf A\dff,\qff \varepsilon\trf)
\off =\off\dff
\image P_{\dff [\dff -\dff \varepsilon\fff,\qff \varepsilon\dff]}\trf(\trf A\trf)
\off \subset\off\dff
H
\pff,
\]

\vspace{-12pt}\vspace{1.5pt}
and\sss a morphism\sss
$(\trf A\dff,\qff \varepsilon\trf)
\qff \ttoo\qff
(\trf A\dff,\qff \varepsilon'\trf)$\sss
to\sss the morphism corresponding\sss to\sss the decomposition\vspace{3pt}
\begin{equation*}
\quad
\image P_{\dff [\dff -\dff \varepsilon'\fff,\qff \varepsilon'\dff]}\trf(\trf A\trf)
\off\qff =\off\qff
\left(\qff
\image P_{\dff [\dff -\dff \varepsilon'\fff,\qff -\dff \varepsilon\dff]}\trf(\trf A\trf)
\qff\right)
\qff \oplus\qff
\left(\qff
\image P_{\dff [\dff -\dff \varepsilon\fff,\qff \varepsilon\dff]}\trf(\trf A\trf)
\qff\right)
\qff \oplus\qff
\left(\qff
\image P_{\dff [\dff \varepsilon\fff,\qff -\dff \varepsilon'\dff]}\trf(\trf A\trf)
\qff\right)
\off.
\end{equation*}

\vspace{-12pt}\vspace{3pt}
Finally,\oss the category $Q$\sss is\dss an abstract\sss analogue 
of\sss $\hat{\mathcal{S}}$\dnsp.\oss
It\dss objects are finitely\sss dimensional\dss Hilbert\sss spaces\sss
belonging\sss to a set\sss containing at\sss least\sss one\sss Hilbert\sss space
of\dss each\sss finite dimension.\oss
The set\sss of\dss objects\dss is\dss equipped\sss with\sss the discrete\sss topology.\oss
A morphism\sss $V\qff \ttoo\qff V\fff'$\sss 
is\dss a\sss triple 
$(\trf U_{\dff -}\dff,\qff U_{\dff +}\dff,\qss \iota\qff)$\nnsp,\oss
where\sss
$\iota\dff \colon\dff
V\qff \ttoo\qff V\fff'$\sss
is\dss an\sss isometric embedding\sss and\dss
$U_{\dff -}\dff,\pff U_{\dff +}$\dss
are subspaces of\dss $V\fff'$
defining an orthogonal\sss decomposition\vspace{1.5pt}
\[
\quad
V\fff'
\off =\off
U_{\dff -}\qff \oplus\qff
\iota\dff(\trf V\trf)\qff \oplus\qff U_{\dff +}
\qff.
\]

\vspace{-12pt}\vspace{1.5pt}
The composition of\dss morphisms\dss is\dss defined\sss in\sss the same
way\sss as in $\hat{\mathcal{S}}$\dnsp,\oss
and\sss the\sss topology on\sss the set\sss of\dss
morphisms\dss is\dss the obvious one.\oss
The obvious functor\sss
$\hat{\mathcal{S}}\qff \ttoo\qff Q$\sss
assigning\sss to a finitely dimensional\sss subspace\sss $V\qff \subset\pff H$\sss
the space $V$ considered as an object\sss of\sss $Q$\nnsp,\oss
is\dss not\sss continuous.\oss
But\sss there\dss is\dss an\sss intermediate category\sss $Q/\fff H$\sss
and continuous functors\sss
$\hat{\mathcal{S}}\off \longleftarrow\off Q/\fff H\qff \ttoo\qff Q$\sss
allowing\sss to compare\sss $\hat{\mathcal{S}}$ and\sss $Q$\nnsp.\oss
See\qss \cite{i2},\oss Section\qss 9.\oss
By\sss passing\sss to\sss the classifying spaces 
we get\sss the following\sss two diagrams.\oss\vspace{-1.5pt}
\[
\quad
\begin{tikzcd}[column sep=boom, row sep=large]
\hat{\mathcal{F}}
&
\protect{\num{\hat{\mathcal{F}}}}
\arrow[l, "\dis ="']
&
\protect{\num{\hat{\mathcal{E}}}}
\arrow[r, "\dis \protect{\num{\hat{\pi}}}"]
\arrow[l, "\dis \protect{\num{\hat{\varphi}}}"']
&
\protect{\num{\hat{\mathcal{S}}}}
\end{tikzcd}
\]

\vspace{-34.5pt}
\[
\quad
\begin{tikzcd}[column sep=boom, row sep=large]
\protect{\num{\hat{\mathcal{S}}}}
&
\protect{\num{Q/\fff H}}
\arrow[l]
\arrow[r]
&
\protect{\num{Q}}
\end{tikzcd}
\]

\vspace{-12pt}
The basic result\sss about\sss these categories and\sss functors\dss is\dss
the fact\sss that\sss all\dss maps in\sss these diagrams are homotopy\sss equivalences.\oss
See\qss \cite{i2},\oss Theorems\qss 9.6\qss and\qss 9.7.\oss

\myuppar{Topological\sss categories defined\sss by coverings.}
We need\sss a classical\sss construction due\sss to\qss Segal\qss \cite{s1}.\oss  
It\sss relates\sss topological\sss categories with
coverings\sss of\dss  topological\sss spaces.\oss
The construction\sss starts with a covering\sss
$U_{\dff a}\dff,\pff a\qff \in\qff \Sigma$\sss
of\dss a space\sss $X$\nnsp.\oss
The idea\dss is,\pss to borrow a phrase of\qss Segal\qss \cite{s3},\oss
to\qss ``disintegrate''\qss the space $X$\sss into pieces 
$U_{\dff a}\dff,\pff a\qff \in\qff \Sigma$\sss and\sss then
assemble\sss these piece into a\qss ``thick''\qss version of\sss $X$\nnsp.\oss
Let\sss $\Sigma^{\dff \fin}$\sss be\sss the set\sss of\dss finite non-empty\sss
subsets\sss of\sss $\Sigma$\sss and\vspace{3pt}
\[
\quad
U_{\dff \sigma}
\off =\off
\bigcap\nolimits_{\pff a\qff \in\qff \sigma}\dff U_{\dff a}
\qff
\]

\vspace{-12pt}\vspace{3pt}
for every\sss $\sigma\qff \in\qff \Sigma^{\dff \fin}$\dnsp.\oss
Following\qss Segal\qss \cite{s1},\oss
let\sss us consider\sss the following\sss topological\sss category\sss $X_{\dff U}$\nsp.\oss
Its space of\dss objects\dss is\dss the disjoint\sss union of\dss the subspaces\sss
$U_{\dff \sigma}\dff,\pff \sigma\qff \in\qff \Sigma^{\dff \fin}$\dnsp.\oss
More formally,\oss the objects of\sss $X_{\dff U}$ are 
pairs\sss $(\trf x\fff,\qff \sigma\trf)$\sss such\sss that\sss
$\sigma\qff \in\qff \Sigma^{\dff \fin}$\sss and\sss
$x\qff \in\qff U_{\dff \sigma}$\nsp.\oss
The set\sss of\dss such\sss pairs\dss is\dss ordered\sss by\sss the relation\sss $\leq$\nnsp,\oss
where\sss $(\trf x\fff,\qff \sigma\trf)\off \leq\off (\trf y\fff,\qff \tau\trf)$\trs
if\qss $x\off =\off y$\sss and\sss $\tau\off \subset\off \sigma$\nnsp.\oss
This order defines\sss a structure of\dss a\sss topological\sss category having
this disjoint\sss union as\sss the space of\dss objects and a single morphism\sss
$(\trf x\fff,\qff \sigma\trf)\qff \ttoo\qff (\trf y\fff,\qff \tau\trf)$\trs
if\trs $x\off =\off y$\sss and\sss $\tau\off \subset\off \sigma$\nnsp.\oss
If\dss we consider $X$ as\sss the\sss topological\sss category\sss
having $X$ as\sss the space of\dss objects and only\sss identity\sss morphisms,\oss
then\sss the rule\sss
$(\trf x\fff,\qff \sigma\trf)
\off \longmapsto\off
x$\sss
defines a continuous functor\sss
$\pr\dff \colon\dff
X_{\dff U}\qff \ttoo\qff X$\nnsp,\oss
and\sss hence defines a map\sss\vspace{3pt}
\[
\quad
\num{\pr}\dff \colon\dff
\num{X_{\dff U}}\qff \ttoo\qff \num{X}
\off =\off
X
\qff.
\]

\vspace{-12pt}\vspace{3pt}
In\sss fact,\oss the definition of\dss the\sss
topological\sss category on $\hat{\mathcal{E}}$\sss
is\dss a slightly\sss modified\sss version of\dss this construction.\oss
The following\sss theorem of\qss Segal\pss \cite{s1}\pss is\dss the basic\sss
result\sss about\sss $X_{\dff U}$\nsp.\oss\vspace{1.5pt}

\mypar{Theorem.}{covering-categories}
\emph{If\qss
$U_{\dff a}\dff,\pff a\qff \in\qff \Sigma$\sss
is\dss a numerable covering,\oss then\sss the map\sss 
$\num{\pr}\dff \colon\dff
\num{X_{\dff U}}\qff \ttoo\qff X$\sss
is\dss a homotopy equivalence.\oss
Moreover,\oss 
there\dss exists\dss a homotopy\sss inverse\qss
$s\qff \colon\qff
X\qff \ttoo\qff \num{X_{\dff U}}$\dss
such\sss that\qss
$\num{\pr}\pff \circ\pff s
\off =\off
\id_{\qff X}$\sss
and\dss the composition\qss
$s\pff \circ\pff \num{\pr}$\sss
is\dss homotopic\sss to\sss the identity\sss
in\sss the class of\dss maps 
$f\dff \colon\dff
\num{X_{\dff U}}\qff \ttoo\qff \num{X_{\dff U}}$\dss
such\sss that\qss 
$\num{\pr}\dff \circ\dff f
\off =\off
\num{\pr}$\nnsp.\oss}  \eproof\vspace{1.5pt}

\myuppar{Enhancing\dss families of\dss operators.}
Let\sss $A_{\dff x}\fff,\qff x\qff \in\qff X$\dss be a family\sss of\dss
self-adjoint\trs Fredholm\sss operators,\oss and\sss let\sss
$\mathbb{A}\dff \colon\dff X\qff \ttoo\qff \hat{\mathcal{F}}$\sss
be\sss the corresponding\sss map.\oss
We would\sss like\sss to\sss lift\sss $\mathbb{A}$\sss to a family of\dss
enhanced operators.\oss
This almost\sss never can\sss be done naively,\oss
by\sss lifting $\mathbb{A}$\sss to a map\sss
$X\qff \ttoo\qff \hat{\mathcal{E}}$\dnsp.\oss
The obvious obstruction\dss is\dss the discrete\sss topology on\sss $\rrr$\sss
used\sss in\sss the definition of\sss $\hat{\mathcal{E}}$\dnsp,\oss
but\sss this wouldn't\sss be possible even\sss if\dss we used\sss the standard\sss
topology of\sss $\rrr$\nnsp,\oss
because of\dss the condition 
$-\qff \varepsilon\fff,\qff \varepsilon\qff \not\in\trf \sigma\dff(\dff A\dff)$
in\sss the definition of\dss enhanced operators.\oss
But\sss there are obvious\sss lifts if\sss $\hat{\mathcal{E}}$
is\dss replaced\dss by $\num{\hat{\mathcal{E}}}$ and $X$
is\dss replaced\sss by $\num{X_{\dff U}}$ 
for an appropriate covering $U$\dnsp.\vspace{1.5pt} 

Appropriate coverings $U$ are defined\sss in\sss terms of\dss
the family\sss $A_{\dff x}\fff,\qff x\qff \in\qff X$\nnsp.\oss
For every $x\qff \in\qff X$\sss the operator\sss $A_{\dff x}$\sss is\trs Fredholm.\oss
By\sss a basic property of\qss Fredholm\dss operators,\oss
there\dss exists\sss
$\varepsilon_{\dff x}\qff >\qff 0$\sss such\sss that\sss
$(\trf A_{\dff x}\dff,\qff \varepsilon_{\dff x}\trf)$\sss
is\dss an enhanced\sss operator.\oss
By\sss another\sss basic property of\trs Fredholm\dss operators\qss
this implies\sss that\sss
$(\trf A_{\dff y}\dff,\qff \varepsilon_{\dff x}\trf)$\sss
is\dss an enhanced\sss operator for every $y$ in some neighborhood of\sss $x$\nnsp.\oss
It\sss follows\sss that\sss
there exists an open covering\sss
$U_{\dff a}\dff,\pff a\qff \in\qff \Sigma$\sss of\dss $X$\sss
and\sss positive numbers\sss
$\varepsilon_{\dff a}\dff,\pff a\qff \in\qff \Sigma$\sss
such\sss that\sss
$(\trf A_{\dff z}\dff,\qff \varepsilon_{\dff a}\trf)$\sss
is\dss an enhanced\sss operator for every\sss $z\qff \in\qff U_{\dff a}$\nsp.\oss
Any such covering\sss will\sss work\sss for us.\oss
For every\sss 
$\sigma\qff \in\qff \Sigma^{\dff \fin}$\dss 
let\sss 
$\varepsilon_{\dff \sigma}
\off =\off
\min\nolimits_{\pff a\qff \in\qff \sigma}\qff \varepsilon_{\dff a}$\nsp.\oss
Then\sss the pair\sss
$(\trf A_{\dff z}\dff,\qff \varepsilon_{\dff \sigma}\trf)$\sss
is\dss an enhanced operator for every\sss 
$z\qff \in\qff U_{\dff \sigma}$\nsp,\oss
and\sss 
$\varepsilon_{\dff \tau}\qff \leq\qff \varepsilon_{\dff \sigma}$\dss
if\qss 
$\tau\qff \supset\qff \sigma$\nnsp.\oss
Therefore\sss the rule\sss
$(\trf z\fff,\qff \sigma\trf)
\off \longmapsto\off
(\dff A_{\dff z}\dff,\qff \varepsilon_{\dff \sigma}\trf)$\sss
defines a continuous\sss functor\sss\vspace{1.5pt}
\[
\quad
\mathbb{A}_{\qff U,\qff \bm{\varepsilon}}\qff \colon\qff
X_{\dff U}\qff \ttoo\qff \hat{\mathcal{E}}
\qff,
\]

\vspace{-12pt}\vspace{1.5pt}
and\sss hence a continuous map\sss\vspace{2.25pt}
\[
\quad\num{\mathbb{A}_{\qff U,\qff \bm{\varepsilon}}}
\qff \colon\qff
\num{X_{\dff U}}\qff \ttoo\qff \num{\hat{\mathcal{E}}}
\pff.
\]

\vspace{-12pt}\vspace{2.25pt}
If\trs the map\sss 
$\mathbb{A}\dff \colon\dff X\qff \ttoo\qff \hat{\mathcal{F}}$\sss 
is\dss also considered as a functor,\oss
we get\sss a commutative diagram\vspace{0pt}
\[
\quad
\hspace{1.2em}
\begin{tikzcd}[column sep=spech, row sep=boom]
X_{\dff U}
\arrow[r, "\dis \mathbb{A}_{\qff U,\qff \bm{\varepsilon}}"]
\arrow[d, "\dis \pr\qff"']
&
\hat{\mathcal{E}}
\arrow[d, "\dis \dff \hat{\varphi}"]
\\
X
\arrow[r, "\dis \mathbb{A}"]
&
\hat{\mathcal{F}}
\end{tikzcd}
\]

\vspace{-9pt}\vspace{-0.75pt}
of\dss topological\sss categories and\sss functors,\oss
and a commutative diagram\vspace{3pt}
\[
\quad
\begin{tikzcd}[column sep=spech, row sep=boom]
\protect{\num{X_{\dff U}}}
\arrow[r, "\dis \protect{\num{\mathbb{A}_{\qff U,\qff \bm{\varepsilon}}}}"]
\arrow[d, "\dis \protect{\num{\pr}}\qff"']
&
\protect{\num{\hat{\mathcal{E}}}}
\arrow[d, "\dis \dff \protect{\num{\hat{\varphi}}}"]
\\
X
\arrow[r, "\dis \mathbb{A}"]
&
\hat{\mathcal{F}}
\end{tikzcd}
\]

\vspace{-9pt}\vspace{-0.75pt} 
of\dss topological\sss spaces and continuous maps.\oss
The maps\sss
$\num{\mathbb{A}_{\qff U,\qff \bm{\varepsilon}}}$\sss
are\sss the promised\dss lifts of\dss $\mathbb{A}$\nnsp.\oss
If\dss the covering\sss
$U_{\dff a}\dff,\pff a\qff \in\qff \Sigma$\sss
is\dss numerable,\oss
in\sss particular,\oss if\dss the space $X$\sss is\dss paracompact,\oss
there exist\sss even\sss lifts\sss
$X\qff \ttoo\qff \num{\hat{\mathcal{E}}}$\sss
well\sss defined up\sss to homotopy,\oss
as we will\sss see now.\oss

Suppose now\sss that\sss $X$\sss is\dss paracompact.\oss
Let\sss
$s\qff \colon\qff
X\qff \ttoo\qff \num{X_{\dff U}}$\dss
be\sss the homotopy\sss inverse of\trs
$\num{\pr}\dff \colon\dff
\num{X_{\dff U}}\qff \ttoo\qff X$\sss
such as\sss in\qss Theorem\qss \ref{covering-categories}.\oss
Since\sss
$\num{\pr}\pff \circ\pff s
\off =\off
\id_{\qff X}$\nnsp,\oss
for every\sss $x\qff \in\qff X$\sss the point\sss $s\trf(\trf x\trf)$
belongs\sss to\sss the geometric realization of\dss the simplex corresponding\sss
to a sequence\vspace{3pt}
\[
\quad
(\trf x\fff,\qff \sigma_{\dff 0}\trf)\off \leq\off
(\trf x\fff,\qff \sigma_{\dff 1}\trf)\off \leq\off
\ldots\off \leq\off
(\trf x\fff,\qff \sigma_{\fff n}\trf)
\pff.
\]

\vspace{-12pt}\vspace{3pt}
The functor\sss
$\mathbb{A}_{\qff U,\qff \bm{\varepsilon}}$\sss
takes every\sss $(\trf x\dff,\qff \sigma_{\fff i}\trf)$\sss
to an object\sss of\dss the form\sss
$(\trf A_{\dff x}\dff,\qff \varepsilon_{\trf i}\trf)$\nnsp.\oss
It\sss follows\sss that\sss\vspace{3pt}
\[
\quad
\num{\hat{\varphi}}\qff \circ\qff
\num{\mathbb{A}_{\qff U,\qff \bm{\varepsilon}}}\qff \circ\qff
s\trf(\trf x\trf)
\off =\off
A_{\dff x}
\]

\vspace{-12pt}\vspace{3pt}
for every\sss $x\qff \in\qff X$\sss
and\sss hence\sss
$\num{\mathbb{A}_{\qff U,\qff \bm{\varepsilon}}}\qff \circ\qff s$\dss
is\dss a\trs lift\dss of\sss $\mathbb{A}$\nnsp,\oss
i.e.\qss
$\num{\hat{\varphi}}\qff \circ\qff
\num{\mathbb{A}_{\qff U,\qff \bm{\varepsilon}}}\qff \circ\qff
s
\off =\off
\mathbb{A}$\nnsp.\oss

\mypar{Theorem.}{lift-fixed}
\emph{Suppose\sss that\sss $X$\sss is\dss paracompact.\oss
The homotopy\sss class
of\trs the\dss lift}\dss\vspace{1.5pt}
\[
\quad
\num{\mathbb{A}_{\qff U,\qff \bm{\varepsilon}}}\qff \circ\qff
s
\qff \colon\qff
X\qff \ttoo\qff \num{\hat{\mathcal{E}}}
\]

\vspace{-12pt}\vspace{1.5pt} 
\emph{does not\sss depend on\sss the choices involved.\oss}

\proof
Since $s$\sss is\dss a\sss homotopy\sss inverse of\dss $\num{\pr}$\nnsp,\oss
the homotopy class of\dss 
$\num{\mathbb{A}_{\qff U,\qff \bm{\varepsilon}}}\qff \circ\qff
s$\dss
does not\sss depend on\sss the choice of\sss $s$\nnsp.\oss
Suppose\sss that\sss
$V_{\fff b}\dff,\pff b\qff \in\qff \Sigma\fff'$\sss
is\dss another open covering of\dss $X$\sss
and\sss that\sss the numbers\sss
$\delta_{\dff b}\dff,\pff b\qff \in\qff \Sigma\fff'$\sss
are positive and\sss
such\sss that\sss
$(\trf A_{\dff z}\dff,\qff \delta_{\dff b}\trf)$\sss
is\dss an enhanced\sss operator for every\sss $z\qff \in\qff V_{\fff b}$\nsp.\oss
Let\sss us\sss choose a common\sss refinement\sss of\dss
the coverings\dss
$U_{\dff a}\dff,\pff a\qff \in\qff \Sigma$\dss and\dss
$V_{\fff c}\dff,\pff c\qff \in\qff \Sigma\fff'$\nnsp,\oss
i.e.\qss a\sss covering\dss
$W_{\fff i}\dff,\pff i\qff \in\qff \Xi$\dss
such\sss that\sss for every\sss 
$i\qff \in\qff \Xi$\sss\vspace{3pt}
\[
\quad
W_{\fff i}\off \subset\off U_{\dff a\trf(\dff i\trf)}
\quad
\mbox{and}\quad
W_{\fff i}\off \subset\off V_{\fff b\trf(\dff i\trf)}
\pff,
\]

\vspace{-12pt}\vspace{3pt}
where\sss
$a\dff \colon\dff \Xi\qff \ttoo\qff \Sigma$\sss
and\sss
$b\dff \colon\dff \Xi\qff \ttoo\qff \Sigma\fff'$\sss
are some maps.\oss
The rules\sss\vspace{1.5pt}
\[
\quad
(\trf x\fff,\qff \sigma\trf)
\off \longmapsto\off
(\trf x\fff,\qff a\trf(\trf \sigma\trf)\trf)
\quad
\mbox{and}\quad
(\trf x\fff,\qff \sigma\trf)
\off \longmapsto\off
(\trf x\fff,\qff b\trf(\trf \sigma\trf)\trf)
\]

\vspace{-12pt}\vspace{1.5pt}
define functors\sss\vspace{1.5pt}
\[
\quad
\mathbf{a}\dff \colon\dff
X_{\trf W}\qff \ttoo\qff X_{\trf U}
\quad
\mbox{and}\quad\dff
\mathbf{b}\dff \colon\dff
X_{\trf W}\qff \ttoo\qff X_{\trf V}
\pff.
\]

\vspace{-12pt}\vspace{1.5pt}
Clearly,\oss both\sss compositions\vspace{-1.5pt}
\[
\quad
\begin{tikzcd}[column sep=large, row sep=huge]
X_{\trf W}
\arrow[r, "\dis \mathbf{a}\vphantom{p}"]
&
X_{\trf U}
\arrow[r, "\dis \pr"]
&
X
\end{tikzcd}
\quad
\mbox{and}\quad
\begin{tikzcd}[column sep=large, row sep=huge]
X_{\trf W}
\arrow[r, "\dis \mathbf{b}\vphantom{p}"]
&
X_{\trf V}
\arrow[r, "\dis \pr"]
&
X
\end{tikzcd}
\]

\vspace{-9pt}
are equal\dss to\sss
$\pr\dff \colon\dff
X_{\trf W}\qff \ttoo\qff X$\sss
and\dss if\dss $s$\sss is\dss a homotopy\sss inverse of\dss 
$\num{\pr}\dff \colon\dff
\num{X_{\trf W}}\qff \ttoo\qff \num{X}$\nnsp,\oss 
then\sss $\num{\mathbf{a}}\qff \circ\qff s$\sss and\sss
$\num{\mathbf{b}}\qff \circ\qff s$\sss are homotopy\sss inverses of,\qss
respectively,\oss the maps\dss\vspace{2.2pt}
\[
\quad
\num{\pr}\dff \colon\dff
\num{X_{\trf U}}\qff \ttoo\qff \num{X}
\quad
\mbox{and}\quad
\num{\pr}\dff \colon\dff
\num{X_{\trf V}}\qff \ttoo\qff \num{X}
\pff.
\]

\vspace{-12pt}\vspace{2.2pt}
For each\sss $i\qff \in\qff \Xi$\qss
let\sss\vspace{2.2pt}
\[
\quad
\gamma_{\dff i}
\off =\off 
\min\qff 
\left\{\qff \varepsilon_{\dff a\trf(\dff i\trf)}\dff,\pff \delta_{\dff b\trf(\dff i\trf)} 
\qff\right\}
\pff.
\]

\vspace{-12pt}\vspace{2.2pt}
Then\sss for every\sss $z\qff \in\qff W_{\fff i}$\dss
the pair\sss
$(\trf A_{\dff z}\dff,\qff \gamma_{\dff i}\trf)$\sss
is\dss an enhanced operator 
and\vspace{1.5pt}
\[
\quad
\left(\trf A_{\dff z}\dff,\qff \gamma_{\dff i}\trf\right)
\off \leq\off
\left(\trf A_{\dff z}\dff,\qff \varepsilon_{\dff a\trf(\dff i\trf)}\trf\right)
\dff,\off
\left(\trf A_{\dff z}\dff,\qff \delta_{\dff b\trf(\dff i\trf)}\trf\right)
\]

\vspace{-12pt}\vspace{1.5pt}
Similar\sss inequalities hold\sss for every\sss
$\sigma\qff \in\qff \Xi^{\trf \fin}$\sss
in\sss the role of\dss $i$\nnsp.\oss
It\sss follows\sss that\sss there are canonical\sss natural\dss
transformations\sss
\[
\quad
\mathbb{A}_{\qff W,\qff \bm{\gamma}}
\qff \ttoo\qff
\mathbb{A}_{\qff U,\qff \bm{\varepsilon}}\qff \circ\qff \mathbf{a}
\quad\off
\mbox{and}\quad\off
\mathbb{A}_{\qff W,\qff \bm{\gamma}}
\qff \ttoo\qff
\mathbb{A}_{\qff V,\qff \bm{\delta}}\qff \circ\qff \mathbf{b}
\pff.
\]

\vspace{-12pt}
of\dss functors\sss $X_{\trf W}\qff \ttoo\qff \hat{\mathcal{E}}$\nnsp.\oss
Therefore\sss the map\dss\vspace{0.6pt}
\[
\quad
\num{\mathbb{A}_{\qff W,\qff \bm{\gamma}}}
\dff \colon\dff
\num{X_{\trf W}}\qff \ttoo\qff \hat{\mathcal{E}}
\]

\vspace{-12pt}\vspace{0.6pt}
is\dss homotopic\sss to\sss each of\dss the maps\dss
$\num{\mathbb{A}_{\qff U,\qff \bm{\varepsilon}}\qff \circ\qff \mathbf{a}}$
and\dss
$\num{\mathbb{A}_{\qff V,\qff \bm{\delta}}\qff \circ\qff \mathbf{b}}$\nnsp.\oss
Therefore\sss the diagram\vspace{-6pt}\vspace{0.6pt}
\[
\quad
\begin{tikzcd}[column sep=huge, row sep=huge]
&
X_{\trf U}
\arrow[rd, "\dis \mathbb{A}_{\qff U,\qff \bm{\varepsilon}}"]
&
\\
X_{\trf W}
\arrow[rr, "\dis \mathbb{A}_{\qff W,\qff \bm{\gamma}}"]
\arrow[ru, "\dis \mathbf{a}"]
\arrow[rd, "\dis \mathbf{b}"']
&
&
\hat{\mathcal{E}}\qff,
\\
&
X_{\trf V}
\arrow[ru, "\dis \mathbb{A}_{\qff V,\qff \bm{\delta}}"']
&
\end{tikzcd}
\]

\vspace{-15pt}\vspace{0.6pt}
turns into\sss a\sss homotopy commutative diagram\sss after\sss passing\sss
to classifying spaces.\oss
By\sss taking\sss the maps\sss
$\num{\mathbf{a}}\qff \circ\qff s$\sss and\sss
$\num{\mathbf{b}}\qff \circ\qff s$\sss as homotopy\sss inverses,\oss
we see\sss the\sss lifts resulting\sss from\sss
$\mathbb{A}_{\qff U,\qff \bm{\varepsilon}}$\sss
and\sss from\sss
$\mathbb{A}_{\qff V,\qff \bm{\delta}}$\sss
are homotopic\sss to\sss the\sss lift\sss resulting\sss from\sss
$\mathbb{A}_{\qff W,\qff \bm{\gamma}}$\nsp.\oss
The\sss theorem\sss follows.\oss  \eproof

\myuppar{Analytical\dss index.} 
Now we can define\sss the\qss \emph{analytical\dss index}\pss of\dss the family\sss 
$A_{\dff x}\fff,\qff x\qff \in\qff X$\dss
as\sss the homotopy class of\dss the\qss \emph{index\dss map},\oss
defined as\sss the composition\vspace{-1.5pt}
\[
\quad
\begin{tikzcd}[column sep=boom, row sep=boom]
X
\arrow[r, "\dis s"]
&
X_{\trf U}
\arrow[r, "\dis \mathbb{A}_{\qff U,\qff \bm{\varepsilon}}"]
&
\protect{\num{\hat{\mathcal{E}}}}
\arrow[r, "\dis \protect{\num{\hat{\pi}}}"]
&
\protect{\num{\hat{\mathcal{S}}}}
\qff.
\end{tikzcd}
\]

\vspace{-12pt}
We say\sss that\sss the analytical\dss index\qss
\emph{vanishes}\pss if\dss the index\sss map\dss is\dss homotopic\sss to a constant\sss map.

Of\dss course,\oss the homotopy class of\dss the index\sss map
carries\sss the same information
as\sss the homotopy class\sss of\dss the map\sss
$\num{\mathbb{A}_{\qff U,\qff \bm{\varepsilon}}}\qff \circ\qff s$\sss
and even as\sss the homotopy class of\dss the map\sss
$\mathbb{A}\dff \colon\dff X\qff \ttoo\qff \hat{\mathcal{F}}$\sss
defining\sss the family.\oss
The point\sss of\dss replacing\sss $\mathbb{A}$\sss and\sss
$\num{\mathbb{A}_{\qff U,\qff \bm{\varepsilon}}}\qff \circ\qff s$\sss
by\sss the composition\vspace{0.6pt}
\[
\quad
\num{\hat{\pi}}\qff \circ\qff
\num{\mathbb{A}_{\qff U,\qff \bm{\varepsilon}}}\qff \circ\qff
s
\]

\vspace{-12pt}\vspace{0.6pt}
is\dss that\sss the category\sss $\hat{\mathcal{S}}$\sss and\sss its classifying space\sss
$\num{\hat{\mathcal{S}}}$\sss are defined\sss in\sss terms of\dss only\sss
finitely dimensional\sss subspaces of\dss $H$\nnsp,\oss
in\sss contrast\sss with\sss $\hat{\mathcal{E}}$\sss and\sss $\num{\hat{\mathcal{E}}}$\nnsp.\oss
This allows\sss to define\sss the index\sss map 
under much weaker assumptions about\sss the family\sss
$A_{\dff x}\fff,\qff x\qff \in\qff X$\sss
than\sss the continuity of\dss the corresponding\sss map\sss
$\mathbb{A}\dff \colon\dff X\qff \ttoo\qff \hat{\mathcal{F}}$\dnsp.\oss

One can\sss go one step further and\sss use\sss the canonical\sss
homotopy equivalence between\sss
$\num{\hat{\mathcal{S}}}$\sss and\sss $\num{Q}$\nnsp.\oss
It\dss leads\sss to a well\sss defined\sss homotopy class\sss of\dss maps\sss
$X\qff \ttoo\qff \num{Q}$\nnsp,\oss which can\sss be also considered
as\sss the\qss \emph{analytical\dss index}\pss of\dss the family\sss 
$A_{\dff x}\fff,\qff x\qff \in\qff X$\nnsp.\oss
Since\sss $Q$\sss and\sss $\num{Q}$\sss are defined\sss in\sss terms of\dss 
only\sss the finitely dimensional\dss linear algebra and\sss involve no\sss Hilbert\sss spaces,\oss 
this version has a definite conceptual\sss advantage.\oss
At\sss the same\sss time\sss the category\sss 
$\hat{\mathcal{S}}$\sss and\sss the space\sss  $\num{\hat{\mathcal{S}}}$\sss
are closer\sss to operators in\sss $H$\sss and\sss 
by\sss this reason are easier\sss to work with.\oss

\myuppar{Fredholm\dss families.}
We will\sss say\sss that\sss a family\sss
$A_{\dff x}\fff,\qff x\qff \in\qff X$\sss
of\dss self-adjoint\sss operators\sss
$H\qff \ttoo\qff H$\dss is\dss a\qss
\emph{Fredholm\dss family}\qss
if\dss all\sss operators $A_{\dff x}$ are\dss Fredholm\dss
and\sss for every\sss $x\qff \in\qff X$\sss
there exists\sss 
$\varepsilon\off =\off \varepsilon_{\dff x}\qff >\qff 0$\sss
and a neighborhood\sss $U_{\dff x}$\sss of\sss $x$
with\sss the following\sss properties.\oss
First,\oss the pair\sss
$(\trf A_{\dff y}\dff,\qff \varepsilon\trf)$\sss
should\sss be an enhanced\sss operator\sss 
for every $y\qff \in\qff U_{\dff x}$.\oss
Second,\oss the subspaces\vspace{1.5pt}
\[
\quad
V_{\fff y}
\off =\off
\image P_{\qff [\dff -\dff \varepsilon\fff,\qff \varepsilon\dff]}\trf(\trf A_{\dff y}\trf)
\]

\vspace{-12pt}\vspace{1.5pt}
should\sss depend continuously on\sss $y\qff \in\qff U_{\dff x}$,\oss
and,\oss finally,\oss the operators\sss
$V_{\fff y}
\qff \ttoo\qff
V_{\fff y}$
induced\sss by\sss $A_{\dff y}$\sss
should\sss also depend continuously on\sss $y\qff \in\qff U_{\dff x}$.\oss
The above construction of\dss the index\sss maps applies without\sss any changes\sss to arbitrary\trs
Fredholm\dss families.\oss

We say\sss that\sss families
$A_{\dff x}\fff,\qff x\qff \in\qff X$
and
$B_{\dff x}\fff,\qff x\qff \in\qff X$
are\qss \emph{Fredholm\dss homotopic}\pss if\dss there exists\sss
a\dss Fredholm\dss family\sss
$H_{\trf x\fff,\dff u}$,\qss 
$(\trf x\fff,\qff u\trf)
\qff \in\qff 
X\dff \times\dff [\dff 0\fff,\qff 1\dff]$\sss
such\sss that\sss
$H_{\trf x\fff,\dff 0}\off =\off A_{\dff x}$\sss
and\dss
$H_{\trf x\fff,\dff 1}\off =\off B_{\dff x}$\sss
for every\sss $x\qff \in\qff X$\nnsp.\oss
Clearly,\oss Fredholm\dss homotopic\sss families\sss have\sss the same analytical\dss index,\oss
and\dss if\dss a\sss family\sss is\trs Fredholm\dss homotopic\sss
to a constant\sss one,\oss then\sss its analytical\dss index\sss vanishes.

\myuppar{Topological\sss categories related\sss to\sss Hilbert\dss bundles.}
Let\sss $\mathbb{H}$\sss be a\sss 
locally\sss trivial\dss Hilbert\dss bundle with separable fibers over $X$\nnsp.\oss
It\dss is\dss convenient\sss to\sss think\sss that\sss $\mathbb{H}$\sss is\dss a family\sss
$H_{\dff x}\dff,\qff x\qff \in\qff X$ of\dss Hilbert\sss spaces parameterized\sss by $X$\nnsp.\oss
We are interested\sss in\sss families of\dss self-adjoint\trs
Fredholm\dss operators\sss 
$A_{\dff x}\dff \colon\dff H_{\dff x}\qff \ttoo\qff H_{\dff x}$\sss
parameterized\sss by\sss $x\qff \in\qff X$\nnsp.\oss

In\sss this context\sss finitely dimensional\sss
subspaces of\dss the fixed\dss Hilbert\sss space $H$\sss
should\sss be replaced\sss by\sss 
finitely dimensional\sss subspaces of\dss the fibers\sss 
$H_{\dff x}$\nsp,\dss $x\qff \in\qff X$\nnsp.\oss
The set\sss $\hat{\mathcal{S}}\dff(\trf \mathbb{H}\trf)$\sss 
of\dss such subspaces has a natural\dss topology.\oss
A\sss local\sss trivialization of\sss $\mathbb{H}$\sss over an open set\sss
$U\qff \subset\qff X$\sss allows\sss to identify\sss $H_{\dff x}$\sss
with\sss $x\qff \in\qff U$\sss with\sss $H$\nnsp,\oss
and\sss hence\sss to identify subspaces contained\sss in\sss these fibers\sss
with\sss pairs\sss $(\trf x\fff,\qff V\trf)$\nnsp,\oss where\sss
$x\qff \in\qff U$\sss and\sss $V$\sss is\dss a finitely dimensional\sss
subspace of\sss $H_{\dff x}$\nsp.\oss
This\sss leads\sss to a\sss topology\sss on\sss the set\sss of\dss subspaces
of\dss fibers\sss $H_{\dff x}$\sss with\sss $x\qff \in\qff U$\nnsp.\oss
Since\sss these subspaces are\sss finitely dimensional,\oss
this\sss topology does not\sss depend\sss on\sss the choice of\dss
trivialization.\oss
Therefore\sss these\sss topologies\sss define a\sss
topology on\sss $\hat{\mathcal{S}}\dff(\trf \mathbb{H}\trf)$\dnsp.\oss

One can\sss turn\sss $\hat{\mathcal{S}}\dff(\trf \mathbb{H}\trf)$\sss
into a\sss topological\sss category as follows.\oss
We\sss take\sss $\hat{\mathcal{S}}\dff(\trf \mathbb{H}\trf)$\sss
as\sss the space of\dss objects.\oss
Morphisms\sss $V\qff \ttoo\qff V\fff'$ exists only\dss if\sss
$V\fff,\qff V\fff'$\sss are contained\sss in\sss the same fiber\sss $H_{\dff x}$\nsp,\oss
and\sss in\sss this case morphisms are defined exactly as morphisms of\sss
$\hat{\mathcal{S}}$\sss with\sss $H$\sss replaced\sss by\sss $H_{\dff x}$\nsp.\oss
If\dss we\sss turn\sss $X$ into\sss topological\sss category\sss 
having\sss only\sss the identity\sss morphisms,\oss
then\sss there\dss is\dss a\sss functor\sss
$\eta\dff \colon\dff 
\hat{\mathcal{S}}\dff(\trf \mathbb{H}\trf)
\qff \ttoo\qff
X$\sss
assigning\sss to an object\sss $V$\sss the point\sss $x\qff \in\qff X$\sss
such\sss that\sss $V\qff \subset\qff H_{\dff x}$\nsp.\oss\vspace{-0.05pt}

There\dss is\dss also an analogue\sss 
$Q\dff(\trf \mathbb{H}\trf)$\sss 
of\sss $Q$\nnsp.\oss
The objects of\sss $Q\dff(\trf \mathbb{H}\trf)$\sss are pairs\sss $(\trf x\fff,\qff V\trf)$\sss
such\sss that\sss $x\qff \in\qff X$\sss and\sss $V$\sss is\dss a finitely\sss dimensional\dss
Hilbert\dss space.\oss
The\sss topology on\sss the set\sss of\dss objects\dss is\dss defined\sss by\sss
the\sss topology of\sss $X$\sss and\sss the discrete\sss topology on\sss
the set\sss of\dss spaces.\oss
Morphisms\sss
$(\trf x\fff,\qff V\trf)
\qff \ttoo\qff 
(\trf x'\fff,\qff V\fff'\trf)$\sss
exist\sss only\sss if\sss $x\off =\off x'$\nnsp,\oss
and\sss in\sss this case\sss they are\sss the same as morphisms\sss
$V\qff \ttoo\qff V\fff'$\sss in\sss the category $Q$\nnsp.\oss
The\sss topology on\sss the set\sss of\dss morphisms\dss is\dss defined\sss by\sss
the\sss topology of\sss $X$\sss and\sss the\sss topology on\sss the set\sss
of\dss morphisms of\sss $Q$\nnsp.\oss
Clearly,\qss $Q\dff(\trf \mathbb{H}\trf)$\sss  
depends only\sss on $X$ and\sss not\sss on\sss the bundle $\mathbb{H}$\nnsp.\oss
There\dss is\dss a functor\sss
$\eta\dff \colon\dff 
Q\dff(\trf \mathbb{H}\trf)
\qff \ttoo\qff
X$\sss
assigning\sss to\sss $(\trf x\fff,\qff V\trf)$\sss the point $x$\nnsp,\oss
and a\sss functor\sss
$Q\dff(\trf \mathbb{H}\trf)
\qff \ttoo\qff
Q$\nnsp.\oss
These functors lead\sss to a homeomorphism\sss
$\num{Q\dff(\trf \mathbb{H}\trf)}
\qff \ttoo\qff
X\dff \times\dff \num{Q}$\dss
({\fff}this depends on\sss the assumption\sss that\sss
$X$\sss is\dss compactly\sss generated\fff).\oss
There\dss is\dss also an\sss intermediate category\sss $Q/\fff \mathbb{H}$\sss
and\sss functors\sss
$\hat{\mathcal{S}}\dff(\trf \mathbb{H}\trf)
\off \longleftarrow\off
Q/\fff \mathbb{H}
\qff \ttoo\qff
Q\dff(\trf \mathbb{H}\trf)$\sss
relating\sss $\hat{\mathcal{S}}\dff(\trf \mathbb{H}\trf)$\sss
with\sss $Q\dff(\trf \mathbb{H}\trf)$\nnsp.\oss
The geometric realizations\sss\vspace{1.5pt}
\[
\quad
\num{\hat{\mathcal{S}}\dff(\trf \mathbb{H}\trf)}
\off \longleftarrow\off
\num{Q/\fff \mathbb{H}}
\qff \ttoo\qff
\num{Q\dff(\trf \mathbb{H}\trf)}
\]

\vspace{-12pt}\vspace{1.5pt}
of\trs these functors are homotopy equivalences.\oss
See\qss \cite{i2},\oss Theorem\qss 17.1.\oss
The obvious advantage of\dss the classifying space\sss 
$\num{Q\dff(\trf \mathbb{H}\trf)}$\sss
is\trs its\dss independence on\sss the bundle\sss $\mathbb{H}$\nnsp.\oss
But,\oss as we will\sss see,\oss the classifying space\sss 
$\num{\hat{\mathcal{S}}\dff(\trf \mathbb{H}\trf)}$\sss
does not\sss depend on\sss $\mathbb{H}$\sss in any essential\sss way\sss either.\oss

\myuppar{Fredholm\dss families of\dss operators in\dss Hilbert\dss bundles.}
The definitions of\qss \emph{Fredholm\dss families}\pss of\dss operators
and of\qss \emph{Fredholm\dss homotopies}\pss apply,\oss
with obvious modifications,\oss to families\sss of\dss operators\sss
$A_{\dff x}\dff \colon\dff H_{\dff x}\qff \ttoo\qff H_{\dff x}$\nsp,\qss
$x\qff \in\qff X$\nnsp.\oss
But\sss the definition of\dss the analytical\dss index should\sss be modified.\oss
In\sss this situation\sss there\dss is\dss no good analogue of\dss
the space\sss $\hat{\mathcal{F}}$\sss and\sss a family cannot\sss be considered
as a continuous map\sss to some space.\oss
There\dss is\dss no good analogue of\dss
the\sss category\sss $\hat{\mathcal{E}}$ either.\oss
But\sss the construction of\dss the index\sss maps\sss
can\sss be easily\sss modified\sss to get\sss a map\sss
$X\qff \ttoo\qff \num{\hat{\mathcal{S}}\dff(\trf \mathbb{H}\trf)}$\sss
well\sss defined\sss up\sss to homotopy,\oss and\dss hence a\sss homotopy class\sss 
$X\qff \ttoo\qff \num{Q\dff(\trf \mathbb{H}\trf)}$\nnsp.

The first\sss steps are\sss the same as for a fixed\dss Hilbert\dss space.\oss
Let\sss
$A_{\dff x}\dff \colon\dff H_{\dff x}\qff \ttoo\qff H_{\dff x}$\nsp,\qss
$x\qff \in\qff X$\dss
be a\dss Fredholm\dss family.\oss 
Then\sss
there exists an open covering\sss
$U_{\dff a}\dff,\pff a\qff \in\qff \Sigma$\sss of\dss $X$\sss
and\sss numbers\sss
$\varepsilon_{\dff a}\qff >\qff 0\dff,\pff a\qff \in\qff \Sigma$\sss
such\sss that\sss
$(\trf A_{\dff z}\dff,\qff \varepsilon_{\dff a}\trf)$\sss
is\dss an enhanced\sss operator\sss in\sss $H_{\dff z}$\sss 
if\dss $z\qff \in\qff U_{\dff \alpha}$\nsp.\oss
As before,\oss for every\sss 
$\sigma\qff \in\qff \Sigma^{\dff \fin}$\dss let\sss 
$\varepsilon_{\dff \sigma}
\off =\off
\min\nolimits_{\pff a\qff \in\qff \sigma}\qff \varepsilon_{\dff a}$\nsp.\oss
Then\sss the pair\sss
$(\trf A_{\dff z}\dff,\qff \varepsilon_{\dff \tau}\trf)$\sss
is\dss an enhanced operator for\sss every $z\qff \in\qff U_{\dff \tau}$\nsp,\oss
and\sss 
$\varepsilon_{\dff \tau}\qff \leq\qff \varepsilon_{\dff \sigma}$\dss
if\qss $\tau\qff \supset\qff \sigma$\nnsp.\oss
Let\sss us\sss construct\sss a continuous functor\sss\vspace{1.5pt}\vspace{1pt}
\[
\quad
\mathbb{A}_{\qff U,\qff \bm{\varepsilon}}\dff \colon\dff
X_{\dff U}\qff \ttoo\qff \hat{\mathcal{S}}\dff(\trf \mathbb{H}\trf)
\pff.
\]

\vspace{-12pt}\vspace{1.5pt}\vspace{1pt}
Suppose\sss that
$(\trf z\fff,\qff \tau\trf)$ 
is\dss an object\sss of\dss $X_{\dff U}$
and\sss let\sss
$\varepsilon\off =\off \varepsilon_{\dff \tau}$.\oss
The\sss functor $\mathbb{A}_{\qff U,\qff \bm{\varepsilon}}$\sss
takes $(\trf z\fff,\qff \tau\trf)$
to\vspace{1.5pt}\vspace{1pt}
\[
\quad
\image P_{\qff [\dff -\dff \varepsilon\fff,\qff \varepsilon\dff]}\trf(\trf A_{\dff z}\trf)
\off \in\off\dff
\ob\dff \hat{\mathcal{S}}\dff(\trf \mathbb{H}\trf)
\pff.
\]

\vspace{-12pt}\vspace{1.5pt}\vspace{1pt}
Suppose now\sss that\sss
$(\trf z\fff,\qff \tau\trf)\qff \ttoo\qff (\trf z\fff,\qff \sigma\trf)$\sss
is\dss a morphism
of\trs $X_{\dff U}$\nsp,\oss
and\dss let\sss
$\varepsilon\off =\off \varepsilon_{\dff \tau}$\sss
and\sss
$\varepsilon'\off =\off \varepsilon_{\dff \sigma}$.\oss
The\sss functor\sss $\mathbb{A}_{\qff U,\qff \bm{\varepsilon}}$\sss
takes\sss this morphism\sss
to\sss the morphism\vspace{1.5pt}\vspace{1pt}
\[
\quad
\image P_{\qff [\dff -\dff \varepsilon\fff,\qff \varepsilon\dff]}\trf(\trf A_{\dff z}\trf)
\off \ttoo\off
\image P_{\qff [\dff -\dff \varepsilon'\fff,\qff \varepsilon'\dff]}\trf(\trf A_{\dff z}\trf)
\off
\]

\vspace{-12pt}\vspace{1.5pt}\vspace{1pt}
of\dss the category\sss $\hat{\mathcal{S}}\dff(\trf \mathbb{H}\trf)$\sss
corresponding\sss to\sss the decomposition\vspace{3pt}
\[
\quad
\image P_{\dff [\dff -\dff \varepsilon'\fff,\qff \varepsilon'\dff]}\trf(\trf A_{\dff z}\trf)
\off\qff =\off\qff
\left(\qff
\image P_{\dff [\dff -\dff \varepsilon'\fff,\qff -\dff \varepsilon\dff]}\trf(\trf A_{\dff z}\trf)
\qff\right)
\qff \oplus\qff
\left(\qff
\image P_{\dff [\dff -\dff \varepsilon\fff,\qff \varepsilon\dff]}\trf(\trf A_{\dff z}\trf)
\qff\right)
\qff \oplus\qff
\left(\qff
\image P_{\dff [\dff \varepsilon\fff,\qff -\dff \varepsilon'\dff]}\trf(\trf A_{\dff z}\trf)
\qff\right)
\pff.
\]

\vspace{-12pt}\vspace{3pt}
Clearly,\qss $\mathbb{A}_{\qff U,\qff \bm{\varepsilon}}$\sss is\dss indeed\sss a\sss functor\sss
$X_{\dff U}\qff \ttoo\qff \hat{\mathcal{S}}\dff(\trf \mathbb{H}\trf)$\nnsp.\oss
Since\sss the subspaces\vspace{1.5pt}
\[
\quad
V_{\fff z}
\off =\off
\image P_{\dff [\dff -\dff \varepsilon\fff,\qff \varepsilon\dff]}\trf(\trf A_{\dff z}\trf)
\pff
\]

\vspace{-12pt}\vspace{1.5pt}
depend continuously on\sss $z\qff \in\qff U_{\dff \tau}$,\oss
it\dss is\dss continuous on objects,\oss
and\dss since\sss the induced operators\sss
$V_{\fff z}\qff \ttoo\qff V_{\fff z}$\sss also depend continuously on $z$\nnsp,\oss
it\dss is\dss continuous on\sss morphisms.\oss
By\sss passing\sss to\sss the geometric realizations\sss the\qss
\emph{index\dss functor}\dss
$\mathbb{A}_{\qff U,\qff \bm{\varepsilon}}$\sss
defines a continuous map\vspace{1.5pt}
\[
\quad
\num{\mathbb{A}_{\qff U,\qff \bm{\varepsilon}}}\qff \colon\qff
\num{X_{\dff U}}
\qff \ttoo\qff 
\num{\hat{\mathcal{S}}\dff(\trf \mathbb{H}\trf)}
\pff.
\]

\vspace{-12pt}
\mypar{Theorem.}{bundle-correct}
\emph{Suppose\sss that\sss $X$\sss is\dss paracompact\sss
and\sss $s\dff \colon\dff X\qff \ttoo\qff \num{X_{\dff U}}$\sss
is\dss a homotopy\sss inverse of\pss $\num{\pr}$
as\sss in\qss Theorem\qss \ref{covering-categories}.\oss
Then\sss the homotopy class of\pss 
$\num{\mathbb{A}_{\qff U,\qff \bm{\varepsilon}}}\qff \circ\qff
s
\trf \colon\dff
X
\qff \ttoo\qff 
\num{\hat{\mathcal{S}}\dff(\trf \mathbb{H}\trf)}$\sss
does not\sss depend on\sss the choices\sss involved.\oss}

\proof
The proof\trs is\dss completely similar\sss to\sss the proof\dss of\trs
Theorem\qss \ref{lift-fixed}.\oss  \eproof

\myuppar{Analytical\dss index\sss of\dss families 
of\dss operators\sss in\dss Hilbert\dss bundles.}
Suppose\sss that\sss $X$\sss is\dss paracompact.\oss
Then\sss any\sss map\sss of\dss the form\dss
$\num{\mathbb{A}_{\qff U,\qff \bm{\varepsilon}}}\qff \circ\qff
s
\trf \colon\dff
X
\qff \ttoo\qff 
\num{\hat{\mathcal{S}}\dff(\trf \mathbb{H}\trf)}$\sss
may\sss be called an\qss \emph{index\sss map}\pss
of\dss the family\sss
$A_{\dff x}\dff,\qff
x\qff \in\qff X$\nnsp.\oss
By\trs Theorem\qss \ref{bundle-correct}\qss the homotopy class
of\dss an index\sss map\dss is\dss independent\sss from choices 
involved\sss in\sss its construction,\oss
and\sss we can define\sss the\qss \emph{index}\pss of\dss the family\sss
$A_{\dff x}\dff,\qff
x\qff \in\qff X$\dss
as\sss the homotopy class of\dss any of\dss its\sss index\sss maps.\oss

A drawback of\dss this definition\dss is\dss that\sss such index\sss maps
are analogues not\sss of\dss the index maps\sss
$X\qff \ttoo\qff \num{\hat{\mathcal{S}}}$\sss of\dss families of\dss operators in a
fixed\dss Hilbert\dss spaces,\oss
but\sss of\dss maps\sss
$X\qff \ttoo\qff X\dff \times\dff \num{\hat{\mathcal{S}}}$\sss
having\sss the identity\sss map\sss $X\qff \ttoo\qff X$\sss as\sss the first\sss
component\sss and an\sss index\sss map as\sss the second.\oss

A simple way\sss to deal\sss with\sss this issue\dss is\dss to use\sss the
homotopy equivalence between\sss
$\num{\hat{\mathcal{S}}\dff(\trf \mathbb{H}\trf)}$\sss 
and\sss 
$\num{Q\dff(\trf \mathbb{H}\trf)}$\nnsp,\oss
which\dss is\dss well\sss defined up\sss to homotopy.\oss
By composing\sss an\sss index\sss map\sss in\sss the above sense with\sss
this homotopy equivalence we get\sss a well\sss defined\sss homotopy class of\dss maps
$X\qff \ttoo\qff \num{Q\dff(\trf \mathbb{H}\trf)}$\nnsp.\oss
Since\sss
$\num{Q\dff(\trf \mathbb{H}\trf)}$\sss
is\dss canonically\sss homeomorphic\sss to\sss
$X\dff \times\dff \num{Q}$\nnsp,\oss
we can consider\sss it\sss as a homotopy class of\dss maps
$X\qff \ttoo\qff X\dff \times\dff \num{Q}$\nnsp.\oss
The construction of\dss the homotopy equivalence between\sss
$\num{\hat{\mathcal{S}}\dff(\trf \mathbb{H}\trf)}$\sss 
and\sss 
$\num{Q\dff(\trf \mathbb{H}\trf)}$\dss
shows\sss that\sss the first\sss component\sss of\dss the homotopy class\sss
$X\qff \ttoo\qff X\dff \times\dff \num{Q}$\sss
is\dss the homotopy class of\dss the identity map\sss
$X\qff \ttoo\qff X$\nnsp.\oss
Therefore\sss the second component\sss
$X\qff \ttoo\qff \num{Q}$\sss
carries\sss the same information as\sss the homotopy classes\sss\vspace{1.5pt}
\[
\quad
X\qff \ttoo\qff \num{Q\dff(\trf \mathbb{H}\trf)}
\quad
\mbox{and}\quad
X\qff \ttoo\qff \num{\hat{\mathcal{S}}}
\pff.
\]

\vspace{-12pt}\vspace{1.5pt}
This suggests\sss to define\sss the\qss \emph{index}\pss of\dss the family\sss
$A_{\dff x}\dff,\qff
x\qff \in\qff X$\dss
as\sss the homotopy class
$X\qff \ttoo\qff \num{Q}$\sss and\sss to call\qss \emph{index\sss maps}\pss
the maps in\sss this homotopy class.\oss
This approach\sss has\sss the advantage of\dss the\sss target\sss space\sss
$\num{Q}$\sss being\dss independent\sss from\sss the bundle\sss $\mathbb{H}$\nnsp.\oss

Nevertheless,\oss the\sss homotopy class of\dss the\sss index\sss maps\sss
$X
\qff \ttoo\qff 
\num{\hat{\mathcal{S}}\dff(\trf \mathbb{H}\trf)}$\sss 
carries\sss the same\sss information as\sss the homotopy class\sss
$X\qff \ttoo\qff \num{Q}$\sss 
and sometimes\dss 
is\dss easier\sss to work\sss with.\oss
Suppose\sss that\sss $\mathbb{H}\fff,\qff \mathbb{K}$\sss are\sss two\dss
Hilbert\dss bundles over\sss $X$\nnsp.\oss 
Then an\sss isomorphism \sss 
$\mathbb{H}\qff \ttoo\qff \mathbb{K}$\sss
covering\sss the identity\sss $X\qff \ttoo\qff X$\sss
induces an\sss isomorphism of\dss topological\sss categories\sss\vspace{1.5pt}
\[
\quad
\hat{\mathcal{S}}\dff(\trf \mathbb{H}\trf)
\qff \ttoo\qff
\hat{\mathcal{S}}\dff(\trf \mathbb{K}\trf)
\pff
\]

\vspace{-12pt}\vspace{1.5pt}
and\sss hence a homeomorphism\sss
$\num{\hat{\mathcal{S}}\dff(\trf \mathbb{H}\trf)}
\qff \ttoo\qff
\num{\hat{\mathcal{S}}\dff(\trf \mathbb{K}\trf)}$\nnsp.\oss
At\sss the same\sss time\sss it\dss is\dss known\sss that\sss
every\dss Hilbert\sss bundle\dss is\dss actually\sss trivial,\oss
i.e.\qss is\dss isomorphic\sss to\sss the bundle\sss
$X\dff \times\dff H\qff \ttoo\qff X$\nnsp.\oss
See\sss the discussion\sss in\dss Section\qss \ref{strictly-fredholm-families}.\oss
It\sss follows\sss that\sss
$\num{\hat{\mathcal{S}}\dff(\trf \mathbb{H}\trf)}$\sss
is\dss homeomorphic\sss to\sss
$X\dff \times\dff \num{\hat{\mathcal{S}}}$\nnsp.\oss
Moreover,\oss the arguments proving\sss the\sss triviality show\sss
that\sss any\sss two homeomorphisms obtained\sss in\sss this way are homotopic\qss
({\fff}in\sss fact,\oss even\sss isotopic).\oss
Therefore\sss the\qss \emph{index}\pss can\sss be defined also 
as\sss the second\sss component\sss
$X\qff \ttoo\qff \num{\hat{\mathcal{S}}}$\sss 
of\dss the homotopy class\sss\vspace{1.5pt}
\[
\quad
X
\qff \ttoo\qff 
\num{\hat{\mathcal{S}}\dff(\trf \mathbb{H}\trf)}
\qff \ttoo\qff 
X\dff \times\dff \num{\hat{\mathcal{S}}}
\pff.
\]

\vspace{-12pt}
\mysection{Analytical\qss index\qss of\pss Fredholm\qss families}{analytic-index-fredholm-non-sa}

\myuppar{Topological\sss categories\sss related\dss to\dss Fredholm\sss operators.}
This section\dss is\dss devoted\sss to\sss families of\dss general\trs Fredholm\dss
operators.\oss
Let\sss $\mathcal{F}$\sss be\sss the space of\trs Fredholm\dss operators\sss $H\qff \ttoo\qff H$\sss
with\sss the norm\sss topology.\oss
Our\sss first\sss goal\dss is\dss to define analogues of\dss
$\hat{\mathcal{E}}$\sss and\sss $\hat{\mathcal{S}}$\dnsp.\oss
Recall\dss that\sss the\qss 
\emph{polar\sss decomposition}\pss of\dss an operator\sss $B$\sss is\dss the
unique presentation\sss
$B\off =\off U\trf \num{B}$\nnsp,\oss
where\vspace{0pt}
\[
\quad
\num{B}
\off =\off
\sqrt{\dff B^*\dff B\dff}
\]

\vspace{-12pt}\vspace{0pt}
and\sss $U$\sss is\dss a partial\sss isometry of\sss $H$\sss
with\sss $\kernel U\off =\off \kernel B$\sss
and\dss
$\image U$\sss is\dss equal\sss to\sss the closure of\sss $\image B$\nnsp.\oss
If\sss $B$\sss is\trs Fredholm,\oss
then\sss $\image B$\sss is\dss automatically\sss closed.\oss

An\qss \emph{enhanced\trs Fredholm\dss operator}\pss is\dss
a pair\sss $(\trf B\dff,\qff \varepsilon\trf)$\nnsp,\oss
where $B\qff \in\qff \mathcal{F}$ and $\varepsilon\qff \in\qff \rrr$\sss
are such\sss that\sss $\varepsilon\qff >\qff 0$\nnsp,\oss the\sss interval\sss
$[\trf 0\fff,\qff \varepsilon\trf]$\sss is\dss disjoint\sss from\sss the essential\sss
spectrum of\sss $\num{B}$\nnsp,\oss
and\sss $\varepsilon\qff \not\in\qff \sigma\dff(\trf \num{B}\trf)$\nnsp.\oss
Let\sss $\mathcal{E}$ be\sss the space\sss of\dss enhanced\trs Fredholm\dss operators.\oss
The\sss topology\dss is\dss defined\sss by\sss
the\sss topology of\sss $\mathcal{F}$\sss
and\sss the discrete\sss topology on\sss $\rrr$\nnsp.\oss
The space\sss $\mathcal{E}$\sss 
is\dss ordered\sss by\sss the relation\sss\vspace{1.5pt}
\[
\quad
(\trf B\dff,\qff \varepsilon\trf)
\off \leq\off
(\trf B'\dff,\qff \varepsilon'\trf)
\quad
\mbox{if}\quad 
B\off =\off B'
\quad
\mbox{and}\quad 
\varepsilon\qff \leq\qff \varepsilon'
\pff.
\]

\vspace{-12pt}\vspace{1.5pt}
This order allows\sss to consider\sss $\mathcal{E}$ as a\sss
topological\sss category.\oss
The obvious forgetting\sss functor\sss
$\varphi\dff \colon\dff \mathcal{E}\qff \ttoo\qff \mathcal{F}$\sss
induces a homotopy equivalence\sss
$\num{\varphi}\dff \colon\dff
\num{\mathcal{E}}\qff \ttoo\qff \num{\mathcal{F}}$\nnsp.\oss
See\qss \cite{i2},\oss Section\qss 14.\oss

The category\sss $\mathcal{S}$\sss is\dss an analogue of\sss
$\hat{\mathcal{S}}$\dnsp.\oss
Its objects are\sss pairs\sss
$(\trf E_{\dff 1}\dff,\pff E_{\trf 2}\trf)$\sss
such\sss that\sss $E_{\dff 1}\dff,\qff E_{\trf 2}$ are\sss finitely\sss
dimensional\sss subspaces of\sss $H$\nnsp,\oss
with\sss morphisms\sss
$(\trf E_{\dff 1}\fff,\qff E_{\dff 2}\trf)
\qff \ttoo\qff 
(\trf E\fff'_{\dff 1}\fff,\qff E\fff'_{\dff 2}\trf)$\sss
being\sss pairs of\dss subspaces\sss
$(\trf F_{\dff 1}\fff,\qff F_{\dff 2}\trf)$
such\sss that\vspace{1.5pt}
\[
\quad
E_{\dff 1}\dff \oplus\dff F_{\dff 1}
\off =\off
E\fff'_{\dff 1}
\quad
\mbox{and}\quad
E_{\dff 2}\dff \oplus\dff F_{\dff 2}
\off =\off
E\fff'_{\dff 2}
\]

\vspace{-12pt}\vspace{1.5pt}
together with an\sss isometry\sss
$f\dff \colon\dff
F_{\dff 1}\qff \ttoo\qff F_{\dff 2}$\nsp.\oss
The composition\dss is\dss defined\sss by\sss taking\sss the direct\sss sums
of\dss the corresponding subspaces\sss $F_{\dff 1}\fff,\qff F_{\dff 2}$
and of\dss the isometries.\oss
The category $\mathcal{S}$\sss is\dss a\sss topological\sss category\sss 
in\sss an obvious way.\oss
Like\sss $\hat{\mathcal{S}}$\dnsp,\oss
it\dss is\dss not\sss associated\sss
with a\sss partial\sss order.\oss
Let\sss $\mathcal{S}_{\qff 0}$\sss
be\sss the full\sss subcategory of\dss $\mathcal{S}$ 
defined\sss by\sss the condition\sss
$\dim\dff E_{\dff 1}\off =\off \dim\dff E_{\trf 2}$\nsp.\oss

There\dss is\dss functor\sss
$\pi\dff \colon\dff
\mathcal{E}
\qff \ttoo\qff 
\mathcal{S}$\trs
taking\sss an object\sss 
$(\trf B\fff,\qff \varepsilon\trf)$\sss of\sss $\mathcal{E}$\sss to\sss
the pair\sss 
$(\trf E_{\dff 1}\dff,\pff E_{\trf 2}\trf)$\nnsp,\oss
where\vspace{1.5pt}
\[
\quad
E_{\dff 1}
\off =\off
P_{\qff [\dff 0\fff,\qff \varepsilon\dff]}\trf(\trf \num{B}\trf)\fff,\qquad
E_{\trf 2}
\off =\off
P_{\qff [\dff 0\fff,\qff \varepsilon\dff]}\trf(\trf \num{B^*}\trf)
\qff.
\] 

\vspace{-12pt}\vspace{1.5pt}
The functor\sss $\pi$\sss takes
a morphism\sss
$(\trf B\dff,\qff \varepsilon\trf)
\qff \ttoo\qff
(\trf B\dff,\qff \varepsilon'\trf)$\nnsp,\oss
where\sss $\varepsilon\qff \leq\qff \varepsilon'$\nnsp,\oss
to\sss the morphism\vspace{1.5pt}\vspace{-0.625pt}
\[
\quad
\pi\qff(\trf B\dff,\qff \varepsilon\trf)
\off \ttoo\off
\pi\qff(\trf B\dff,\qff \varepsilon'\trf)
\] 

\vspace{-12pt}\vspace{1.5pt}\vspace{-0.625pt}
of\sss $\mathcal{S}$\sss
defined\sss by\sss the pair of\dss subspaces\vspace{1.5pt}
\[
\quad
(\trf F_{\dff 1}\fff,\qff F_{\dff 2}\trf)
\off =\off
\left(\qff
P_{\qff [\dff \varepsilon\fff,\qff \varepsilon'\dff]}\trf(\trf \num{B}\trf)\fff,\off
P_{\qff [\dff \varepsilon\fff,\qff \varepsilon'\dff]}\trf(\trf \num{B^*}\trf)
\qff\right)
\] 

\vspace{-12pt}\vspace{1.5pt}
together\sss with\sss the isometry\sss\vspace{1.5pt}
\[
\quad
P_{\qff [\dff \varepsilon\fff,\qff \varepsilon'\dff]}\trf(\trf \num{B}\trf)
\off \ttoo\off
P_{\qff [\dff \varepsilon\fff,\qff \varepsilon'\dff]}\trf(\trf \num{B^*}\trf)
\] 

\vspace{-12pt}\vspace{1.5pt}
induced\sss by\sss $U$\nnsp,\oss where\sss
$B\off =\off U\trf \num{B}$\sss is\dss the polar decomposition of\sss $B$\nnsp.\oss
A routine check shows\sss that\sss these rules define a functor\sss
$\pi\trf \colon\dff
\mathcal{E}\qff \ttoo\trf \mathcal{S}$\nsp\dnsp.\oss
The geometric realization\sss
$\num{\pi}\dff \colon\dff
\num{\mathcal{E}}
\qff \ttoo\qff 
\num{\mathcal{S}}$\sss
is\dss a homotopy equivalence,\oss 
as\sss it\sss easily\sss follows\sss from\sss the results of\qss \cite{i2},\oss
Sections\qss 14\qss and\qss 15.\oss

There\dss is\dss also an analogue\sss $S$\sss of\dss the category\sss $Q$\nnsp.\oss
The objects of\sss $S$
are\sss pairs\sss 
$(\trf V_{\fff 1}\dff,\pff V_{\fff 2}\trf)$ 
of\dss finite\-ly\sss dimensional\dss Hilbert\sss spaces,\oss
with\sss morphisms\sss
$(\trf V_{\dff 1}\dff,\pff V_{\dff 2}\trf)
\qff \ttoo\qff 
(\trf W_{\dff 1}\dff,\pff W_{\dff 2}\trf)$\sss
being\sss triples\sss
$(\trf i_{\trf 1}\dff,\pff i_{\trf 2}\dff,\pff g\trf)$\nnsp,\oss
where\dss 
$i_{\trf 1}\dff \colon\dff V_{\fff 1}\qff \ttoo\qff W_{\fff 1}$\dss
and\dss
$i_{\trf 2}\dff \colon\dff V_{\fff 2}\qff \ttoo\qff W_{\fff 2}$\dss
are isometric embeddings and\vspace{1.5pt}
\[
\quad
g\dff \colon\dff
W_{\fff 1}\qff \ominus\qff i_{\trf 1}\dff(\trf V_{\fff 1}\trf)
\qff \ttoo\qff
W_{\fff 2}\qff \ominus\qff i_{\trf 2}\dff(\trf V_{\fff 2}\trf)
\]

\vspace{-12pt}\vspace{1.5pt}
is\dss an\sss isometry.\oss
The composition\dss is\dss defined\sss in\sss an obvious way
and amounts\sss to\sss taking\sss the direct\sss sum of\dss
the isometric\sss isomorphisms $g$\nnsp.\oss
As in\sss the case of\dss self-adjoint\sss operators,\oss
there\dss is\dss an\sss intermediate category\sss $S/\fff H$\sss
and continuous functors\sss
$\mathcal{S}\off \longleftarrow\off S/\fff H\qff \ttoo\qff S$\sss
inducing\sss homotopy equivalences\sss
$\num{\mathcal{S}}
\off \longleftarrow\off 
\num{S/\fff H}\qff \ttoo\qff \num{S}$\nnsp.\oss

\myuppar{Families of\trs Fredholm\dss operators.}
In\dss this section\sss we will\sss say\sss that\sss a family\sss
$B_{\dff x}\fff,\qff x\qff \in\qff X$\sss
of\dss operators\sss
$H\qff \ttoo\qff H$\dss is\dss a\qss
\emph{Fredholm\dss family}\pss
if\dss operators $B_{\dff x}$ are\dss Fredholm\qss
({\halfff}but\sss are not\sss assumed\sss to be self-adjoint\halfff)\qss
and\sss for every\sss $x\qff \in\qff X$\sss
there exists\sss 
$\varepsilon\off =\off \varepsilon_{\dff x}\qff >\qff 0$\sss
and a neighborhood\sss $U_{\dff x}$\sss of\sss $x$
with\sss the following\sss properties.\oss
First,\oss the pair\sss
$(\trf B_{\dff y}\dff,\qff \varepsilon\trf)$\sss
should\sss be an enhanced\dss Fredholm\dss operator\sss 
for every $y\qff \in\qff U_{\dff x}$.\oss
Second,\oss the subspaces\vspace{1.5pt}\vspace{1.05pt}
\[
\quad
E_{\dff 1}\dff(\trf y\trf)
\off =\off
\image\dff P_{\qff [\dff 0\fff,\qff \varepsilon\dff]}\trf(\trf \num{B_{\dff y}}\trf)\fff,\qquad
E_{\trf 2}\dff(\trf y\trf)
\off =\off
\image\dff P_{\qff [\dff 0\fff,\qff \varepsilon\dff]}\trf(\trf \num{B^*_{\dff y}}\trf)
\]

\vspace{-12pt}\vspace{1.5pt}\vspace{1.05pt}
should\sss depend continuously on\sss $y\qff \in\qff U_{\dff x}$,\oss
and,\oss finally,\oss the operators\sss
\[
E_{\dff 1}\dff(\trf y\trf)
\qff \ttoo\qff
E_{\trf 2}\dff(\trf y\trf)
\]
induced\sss by\sss $B_{\dff y}$
should\sss also depend continuously on\sss $y\qff \in\qff U_{\dff x}$.\oss

\myuppar{Index of\qss Fredholm\dss families.}
Let\sss $B_{\dff x}\fff,\qff x\qff \in\qff X$\dss be a family\sss of\trs
Fredholm\sss operators,\oss and\sss let\sss
$\mathbb{B}\dff \colon\dff X\qff \ttoo\qff \mathcal{F}$\sss
be\sss the corresponding\sss map.\oss
As in\sss the self-adjoint\sss case,\oss
the map\sss $\mathbb{B}$\sss can\sss be\qss ``lifted''\qss to a map\sss
to\sss $\num{\mathcal{E}}$\nnsp.\oss
In\sss more details,\oss there exists an open covering\sss
$U_{\dff a}\dff,\pff a\qff \in\qff \Sigma$\sss of\dss $X$\sss
and\sss positive numbers\sss
$\varepsilon_{\dff a}\dff,\pff a\qff \in\qff \Sigma$\sss
such\sss that\sss
$(\trf B_{\dff z}\dff,\qff \varepsilon_{\dff a}\trf)$\sss
is\dss an enhanced\trs Fredholm\dss operator for every\sss $z\qff \in\qff U_{\dff a}$\nsp.\oss
The families\sss 
$U_{\dff a}\dff,\pff a\qff \in\qff \Sigma$\dss 
and\dss 
$\varepsilon_{\dff a}\dff,\pff a\qff \in\qff \Sigma$\sss  
define a continuous functor\vspace{1.5pt}
\[
\quad
\mathbb{B}_{\qff U,\qff \bm{\varepsilon}}\qff \colon\qff
X_{\dff U}\qff \ttoo\qff \mathcal{E}
\qff
\]

\vspace{-12pt}\vspace{1.5pt}
exactly as in\dss Section\qss \ref{analytic-index-section}.\oss
The geometric realization of\dss $\mathbb{B}_{\qff U,\qff \bm{\varepsilon}}$\sss
is\dss a continuous map\sss\vspace{1.5pt}
\[
\quad\num{\mathbb{B}_{\qff U,\qff \bm{\varepsilon}}}
\qff \colon\qff
\num{X_{\dff U}}\qff \ttoo\qff \num{\mathcal{E}}
\pff.
\]

\vspace{-12pt}\vspace{1.5pt}
If\sss $X$\sss is\dss paracompact,\oss
then one can\sss take\sss the homotopy\sss inverse\sss
$s\qff \colon\qff
X\qff \ttoo\qff \num{X_{\dff U}}$\dss
of\trs the canonical\sss map\sss
$\num{\pr}\dff \colon\dff
\num{X_{\dff U}}\qff \ttoo\qff X$\sss
and\sss consider\sss the composition\vspace{1.5pt}
\[
\quad\num{\mathbb{B}_{\qff U,\qff \bm{\varepsilon}}}\qff \circ\qff
s
\qff \colon\qff
X\qff \ttoo\qff \num{\mathcal{E}}
\qff.
\]

\vspace{-12pt}\vspace{1.5pt}
Its homotopy class does not\sss depend on\sss the choices involved.\oss
The proof\dss of\dss this fact\dss is\dss completely\sss similar\sss to\sss
the proof\dss of\trs Theorem\qss \ref{lift-fixed}.\oss
Now one can define\sss the\qss \emph{analytical\dss index}\pss of\dss the family\sss 
$B_{\dff x}\fff,\qff x\qff \in\qff X$\dss
as\sss the homotopy class of\dss the\qss \emph{index\dss map},\oss
defined as\sss the composition\vspace{-4.5pt}\vspace{-0.125pt}
\[
\quad
\begin{tikzcd}[column sep=boom, row sep=boom]
X
\arrow[r, "\dis s"]
&
X_{\trf U}
\arrow[r, "\dis \mathbb{B}_{\qff U,\qff \bm{\varepsilon}}"]
&
\protect{\num{\mathcal{E}}}
\arrow[r, "\dis \protect{\num{\pi}}"]
&
\protect{\num{\mathcal{S}}}
\qff.
\end{tikzcd}
\]

\vspace{-12pt}
As in\sss the case of\dss self-adjoint\sss operators,\oss
one can\sss go one step further and\sss use\sss the canonical\sss
homotopy equivalence between\sss
$\num{\mathcal{S}}$\sss and\sss $\num{S}$\nnsp.\oss
It\dss leads\sss to a homotopy class\sss of\dss maps\sss
$X\qff \ttoo\qff \num{S}$\nnsp,\oss which can\sss be also considered
as\sss the\qss \emph{analytical\dss index}\pss of\dss the family\sss 
$B_{\dff x}\fff,\qff x\qff \in\qff X$\nnsp.\oss

\myuppar{Topological\sss categories related\sss to\dss 
Fredholm\dss operators\sss in\dss Hilbert\dss bundles.}
As in\dss Section\qss \ref{analytic-index-section},\oss
let\sss $H_{\dff x}\dff,\qff x\qff \in\qff X$\nnsp,\oss
denoted also by\sss $\mathbb{H}$\nnsp,\oss 
be a\sss locally\sss trivial\dss Hilbert\dss bundle with separable fibers over $X$\nnsp.\oss
In\sss this section we are interested\sss in\sss families of\trs
Fredholm\dss operators\sss 
$B_{\dff x}\dff \colon\dff H_{\dff x}\qff \ttoo\qff H_{\dff x}$\sss
parameterized\sss by\sss $x\qff \in\qff X$\sss
and\sss not\sss assumed\sss to be self-adjoint.\oss

Let\sss $\mathcal{S}\dff(\trf \mathbb{H}\trf)$\sss be\sss the following category.\oss
The objects of\dss $\mathcal{S}\dff(\trf \mathbb{H}\trf)$\sss are pairs\sss
$(\trf E_{\dff 1}\dff,\pff E_{\trf 2}\trf)$\sss
such\sss that\sss $E_{\dff 1}\dff,\pff E_{\trf 2}$\sss are
finitely dimensional\sss subspaces of\dss the same fiber\sss $H_{\dff x}$\sss
of\dss $\mathbb{H}$\nnsp.\oss
The\sss topology on\sss the set\sss of\dss objects\dss is\dss defined\sss
in\sss the same way as for\sss $\hat{\mathcal{S}}\dff(\trf \mathbb{H}\trf)$\nnsp.\oss
Morphisms\sss
$(\trf E_{\dff 1}\fff,\qff E_{\dff 2}\trf)
\qff \ttoo\qff 
(\trf E\fff'_{\dff 1}\fff,\qff E\fff'_{\dff 2}\trf)$\sss
exists only\sss if\dss the four\sss involved\sss subspaces
are subspaces of\dss the same fiber\sss $H_{\dff x}$\nsp,\oss
and\sss in\sss this case morphisms are defined exactly as morphisms of\sss
$\mathcal{S}$\sss with\sss $H$\sss replaced\sss by\sss $H_{\dff x}$\nsp.\oss
If\dss we\sss turn\sss $X$ into\sss topological\sss category\sss 
having\sss only\sss the identity\sss morphisms,\oss
then\sss there\dss is\dss a\sss functor\sss
$\eta\dff \colon\dff 
\mathcal{S}\dff(\trf \mathbb{H}\trf)
\qff \ttoo\qff
X$\sss
assigning\sss to an object\sss $(\trf E_{\dff 1}\dff,\pff E_{\trf 2}\trf)$\sss 
the point\sss $x\qff \in\qff X$\sss
such\sss that\sss $E_{\dff 1}\dff,\pff E_{\trf 2}\qff \subset\qff H_{\dff x}$\nsp.\oss

There\dss is\dss also an analogue\sss 
$S\dff(\trf \mathbb{H}\trf)$\sss 
of\sss $S$\nnsp.\oss
The objects of\sss $S\dff(\trf \mathbb{H}\trf)$\sss are\sss triples\sss 
$(\trf x\fff,\qff V_{\dff 1}\dff,\pff V_{\trf 2}\trf)$\sss
such\sss that\sss $x\qff \in\qff X$\sss and\sss 
$(\trf V_{\dff 1}\dff,\pff V_{\trf 2}\trf)$\sss 
is\dss an object\sss of\dss $S$\nnsp.\oss
The\sss topology on\sss the set\sss of\dss such\sss triples\dss is\dss defined\sss by\sss
the\sss topology of\sss $X$\sss and\sss the discrete\sss topology on\sss
the set\sss of\dss spaces.\oss
Morphisms\sss
$(\trf x\fff,\qff V_{\dff 1}\dff,\pff V_{\trf 2}\trf)
\qff \ttoo\qff 
(\trf x\fff,\qff V\fff'_{\fff 1}\dff,\pff V\fff'_{\dff 2}\trf)$\sss
exist\sss only\sss if\sss $x\off =\off x'$\nnsp,\oss
and\sss in\sss this case\sss they are\sss the same as morphisms\sss
$(\trf V_{\dff 1}\dff,\pff V_{\trf 2}\trf)
\qff \ttoo\qff 
(\trf V\fff'_{\fff 1}\dff,\pff V\fff'_{\dff 2}\trf)$\sss in\sss the category $S$\nnsp.\oss
The\sss topology on\sss the set\sss of\dss morphisms\dss is\dss defined\sss by\sss
the\sss topology of\sss $X$\sss and\sss the\sss topology on\sss the set\sss
of\dss morphisms of\sss $S$\nnsp.\oss
Clearly,\qss $S\dff(\trf \mathbb{H}\trf)$\sss  
depends only\sss on $X$\nnsp.\oss
There are functors\sss
$\eta\dff \colon\dff 
S\dff(\trf \mathbb{H}\trf)
\qff \ttoo\qff
X$\sss
and\sss
$S\dff(\trf \mathbb{H}\trf)
\qff \ttoo\qff
S$\sss
leading\sss to a homeomorphism\sss
$\num{S\dff(\trf \mathbb{H}\trf)}
\qff \ttoo\qff
X\dff \times\dff \num{S}$\sss
for compactly\sss generated\sss $X$\nnsp.\oss
There\dss is\dss also an\sss intermediate category\sss $S/\fff \mathbb{H}$\sss
and\sss functors\sss
$\mathcal{S}\dff(\trf \mathbb{H}\trf)
\off \longleftarrow\off
S/\fff \mathbb{H}
\qff \ttoo\qff
S\dff(\trf \mathbb{H}\trf)$\sss
such\sss that\sss the geometric realizations\sss
$\num{\mathcal{S}\dff(\trf \mathbb{H}\trf)}
\off \longleftarrow\off
\num{S/\fff \mathbb{H}}
\qff \ttoo\qff
\num{S\dff(\trf \mathbb{H}\trf)}$\sss
are homotopy equivalences.\oss

\myuppar{Fredholm\dss families of\dss operators in\dss Hilbert\dss bundles.}
The definitions of\qss \emph{Fredholm\dss families}\pss of\dss operators
and of\qss \emph{Fredholm\dss homotopies}\pss apply,\oss
with obvious modifications,\oss to families\sss of\dss operators\sss
$B_{\dff x}\dff \colon\dff H_{\dff x}\qff \ttoo\qff H_{\dff x}$\nsp,\qss
$x\qff \in\qff X$\nnsp.\oss
But,\oss as in\sss the self-adjoint\sss case,\oss 
the definition of\dss the analytical\dss index should\sss be modified.\oss
Let\sss $B_{\dff x}\fff,\qff x\qff \in\qff X$\dss be a family\sss of\trs
Fredholm\sss operators,\oss and\sss let\sss
$U_{\dff a}\dff,\pff a\qff \in\qff \Sigma$\sss of\dss $X$\sss
and\sss $\varepsilon_{\dff a}\dff,\pff a\qff \in\qff \Sigma$\sss
be as above.\oss
Then one can define a continuous functor\vspace{1.5pt}\vspace{0.25pt}
\[
\quad
\mathbb{B}_{\qff U,\qff \bm{\varepsilon}}\dff \colon\dff
X_{\dff U}\qff \ttoo\qff \mathcal{S}\dff(\trf \mathbb{H}\trf)
\pff
\]

\vspace{-12pt}\vspace{1.5pt}\vspace{0.25pt}
as follows.\oss
Let
$(\trf z\fff,\qff \tau\trf)$ 
be object\sss of\dss $X_{\dff U}$
and\sss $\varepsilon\off =\off \varepsilon_{\dff \tau}$.\oss
The\sss functor $\mathbb{B}_{\qff U,\qff \bm{\varepsilon}}$\sss
takes $(\trf z\fff,\qff \tau\trf)$
to\vspace{3pt}
\[
\quad
(\trf E_{\dff 1}\fff,\qff E_{\dff 2}\trf)
\off =\off
\left(\qff
\image\dff P_{\qff [\dff 0\fff,\qff \varepsilon\dff]}\trf(\trf \num{B}\trf)\fff,\off
\image\dff P_{\qff [\dff 0\fff,\qff \varepsilon\dff]}\trf(\trf \num{B^*}\trf)
\qff\right)
\qff,
\]

\vspace{-12pt}\vspace{3pt}
and\sss this rule naturally extends\sss to morphisms.\oss
The functor\sss
$\mathbb{B}_{\qff U,\qff \bm{\varepsilon}}$\sss
leads\sss to a continuous map\vspace{1.5pt}
\[
\quad
\num{\mathbb{B}_{\qff U,\qff \bm{\varepsilon}}}\qff \colon\qff
\num{X_{\dff U}}
\qff \ttoo\qff 
\num{\mathcal{S}\dff(\trf \mathbb{H}\trf)}
\pff.
\]

\vspace{-12pt}\vspace{1.5pt}
\mypar{Theorem.}{bundle-correct-non-sa}
\emph{Suppose\sss that\sss $X$\sss is\dss paracompact\sss
and\sss $s\dff \colon\dff X\qff \ttoo\qff \num{X_{\dff U}}$\sss
is\dss a homotopy\sss inverse of\pss $\num{\pr}$
as\sss in\qss Theorem\qss \ref{covering-categories}.\oss
Then\sss the homotopy class of\pss 
$\num{\mathbb{B}_{\qff U,\qff \bm{\varepsilon}}}\qff \circ\qff
s
\trf \colon\dff
X
\qff \ttoo\qff 
\num{\mathcal{S}\dff(\trf \mathbb{H}\trf)}$\sss
does not\sss depend on\sss the choices\sss involved.\oss}

\proof
The proof\trs is\dss completely similar\sss to\sss the proof\dss of\trs
Theorem\qss \ref{lift-fixed}.\oss  \eproof

\myuppar{Analytical\dss index\sss of\dss families 
of\dss operators\sss in\dss Hilbert\dss bundles.}
Suppose\sss that\sss $X$\sss is\dss paracompact.\oss
Then\sss any\sss map\sss of\dss the form\dss
$\num{\mathbb{A}_{\qff U,\qff \bm{\varepsilon}}}\qff \circ\qff
s
\trf \colon\dff
X
\qff \ttoo\qff 
\num{\mathcal{S}\dff(\trf \mathbb{H}\trf)}$\sss
may\sss be called an\qss \emph{index\sss map}\pss
of\dss the family\sss
$B_{\dff x}\dff,\qff
x\qff \in\qff X$\nnsp.\oss
By\trs Theorem\qss \ref{bundle-correct-non-sa}\qss the homotopy class
of\dss an index\sss map\dss is\dss independent\sss from choices 
involved\sss in\sss its construction,\oss
and\sss we can define\sss the\qss \emph{index}\pss of\dss the family\sss
$B_{\dff x}\dff,\qff
x\qff \in\qff X$\dss
as\sss the homotopy class of\dss any of\dss its\sss index\sss maps.\oss
One can also use\sss the
homotopy equivalence between\sss
$\num{\mathcal{S}\dff(\trf \mathbb{H}\trf)}$\sss 
and\sss 
$\num{S\dff(\trf \mathbb{H}\trf)}$\sss
and\sss get\sss a map\sss
$X\qff \ttoo\qff
\num{S\dff(\trf \mathbb{H}\trf)}
\qff \ttoo\qff
X\dff \times\dff \num{S}$\sss
well\sss defined up\sss to homotopy.\oss
The first\sss component\sss of\dss this map\dss is\dss
homotopic\sss to\sss the identity\sss $X\qff \ttoo\qff X$\nnsp,\oss
and one can define\qss \emph{index}\pss as\sss the homotopy class
of\dss the second component\sss
$X\qff \ttoo\qff \num{S}$\nnsp.\oss
Also,\oss as in\sss the self-adjoint\sss case,\oss one can use\sss the fact\sss that\dss 
Hilbert\dss bundles are\sss trivial\sss and define\qss \emph{index}\pss
as\sss the second\sss component\sss
$X\qff \ttoo\qff \num{\mathcal{S}}$\sss 
of\dss the homotopy class\sss\vspace{1.5pt}
\[
\quad
X
\qff \ttoo\qff 
\num{\mathcal{S}\dff(\trf \mathbb{H}\trf)}
\qff \ttoo\qff 
X\dff \times\dff \num{\mathcal{S}}
\pff.
\]

\mysection{Hilbert\qss bundles\qss and\qss strictly\qss Fredholm\qss families}{strictly-fredholm-families}

\myuppar{Compactly\sss generated spaces and\sss the compact-open\sss topology.}
Recall\sss that\sss we\sss tacitly\sss assume all\sss topological\sss spaces\sss to be
compactly\sss generated.\oss
In\sss this section\sss this assumption will\sss play a more prominent\sss role\sss
than\sss usual,\oss and\sss by\sss this reason\sss it\sss will\sss be sometimes stated explicitly.\oss
For\sss two\sss topological\sss spaces\sss $Y\fff,\qff Z$\sss
we will\sss denote by\sss $\map\trf(\trf Y\fff,\qff Z\trf)$\sss the space
of\dss continuous maps\sss $Y\qff \ttoo\qff Z$\sss
with\sss the\qss \emph{compact-open\sss topology}.\oss
For\sss a compactly\sss generated space $X$\sss
the continuity of\dss a map\sss
$X\qff \ttoo\qff \map\trf(\trf Y\fff,\qff Z\trf)$\sss
is\dss equivalent\sss to\sss the continuity of\dss the corresponding\qss
(adjoint{\fff})\qss map\sss
$X\dff \times\dff Y\qff \ttoo\qff Z$\nnsp.\oss
This\dss is\dss one of\dss the main\sss reasons for working\sss with\sss
compactly\sss generated spaces and\sss the compact-open\sss
topology\sss in\sss the context\sss of\dss Hilbert\sss bundles.\oss

\mypar{Lemma.}{composition}
\emph{Suppose\sss that\dss $X$\sss is\dss a compactly\sss generated\sss space.\oss
Let}\sss\vspace{1.5pt} 
\[
\quad
f\dff \colon\dff X\qff \ttoo\qff \map\trf(\trf Y\fff,\qff Z\trf)
\quad
\mbox{\emph{and}}\dff\quad
g\dff \colon\dff X\qff \ttoo\qff \map\trf(\trf Z\fff,\qff W\trf)
\]

\vspace{-12pt}\vspace{1.5pt}
\emph{be\sss continuous maps.\oss
Then\sss the map}\sss\vspace{1.5pt}
\[
\quad
h\dff \colon\dff
X
\qff \ttoo\qff 
\map\trf(\trf Y\fff,\qff Z\trf)
\]

\vspace{-12pt}\vspace{1.5pt}
\emph{defined\dss by\sss 
$h\dff(\dff x\trf)
\off =\off 
g\dff(\dff x\trf)\dff \circ\dff f\dff(\dff x\trf)$\sss
is\dss continuous.\oss}

\proof
The continuity of\sss $f$\sss
is\dss equivalent\sss to\sss the continuity of\dss the map\sss
$X\dff \times\dff Y\qff \ttoo\qff X\dff \times\dff Z$\sss
defined\sss by\sss
$(\dff x\fff,\qff y\trf)
\off \longmapsto\off
\left(\dff x\fff,\qff f\dff(\dff x\trf)\dff (\trf y\trf)\trf\right)$\nnsp,\oss
and\sss similarly\sss for\sss $g$ and\sss $h$\nnsp.\oss
Therefore\sss the\sss continuity\sss of\dss the
composition\sss
$X\dff \times\dff Y\qff \ttoo\qff X\dff \times\dff Z\qff \ttoo\qff X\dff \times\dff W$\sss
implies\sss the continuity of\sss $h$\nnsp.\oss  \eproof

\myuppar{Local\sss trivializations of\qss Hilbert\dss bundles.}
Suppose\sss that\sss $\mathbb{K}$\sss is\dss a\sss locally\sss trivial\dss
Hilbert\dss bundle over a\sss compactly\sss generated\sss topological\sss space\sss $Y$\sss
with\sss fibers isomorphic\sss to a separable infinitely dimensional\dss 
Hilbert\dss space $K$\nnsp.\oss
As before,\oss we will\dss treat\sss $\mathbb{K}$\sss
as a family\sss
$K_{\dff y}\dff,\qff y\qff \in\qff Y$ of\dss Hilbert\sss spaces parameterized\sss by $X$\nnsp.\oss
A\qss \emph{local\dss trivialization}\pss of\dss $\mathbb{K}$\sss can\sss be considered as a 
family of\dss Hilbert\sss space isomorphisms\sss
$t_{\trf U}\dff(\dff x\trf)\dff \colon\dff K_{\dff x}\qff \ttoo\qff K$\nnsp,\oss 
where $x$ runs over an open subset\qss $U\qff \subset\qff Y$\nnsp.\oss
Two\sss trivialization\dss $t_{\trf U}$\sss and\sss $t_{\qff V}$\sss 
are\sss related\sss by continuous
\emph{transition\sss map}\dss
\[
\quad
(\trf U\dff \cap\dff V\trf)\dff \times\dff K
\qff \ttoo\qff
(\trf U\dff \cap\dff V\trf)\dff \times\dff K
\]

\vspace{-12pt}
given\sss by\sss the formula
\[
\quad
(\dff x\fff,\qff v\trf)
\off \longmapsto\off
\left(\trf
x\dff,\qff
t_{\trf U}\dff(\dff x\trf)
\dff \circ\dff
t_{\qff V}\dff(\dff x\trf)^{\dff -\dff 1}\dff (\dff v\trf)
\trf\right)
\pff.
\]

\vspace{-12pt}
Since\sss $Y$\sss is\dss assumed\sss to be compactly\sss generated,\oss
the continuity\sss of\dss the above\sss transition map\sss is\dss equivalent\sss to\sss
the continuity of\dss the\qss \emph{transition\dss function}\vspace{1.25pt}
\[
\quad
s_{\trf U\fff V}\dff \colon\dff
x
\off \longmapsto\off
t_{\trf U}\dff(\dff x\trf)
\dff \circ\dff
t_{\qff V}\dff(\dff x\trf)^{\dff -\dff 1}
\]

\vspace{-12pt}
from\sss $U\dff \cap\trf V$\sss to\sss the group of\dss 
isometries\sss $K\qff \ttoo\qff K$\sss
equipped\sss with\sss the com\-pact-open\dss topology.\oss
Since\sss
$s_{\qff V\fff U}\off =\off s_{\trf U\fff V}^{\dff -\dff 1}$\sss
is\dss also continuous,\oss
the maps\sss $s_{\trf U\fff V}$\sss will\sss be continuous even\sss if\dss
the group of\dss isometries of\sss $K$\sss is\dss equipped\sss with\sss a stronger\sss topology,\oss
namely\sss with\sss the\sss topology\sss
induced\sss from\sss the product\sss of\dss the compact-open\sss topologies
by\sss the map\vspace{1.5pt}
\[
\quad
g
\off \longmapsto\off 
\left(\trf g\fff,\qff g^{\dff -\dff 1}\trf\right)
\qff.
\]

\vspace{-12pt}\vspace{1.5pt}
We will\sss denote by $\mathcal{U}\dff(\trf K\trf)$\sss the groups of\dss
isometries of $K$ with\sss the\sss latter\sss topology.\oss
Actually,\oss this\sss topology coincides with\sss the strong operator\sss topology,\oss
but\sss the author prefers\sss to ignore\sss this fact.\oss
We\sss will\sss reserve\sss the notation $U\dff(\trf K\trf)$\sss
for\sss the same group with\sss the norm\sss topology.\oss

\myuppar{Trivializations\sss of\trs Hilbert\dss bundles.}
Let\sss $\mathbb{K}$\sss be a\dss Hilbert\dss bundle as above.\oss 
A\qss \emph{trivialization}\pss of\dss the bundle\sss $\mathbb{K}$\sss is\dss
an\sss isomorphism\sss between $\mathbb{K}$\sss and\sss the\sss trivial\sss
bundle\sss $Y\dff \times\dff K\qff \ttoo\qff Y$\nnsp.\oss
A\sss trivialization\sss $t$\sss can\sss be\sss thought\sss 
as a family\sss of\dss isometries\sss
$t_{\trf y}\dff \colon\dff K_{\dff y}\qff \ttoo\qff K$,\trs
$y\qff \in\qff Y$\dss 
defining a homeomorphism
between\sss the\sss total\sss space\dss
$\bigcup_{\qff y\qff \in\qff Y}\qff K_{\dff y}$\sss
of\dss the bundle\sss $\mathbb{K}$\sss and\sss $Y\dff \times\dff K$\nnsp.\oss
Two\sss trivializations\sss $t\fff,\qff u$\sss 
differ\sss by a\sss family\sss\vspace{1.5pt}
\[
\quad
t_{\dff y}\dff \circ\dff u_{\dff y}^{\dff -\dff 1}
\dff \colon\dff
K\qff \ttoo\qff K
\] 

\vspace{-12pt}\vspace{1.5pt}
of\dss isometries.\oss
Somewhat\sss disappointingly,\oss the map\sss
$s\dff \colon\dff
y\off \longmapsto\off t_{\dff y}\dff \circ\dff u_{\dff y}^{\dff -\dff 1}$\sss
from\sss $Y$\sss to\sss the unitary\sss group of\sss $K$\sss
is\dss continuous only\sss in a\sss relatively\sss weak sense.\oss
Assuming\sss that\sss $Y$\sss is\dss compactly generated,\pss
$s$\sss is\dss continuous as a map\sss to\sss $\mathcal{U}\dff(\trf K\trf)$\nnsp,\oss
but\sss usually\sss not\sss as a map\sss to\sss $U\dff(\trf K\trf)$\nnsp.\oss

In fact,\oss Hilbert\dss bundles with\sss fibers\sss isomorphic\sss to $K$\sss
over compactly\sss generated spaces
are\sss bundles with\sss the structure group\sss $\mathcal{U}\dff(\trf K\trf)$\nnsp.\oss
The group\sss $\mathcal{U}\dff(\trf K\trf)$\sss is\dss known\sss to be contractible.\oss
This result\dss is\dss essentially due\sss to\sss to\dss Dixmier\sss and\dss Douady\qss \cite{dd},\oss
and was put\sss in\sss the present\sss form\sss by\dss Atiyah\dss and\trs Segal\qss \cite{ase}.\oss
The contractibility of\dss $\mathcal{U}\dff(\trf K\trf)$\sss implies\sss
that\dss Hilbert\dss bundles over\sss 
paracompact\sss compactly\sss generated spaces are\sss trivial.\oss
Moreover,\oss given a\sss bundle with\sss the base $Y$\sss and\sss a
closed subset\sss $Z\qff \subset\qff Y$\nnsp,\oss 
a\sss trivialization over\sss $Z$\sss can\sss be extended\sss to a\sss
trivialization over\sss $Y$\dnsp.\oss
This implies\sss
that\sss any\sss two\sss trivializations are homotopic\sss in an obvious sense.\oss

\myuppar{Adapted\sss pairs and\sss local\dss trivializations.}
As in\dss Section\qss \ref{analytic-index-section},\oss
let\sss $\mathbb{H}$\sss be a\sss 
locally\sss trivial\dss Hilbert\dss bundle with separable fibers over $X$\nnsp,\oss
thought\sss as a family\sss
$H_{\dff x}\dff,\qff x\qff \in\qff X$ of\dss Hilbert\sss spaces parameterized\sss by $X$\nnsp.\oss
Let\sss 
$A_{\dff x}\dff \colon\dff H_{\dff x}\qff \ttoo\qff H_{\dff x}$,\trs
$x\qff \in\qff X$\dss
be a\dss family of\dss self-adjoint\trs Fredholm\dss operators,\oss
also denoted\sss by\sss $\mathbb{A}$\nnsp.\oss
Let\sss us\sss slightly\sss rephrase\sss the definition of\trs
Fredholm\dss families.\oss

Suppose\sss that\sss $U\qff \subset\qff X$\sss 
and\sss $\varepsilon\qff >\qff 0$\nnsp.\oss
We will\sss say\sss that\sss the pair $(\trf U\fff,\qff \varepsilon\trf)$\sss
is\qss \emph{adapted\sss to\sss the family\sss $\mathbb{A}$\nnsp,}\oss
or\sss simply\qss \emph{adapted},\oss if\trs
for\sss every\sss $x\qff \in\qff U$\dss
the\sss image\dss 
$\image
P_{\qff [\dff -\dff \varepsilon\fff,\dff \varepsilon\dff]}\dff(\dff A_{\dff x}\dff)$\sss
is\trs finitely\sss dimensional,\pss
${}-\qff \varepsilon\dff,\qff \varepsilon
\pff \not\in\pff \sigma\dff(\dff A_{\dff x}\dff)$\nnsp,\oss
the family\sss of\dss the images\dss\vspace{1.5pt}\vspace{-0.575pt}
\[
\quad
V_{\dff x}
\off =\off
\image
P_{\qff [\dff -\dff \varepsilon\fff,\dff \varepsilon\dff]}\dff(\dff A_{\dff x}\dff)\fff,\qff
x\qff \in\qff U
\]

\vspace{-12pt}\vspace{1.5pt}\vspace{-0.575pt}
is\dss continuous,\oss
as also\sss the family of\dss operators\sss
$V_{\dff x}\qff \ttoo\qff V_{\dff x}$\sss
induced\sss by\sss $A_{\dff x}$\nsp.\oss 
Clearly,\pss $\mathbb{A}$\sss is\dss a\dss Fredholm\dss family\trs
if\trs and\dss only\trs if\trs there exists a\sss family\sss of\dss adapted\sss pairs\sss
$(\trf U\fff,\qff \varepsilon\trf)$\sss such\sss that\sss the sets $U$ are open and
cover $X$\nnsp.\oss
We will\sss call\sss such a\sss family of\dss pairs\sss
$(\trf U\fff,\qff \varepsilon\trf)$\sss
an\qss \emph{atlas\sss for}\sss $\mathbb{A}$\nnsp.\oss

Suppose\sss that\sss the pair\sss $(\trf U\fff,\qff \varepsilon\trf)$\sss
is\dss adapted\dss to\sss $\mathbb{A}$\nnsp.\oss
A\sss local\dss trivialization\sss $t_{\trf U}$ of\dss $\mathbb{H}$\sss defined over $U$\sss
is\dss said\sss to be\qss \emph{strictly\sss adapted\dss to}\dss $(\trf U\fff,\qff \varepsilon\trf)$\trs
({\fff}and\sss the family\sss $\mathbb{A}$\nsp)\pss
if\trs the family of\dss projections\vspace{1.5pt}\vspace{-0.95pt}
\[
\quad
t_{\trf U}\dff(\dff x\trf)
\qff \circ\qff
P_{\qff \geq\qff \varepsilon}\qff
\bigl(\dff A_{\dff x}\dff\bigr)
\qff \circ\qff
t_{\trf U}\dff(\dff x\trf)^{\dff -\dff 1}
\off =\off
P_{\qff \geq\qff \varepsilon}\qff
\left(\trf 
t_{\trf U}\dff(\dff x\trf)
\qff \circ\qff
A_{\dff x}
\qff \circ\qff
t_{\trf U}\dff(\dff x\trf)^{\dff -\dff 1}
\trf\right)
\]

\vspace{-12pt}\vspace{1.5pt}\vspace{-0.95pt}
is\dss norm continuous.\oss
Equivalently,\oss the family of\dss subspaces\dss\vspace{1.5pt}\vspace{-0.95pt}
\[
\quad
H_{\qff \geq\qff \varepsilon}^{\dff U}\dff(\dff x\trf)
\off =\off
t_{\trf U}\dff(\dff x\trf)\qff
\bigl(\trf
\image
P_{\qff \geq\qff \varepsilon}\trf(\dff A_{\dff x}\dff)
\trf\bigr)
\]

\vspace{-12pt}\vspace{1.5pt}\vspace{-0.95pt}
is\dss norm continuous.\oss
The next\sss lemma shows\sss that\sss the property of\dss being adapted\sss
to\sss $(\trf U\fff,\qff \varepsilon\trf)$\sss is\dss independent\sss of\dss
the choice of\sss $\varepsilon\qff >\qff 0$\sss such\sss that\sss
$(\trf U\fff,\qff \varepsilon\trf)$\sss is\dss adapted.\oss
We will\sss say\sss that\sss a\sss local\dss trivialization\dss is\qss
\emph{strictly\sss adapted\dss to}\trs $U$\sss if\dss it\dss is\dss
strictly adapted\sss to\sss $(\trf U\fff,\qff \varepsilon\trf)$\sss
for some\sss $\varepsilon\qff >\qff 0$\nnsp.\oss

\mypar{Lemma.}{control-independence}
\emph{Let\sss $U\qff \subset\qff X$\sss
and\sss $(\trf U\fff,\qff \varepsilon\trf)$\nnsp,\qss
$(\trf U\fff,\qff \delta\trf)$\sss are pairs adapted\sss to\sss $\mathbb{A}$\nnsp.\oss
If\dss a\sss local\dss trivialization\sss $t_{\trf U}$\sss defined over $U$\sss
is\dss strictly\sss adapted\sss to\sss $(\trf U\fff,\qff \varepsilon\trf)$\nnsp,\oss
then\sss $t_{\trf U}$\sss is\dss strictly\sss adapted\dss to $(\trf U\fff,\qff \delta\trf)$\nnsp.}

\proof
Suppose\sss that\sss $\delta\qff <\qff \varepsilon$\nnsp.\oss
Then\oss\vspace{1.5pt}\vspace{-0.95pt}
\[
\quad
\image
P_{\qff \geq\qff \varepsilon}\trf(\dff A_{\dff x}\dff)
\off =\off
\image
P_{\qff \geq\qff \delta}\trf(\dff A_{\dff x}\dff)
\qff \oplus\qff
\image
P_{\qff [\dff \delta\fff,\dff \varepsilon\dff]}\dff(\dff A_{\dff x}\dff)
\]

\vspace{-12pt}\vspace{1.5pt}\vspace{-0.95pt}
for\sss every\sss $x\qff \in\qff U$\nnsp.\oss
Since\sss $\varepsilon\fff,\qff \delta\pff \not\in\pff \sigma\dff(\dff A_{\dff x}\dff)$\nnsp,\oss
the family of\dss subspaces\sss\vspace{1.5pt}\vspace{-0.95pt}
\[
\quad
t_{\trf U}\dff(\dff x\trf)\qff
\bigl(\trf
\image
P_{\qff [\dff \delta\fff,\dff \varepsilon\dff]}\dff(\dff A_{\dff x}\dff)
\trf\bigr)\qff,\off\off
x\qff \in\qff U
\]

\vspace{-12pt}\vspace{1.5pt}\vspace{-0.95pt}
is\dss norm continuous.\oss
Hence\sss the continuity of\dss the family\sss 
$H_{\qff \geq\qff \varepsilon}^{\dff U}\dff(\dff x\trf)$\sss
is\dss equivalent\sss to\sss the continuity of\dss the family\sss
$H_{\qff \geq\qff \delta}^{\dff U}\dff(\dff x\trf)$\nnsp.\oss
The case of\dss $\varepsilon\qff <\qff \delta$\sss is\dss similar.\oss
This proves\sss the\sss lemma.\oss  \eproof

\mypar{Corollary.}{adapted-to-two}
\emph{Let\trs
$(\trf U\fff,\qff \varepsilon\trf)${\nsp}
and\dss
$(\trf V\fff,\qff \delta\trf)$
be\sss adapted\sss pairs\sss
and\sss $t_{\trf U}$\sss and\sss $t_{\qff V}$\sss
be\sss local\sss trivializations strictly\sss adapted\sss to\dss 
$(\trf U\fff,\qff \varepsilon\trf)$
and\trs
$(\trf V\fff,\qff \delta\trf)$
respectively.\oss
Then\sss the pair\sss $(\trf U\dff \cap\dff V\fff,\qff \varepsilon\trf)$\sss
is\dss adapted and\sss the restrictions of\qss $t_{\trf U}$\sss and\sss $t_{\qff V}$\sss to\sss
$U\dff \cap\dff V$\sss are strictly\sss adapted\sss to
$(\trf U\dff \cap\dff V\fff,\qff \varepsilon\trf)$\nnsp.\oss}

\proof
The first\sss claim\dss is\dss trivial,\oss
and\dss the second one follows\sss from\trs Lemma\qss \ref{control-independence}.\oss  \eproof

\myuppar{Strictly\dss Fredholm\dss families and adapted\dss trivializations.}
We will\sss say\sss that\sss $\mathbb{A}$\sss is\dss a\qss
\emph{strictly\trs Fredholm\dss family}\pss if\trs
there exists an atlas\sss for\sss $\mathbb{A}$\sss
such\sss that\sss for each\sss pair\sss
$(\trf U\fff,\qff \varepsilon\trf)$\sss
from\sss this atlas\sss there\dss exists\dss
a\sss local\dss trivialization\sss
$t_{\trf U}$\sss strictly\sss adapted\sss to\sss
$(\trf U\fff,\qff \varepsilon\trf)$\nnsp,\oss
or simply\sss to $U$\nnsp,\oss
and\sss the family\sss 
$\mathbb{A}$\nnsp.\oss
We will\sss say\sss that\sss a\sss trivialization\sss 
$t$\sss of\dss the bundle\sss $\mathbb{H}$\sss
is\qss \emph{strictly\sss adapted\dss to\sss the family\dss
$\mathbb{A}$}\qss
if\trs there exists an atlas for\sss $\mathbb{A}$\sss
such\sss that\sss for every\sss
pair\sss $(\trf U\fff,\qff \varepsilon\trf)$\sss
from\sss this atlas\sss
the restriction of\sss $t$\sss to $U$\sss is\dss strictly\sss adapted\sss to\sss
$(\trf U\fff,\qff \varepsilon\trf)$\sss
and\sss $\mathbb{A}$\nnsp.\oss
Clearly,\oss if\dss there exists a\sss trivialization strictly\sss adapted\dss to $\mathbb{A}$\nnsp,\oss
then\sss $\mathbb{A}$\sss is\dss a strictly\dss Fredholm\dss family.\oss

\mypar{Lemma.}{atlas-independence}
\emph{If\qss $t$\dss is\dss a\sss trivialization strictly\sss adapted\sss to\sss $\mathbb{A}$\sss
and\sss $(\trf U\fff,\qff \varepsilon\trf)$\sss is\dss a pair adapted\sss to\sss $\mathbb{A}$\nnsp,\oss
then\sss the restriction of\qss $t$\dss to\sss $U$\sss is\dss 
strictly\sss adapted\dss to\sss $\mathbb{A}$\nnsp.\oss}

\proof
Since\sss $t$\sss is\dss strictly\sss adapted\sss to $\mathbb{A}$\nnsp,\oss
for every\sss $x\qff \in\qff X$\sss there exists an adapted\sss pair\sss
$(\trf V\fff,\qff \delta\trf)$ such\sss that\sss $V$\sss is\dss open,\pss 
$x\qff \in\qff V$\dnsp,\oss
and\sss the restriction\sss $t_{\qff V}$\sss of\dss $t$\sss to $V$\sss 
is\dss strictly\sss adapted\sss to $(\trf V\fff,\qff \delta\trf)$\nnsp.\oss
By\trs Corollary\qss \ref{adapted-to-two}\qss the restriction of\sss $t$\sss
to\sss $U\dff \cap\dff V$\sss is\dss strictly\sss adapted\sss to\sss
$(\trf U\dff \cap\dff V\fff,\qff \varepsilon\trf)$\nnsp.\oss
By\sss the definition,\oss the property of\dss being strictly\sss adapted\sss to\sss
$(\trf U\fff,\qff \varepsilon\trf)$\sss is\dss a continuity\sss property.\oss
Therefore\sss it\sss holds\sss if\dss it\dss holds in a neighborhood of\dss each\sss point\sss
$x\qff \in\qff V$\dnsp.\oss
The\sss lemma\sss follows.\oss  \eproof

\mypar{Theorem.}{adapted}
\emph{If\pss $\mathbb{A}$\sss
is\dss a strictly\trs Fredholm\dss family\sss
and\dss $X$ is\dss a\sss triangulated space,\oss
then\sss there exists a\sss trivialization of\trs the\sss bundle\qss
$\mathbb{H}$\dss 
strictly\sss adapted\dss to\sss $\mathbb{A}$\nnsp.\oss}\vspace{-0.125pt}

\proof
Let\sss us\sss choose an atlas of\dss pairs\sss
$(\trf U\fff,\qff \varepsilon\trf)$\sss
adapted\sss to\sss $\mathbb{A}$\sss and\dss
and\dss local\dss trivializations\sss $t_{\trf U}$\sss
strictly\sss adapted\dss to\sss $(\trf U\fff,\qff \varepsilon\trf)$\sss
and\sss $\mathbb{A}$\nnsp.\oss
Replacing\sss the given\sss triangulation of\sss $X$\sss
by a subdivision,\oss if\dss necessary,\oss
we can assume\sss that\sss every\sss simplex $\sigma$ of\dss
this\sss triangulation\dss is\dss contained\sss in\sss
the set\sss $U$\sss for some adapted\sss pair\sss 
$(\trf U\fff,\qff \varepsilon\trf)$\sss from our atlas.\oss

Let\sss us\sss begin\sss by\sss choosing isomorphisms\sss
$t\trf(\dff v\trf)\dff \colon\dff H_{\dff v}\qff \ttoo\qff H$\sss
for\sss the vertices $v$ of\dss the\sss triangulation.\oss
Such\sss isomorphisms define a\sss trivialization of\dss the restriction
of\dss the bundle\sss to\sss the $0$\dnsp-skeleton $\ske_{\trf 0}\dff X$,\oss
which\dss is\dss obviously strictly\sss adapted\sss to\sss the restriction of\dss the family\sss
$A_{\dff x}\dff,\qff x\qff \in\qff X$\sss to\sss 
$\ske_{\trf 0}\dff X$\nnsp,\oss
i.e.\qss to\sss the family\sss
$A_{\dff x}\dff,\qff x\qff \in\qff \ske_{\trf 0}\dff X$\nnsp.\oss

Suppose\sss that\sss we already\sss found a\sss trivialization\sss $t$\sss
of\dss the restriction of\dss the bundle\sss $\mathbb{H}$\sss to\sss the 
$n$\dnsp-skeleton $\ske_{\trf n}\dff X$\sss
which\dss is\dss strictly\sss adapted\sss to\sss the family\sss
$A_{\dff x}\dff,\qff x\qff \in\qff \ske_{\trf n}\dff X$
for some $n\qff \geq\qff 0$\nnsp.\oss
Let\sss $\sigma$\sss be an $(\dff n\qff +\qff 1\dff)$\dnsp-simplex of\dss
the\sss triangulation,\oss and\sss let\sss $\partial\dff \sigma$ be its boundary.\oss
By\sss our assumptions,\pss $\sigma\qff \subset\qff U$\sss for an adapted\sss
pair\sss $(\trf U\fff,\qff \varepsilon\trf)$\nnsp.\oss
Recall\sss that\sss $t_{\trf U}$\sss is\dss 
strictly\sss adapted\sss to 
$(\trf U\fff,\qff \varepsilon\trf)$\nnsp.\oss
Let\sss $U_{\dff n}\off =\off \ske_{\trf n}\dff X\dff \cap\dff U$\dnsp.\oss
Then\sss $\partial\dff \sigma\qff \subset\qff U_{\dff n}$\nsp,\oss
the pair\sss
$(\trf U_{\dff n}\fff,\qff \varepsilon\trf)$\sss
is\dss adapted\dss to\sss
$A_{\dff x}\dff,\qff x\qff \in\qff \ske_{\trf n}\dff X$\nnsp,\oss 
and\sss the restriction of\sss $t$\sss to\sss
$U_{\dff n}$\sss is\dss strictly\sss adapted\sss to 
$(\qff U_{\dff n}\fff,\qff \varepsilon\trf)$\nnsp.\oss
Let\sss\vspace{1.5pt}
\[
\quad 
H_{\qff <\qff \varepsilon}^{\dff U}\dff(\dff x\trf)
\off =\off
H\qff \ominus\qff H_{\qff \geq\qff \varepsilon}^{\dff U}\dff(\dff x\trf)
\]

\vspace{-12pt}\vspace{1.5pt}
for\sss $x\qff \in\qff U$\nnsp,\oss and\dss
for\sss $x\qff \in\qff U_{\dff n}$\sss let\vspace{1.5pt}
\[
\quad
H_{\qff \geq\qff \varepsilon}\dff(\dff x\trf)
\off =\off
t\trf(\trf x\trf)\qff
\bigl(\trf
\image
P_{\qff \geq\qff \varepsilon}\trf(\dff A_{\dff x}\dff)
\trf\bigr)
\quad
\mbox{and}\quad
H_{\qff <\qff \varepsilon}\dff(\dff x\trf)
\off =\off
H\qff \ominus\qff
H_{\dff\geq\dff \varepsilon}\dff(\dff x\trf)
\qff.
\]

\vspace{-12pt}\vspace{1.5pt}
Recall\dss that\sss a\qss \emph{polarization}\pss of\dss a\dss Hilbert\dss space\sss $K$\dss
is\dss an orthogonal\sss decomposition of\dss the form\sss
$K\off =\off K_{\dff -}\dff \oplus\dff K_{\dff +}$\sss such\sss that\sss both subspaces\sss
$K_{\dff -}$ and\sss $K_{\dff +}$ are\sss infinitely dimensional.\oss
For each\sss $x\qff \in\qff U_{\dff n}$\sss we have\sss two 
polarizations of\dss $H$\nnsp,\oss namely\vspace{1.5pt}
\[
\quad
H
\off =\off
H_{\qff <\qff \varepsilon}^{\dff U}\dff(\dff x\trf)
\qff \oplus\pff
H_{\qff \geq\qff \varepsilon}^{\dff U}\dff(\dff x\trf)
\quad
\mbox{and}\quad
H
\off =\off
H_{\qff <\qff \varepsilon}\dff(\dff x\trf)
\qff \oplus\qff
H_{\qff \geq\qff \varepsilon}\dff(\dff x\trf)
\off.
\] 

\vspace{-12pt}\vspace{1.5pt}
Since\sss these polarization arise\sss from strictly\sss adapted\dss trivializations,\oss
these\sss polarizations continuously depend on $x\qff \in\qff U_{\dff n}$\sss
in\sss the norm\sss topology.\oss
The contractibility of\dss groups\sss $U\dff(\trf K\trf)$\sss with\sss
the norm\sss topology\sss implies\sss that\sss the space of\dss polarizations\dss
is\dss contractible,\oss 
and\sss hence\sss the first\sss family\dss is\dss homotopic\sss
to\sss the second one by a norm continuous homotopy.\oss

The first\sss family\dss is\dss defined also for\sss
$x\qff \in\qff U$\nnsp,\oss
and we can extend\sss this homotopy\sss to $U$\nnsp.\oss
By\sss taking\sss the composition of\dss this homotopy\sss with\sss
$t_{\trf U}$\sss we get,\oss in\sss particular,\oss
a new\sss local\dss trivialization\sss $t\fff'_{\dff U}$\sss
strictly\sss adapted\sss to\sss $(\trf U\fff,\qff \varepsilon\trf)$
and such\sss that\sss for\sss every\sss $x\qff \in\qff U_{\dff n}$\sss
the images of\dss\vspace{2.25pt}
\[
\quad
\image
P_{\qff \geq\qff \varepsilon}\trf(\dff A_{\dff x}\dff)
\] 

\vspace{-12pt}\vspace{2.25pt}
under\sss $t\fff'_{\dff U}\trf(\trf x\trf)$
and\sss
$t\trf(\trf x\trf)$\sss
are equal.\oss
The family of\dss these images\dss is\dss a\dss Hilbert\dss subbundle of\dss
the\sss trivial\sss bundle\sss $U_{\dff n}\dff \times\dff H$\nnsp.\oss
The contractibility of\dss the unitary\sss groups\sss $\mathcal{U}\dff(\trf K\trf)$\sss
implies\sss first\sss that\sss this subbundle\dss is\dss trivial\sss as a\dss
Hilbert\dss bundle,\oss
and\sss then\sss that\sss the family
$t\fff'_{\dff U}\trf(\trf x\trf)\fff,\qff x\qff \in\qff U_{\dff n}$\sss
is\dss homotopic\sss to\sss the family\sss
$t_{\trf U}\trf(\trf x\trf)\fff,\qff x\qff \in\qff U_{\dff n}$\sss
in\sss the class of\dss families with\sss the same images of\dss
subspaces\sss
$\image
P_{\qff \geq\qff \varepsilon}\trf(\dff A_{\dff x}\dff)$\nnsp.\oss
By extending\sss this homotopy\sss to\sss $x\qff \in\qff U$\sss
we get\sss a\sss trivialization\sss 
$t\fff''_{\dff U}\trf(\trf x\trf)\fff,\qff x\qff \in\qff U$\sss
strictly\sss adapted\sss to\sss $(\trf U\fff,\qff \varepsilon\trf)$ 
and such\sss that\sss
$t\fff''_{\dff U}\trf(\trf x\trf)
\off =\off
t\trf(\trf x\trf)$\sss
for\sss $x\qff \in\qff U_{\dff n}$.\oss

Let\sss us extend $t$\dss from\sss $\ske_{\trf n}\dff X$\sss
to\sss $\ske_{\trf n}\dff X\dff \cup\dff \sigma$\sss by\sss $t\fff''_{\dff U}$\sss
and denote\sss the resulting map again\sss by $t$\nnsp.\oss
Then\sss the restriction of\sss $t$\sss to\sss $U_{\dff n}\dff \cup\dff \sigma$\sss
is\dss strictly\sss adapted\sss to\sss
$(\trf U_{\dff n}\dff \cup\dff \sigma\fff,\qff \varepsilon\trf)$\nnsp.\oss
Lemma\qss \ref{atlas-independence}\qss implies\sss that\sss the
extended\sss $t$\sss is\dss strictly\sss adapted.\oss
It\sss follows\sss that\sss we can extend\sss the\sss trivialization\sss $t$\sss
to $\sigma$\nnsp.\oss
By doing\sss this for all $(\dff n\qff +\qff 1\dff)$\dnsp-simplices simultaneously,\oss
we can extend\sss $t$\sss to\sss the next\sss skeleton\sss
$\ske_{\trf n\dff +\dff 1}\dff X$\nnsp.\oss
By continuing\sss in\sss this way\sss we will\sss get\sss
a strictly\sss adapted\sss trivialization.\oss  \eproof

\mypar{Theorem.}{adapted-paracompact}
\emph{If\pss $\mathbb{A}$\sss
is\dss a strictly\trs Fredholm\dss family\sss
and\dss $X$ is\dss a\sss paracompact\sss space,\oss
then\sss there exists a\sss trivialization of\trs the\sss bundle\qss
$\mathbb{H}$\dss
strictly\sss adapted\dss to\sss $\mathbb{A}$\nnsp.\oss}

\proof
Let\sss us\sss choose an atlas\sss
$(\trf U_{\dff a}\dff,\qff \varepsilon_{\dff a}\trf)\fff,\off a\qff \in\qff \Sigma$\sss
adapted\sss to\sss $\mathbb{A}$\sss and\dss
and\dss local\dss trivializations\sss $t_{\trf a}$\sss
strictly\sss adapted\dss to\sss $(\trf U_{\dff a}\dff,\qff \varepsilon_{\dff a}\trf)$\sss
and\sss $\mathbb{A}$\nnsp.\oss
Let\sss $\Sigma\ffin$\sss and\sss $U_{\dff \sigma}\dff,\qff \varepsilon_{\dff \sigma}$\sss 
for\sss $\sigma\qff \in\qff \Sigma\ffin$\sss
be defined as in\dss Section\qss \ref{analytic-index-section},\oss 
and\sss let\sss $X_{\dff U}$\sss
be\sss the\sss topological\sss category\sss from\dss Section\qss \ref{analytic-index-section}.\oss
The space of\dss objects of\sss $X_{\dff U}$\sss is\dss the disjoint\sss union\sss
of\dss subspaces\sss $U_{\dff \sigma}$\nsp,\oss
and\dss the category\sss $X_{\dff U}$\sss is\dss defined\sss by an order\sss $\leq$\sss
on\sss $\ob\trf X_{\dff U}$\nsp.\oss
Clearly,\oss the set\sss pairs\sss $(\trf u\fff,\qff u\trf)$\nnsp,\oss
where\sss $u\qff \in\qff \ob\trf X_{\dff U}$\nsp,\oss
is\dss the union of\dss several\sss components of\dss the set\sss of\dss pairs\sss
$(\trf u\fff,\qff v\trf)$\sss of\dss comparable objects of\sss $X_{\dff U}$\nsp,\oss
i.e.\qss such\sss that\sss either\sss $u\qff \leq\qff v$\sss or\sss $v\qff \leq\qff u$\nnsp.\oss
In other\sss words,\oss the subspace of\dss the identity\sss morphisms\dss
is\dss open and closed\sss in\sss the space of\dss all\sss morphisms.\oss
In\sss the\sss terminology of\qss \cite{i2}\qss this means\sss that\sss
$\ob\trf X_{\dff U}$\sss is\dss a partially ordered space with\qss 
\emph{free equalities},\oss and\sss hence\sss $X_{\dff U}$\sss
can\sss be\sss treated as a\sss topological\sss simplicial\sss complex.\oss
In\sss particular,\oss the standard\sss geometric realization $\num{X_{\dff U}}$\sss 
is\dss equal\sss to\sss the\qss ``naive''\qss
geometric realization\sss $\bbnum{X_{\dff U}}$\nnsp.\oss
We refer\sss to\qss \cite{i2},\oss Section\qss 5,\oss for\sss the details.\oss

Let\sss $\mathbb{H}^{\dff U}$ be\sss the bundle induced\sss from\sss $\mathbb{H}$\sss by\sss 
$\num{\pr}\dff \colon\dff
\num{X_{\dff U}}\qff \ttoo\qff X$\nnsp.\oss
The family\sss $\mathbb{A}$\sss defines a family\sss $\mathbb{A}^{\fff U}$\sss of\dss
operators in\sss the bundle\sss $\mathbb{H}^{\dff U}$\dnsp.\oss
If\sss $s$\sss is\dss a homotopy\sss inverse of\sss $\num{\pr}$\sss as
in\trs Theorem\qss \ref{covering-categories},\oss
then\sss $\mathbb{H}$\sss and\sss $\mathbb{A}$\sss are equal\sss to\sss
the bundle and\sss the family of\dss operators induced\sss by\sss $s$\sss 
from\sss $\mathbb{H}^{\dff U}$\sss and\sss $\mathbb{A}^{\fff U}$ respectively.\oss
Therefore\sss it\dss is\dss sufficient\sss to prove\sss that\sss
there exists a\sss trivialization of\sss $\mathbb{H}^{\dff U}$\sss
strictly adapted\sss to\sss $\mathbb{A}^{\fff U}$\dnsp.\oss
As in\sss the proof\dss of\trs Theorem\qss \ref{adapted},\oss
we will\sss prove\sss this using an\sss induction by\sss skeletons.\oss
If\dss $\sigma\qff \in\qff \Sigma\ffin$\sss and\sss $a\qff \in\qff \sigma$\nnsp,\oss
then\sss the restriction $t_{\dff \sigma}$\sss of\sss $t_{\dff a}$\sss to\sss 
$U_{\dff \sigma}\qff \subset\qff U_{\dff a}$\sss is\dss strictly adapted\sss to $\mathbb{A}$\nnsp.\oss
Hence\sss the\sss trivializations\sss $t_{\trf \sigma}$\sss
define a\sss trivialization of\dss the restriction of\dss $\mathbb{H}^{\dff U}$\sss
to\sss the $0${\nsp}th skeleton\sss $\bbnum{\nsp\ssk_{\trf 0}\dff X_{\dff U}}$\sss
strictly adapted\sss to\sss the restriction of\sss $\mathbb{A}^{\fff U}$\dnsp.\oss

Suppose\sss that\sss we already constructed a\sss trivialization $t$ of\dss
the restriction of\sss $\mathbb{H}^{\dff U}$ to\sss the 
$n$\dnsp-skele\-ton $\bbnum{\nsp\ske_{\trf n}\dff X_{\dff U}}$\sss
strictly\sss adapted\sss to\sss the family\sss
$A^{\fff U}_{\dff u}\dff,\qff u\qff \in\qff \bbnum{\nsp\ske_{\trf n}\dff X_{\dff U}}$
for some $n\qff \geq\qff 0$\nnsp.\oss
An $(\dff n\qff +\qff 1\dff)$\dnsp-simplex of\dss $X_{\dff U}$\sss
is\dss defined\sss by a point\sss $x\qff \in\qff X$\sss and a strictly
decreasing sequence\vspace{0.75pt}
\[
\quad
\sigma\trf(\dff n\qff +\qff 1\dff)
\off \supsetneq\off
\sigma\trf(\dff n\trf)
\off \supsetneq\off
\ldots
\off \supsetneq\off
\sigma\trf(\trf 0\trf)
\]

\vspace{-12pt}\vspace{0.75pt}
of\dss elements of\sss $\Sigma\ffin$\dnsp.\oss
Therefore\sss the space of\dss $(\dff n\qff +\qff 1\dff)$\dnsp-simplices of\dss $X_{\dff U}$\sss
can\sss be identified\sss with\sss the disjoint\sss union\sss
of\dss the sets\sss $U_{\dff \sigma\trf(\dff n\qff +\qff 1\dff)}$\sss
over\sss all\sss such sequences.\oss
In\sss particular,\oss every component\sss of\dss the space of\sss
$(\dff n\qff +\qff 1\dff)$\dnsp-simplices
can\sss be identified\sss with\sss $U_{\dff \sigma}$\sss 
for some\sss $\sigma\qff \in\qff \Sigma\ffin$\nnsp.\oss

Let\sss us consider\sss such component\sss $U_{\dff \sigma}$\sss
and\sss the canonical\sss map\sss\vspace{0.75pt}
\begin{equation}
\label{simplex-in}
\quad
\Delta^{n\dff +\dff 1}\dff \times\qff U_{\dff \sigma}
\qff \ttoo\qff
\bbnum{X_{\dff U}}
\off =\off
\num{X_{\dff U}}
\pff.
\end{equation}

\vspace{-12pt}\vspace{0.75pt}
The bundle $\mathbb{H}^{\dff U}$ and\sss the family $\mathbb{A}^{\fff U}$
induce a bundle $\mathbb{H}^{\dff \sigma}$ 
and a\sss family of\dss operators $\mathbb{A}^{\fff \sigma}$ over\sss
$\Delta^{n\dff +\dff 1}\dff \times\qff U_{\dff \sigma}$\nsp,\oss
and\sss the already constructed\sss trivialization\sss $t$\sss
defines a\sss trivialization $t'$ of\sss $\mathbb{H}^{\dff \sigma}$ over\trs
$\partial\trf \Delta^{n\dff +\dff 1}\dff \times\qff U_{\dff \sigma}$\dss
strictly adapted\sss $\mathbb{A}^{\fff \sigma}$\dnsp.\oss
Since\sss the map\sss $\num{\pr}$\sss collapses
each simplex\sss to a point,\oss the fibers of\sss
$\mathbb{H}^{\dff U}$ over a simplex can\sss be\sss treated as equal,\oss
as also\sss the operators\sss $A^{\fff U}_{\dff u}$\nsp.\oss
Hence\sss the\sss trivialization\sss $t_{\trf \sigma}$\sss
defines a\sss trivialization of\sss $\mathbb{H}^{\dff \sigma}$
strictly adapted\sss to\sss $\mathbb{A}^{\fff \sigma}$\dnsp.\oss

Let\sss us\sss use\sss this\sss trivialization\sss to identify\sss
$\mathbb{H}^{\dff \sigma}$\sss with\sss 
the\sss trivial\sss bundle over\sss
$\Delta^{n\dff +\dff 1}\dff \times\qff U_{\dff \sigma}$\nsp.\oss
Now we can use\sss the contractibility of\dss 
the space of\dss polarizations 
and\dss the contractibility\sss of\sss 
the unitary\sss groups\sss $\mathcal{U}\dff(\trf K\trf)$\sss
as in\sss the proof\dss of\trs Theorem\qss \ref{adapted}\qss
to extend $t'$\sss to a\sss trivialization $t''$ 
of\sss $\mathbb{H}^{\dff \sigma}$ over\sss
$\Delta^{n\dff +\dff 1}\dff \times\qff U_{\dff \sigma}$\dss
strictly adapted\sss to\sss $\mathbb{A}^{\fff \sigma}$\dnsp.\oss
Since we are dealing\sss with a\sss topological\sss simplicial\sss
complex,\oss the map\qss (\ref{simplex-in})\qss is\dss injective
and we can extend $t$\sss to\sss the image of\dss the map\qss (\ref{simplex-in})\qss 
by\sss using $t''$\dnsp.\oss
By\sss doing\sss this simultaneously\sss for every component\sss
of\dss the space of\sss $(\dff n\qff +\qff 1\dff)$\dnsp-simplices
we can extend\sss $t$\sss to\sss next\sss skeleton\sss
$\bbnum{\nsp\ske_{\trf n\dff +\dff 1}\dff X_{\dff U}}$\nnsp.\oss
The\sss topology of\sss 
$\num{X_{\dff U}}\off =\off \bbnum{X_{\dff U}}$\sss is,\oss
by\sss the definition,\sss the direct\sss limit\sss of\dss the\sss
topologies of\dss the skeletons.\oss
Therefore,\oss by continuing\sss in\sss this way\sss
we will\sss get\sss a\sss trivialization of\dss $\mathbb{H}^{\dff U}$
strictly adapted\sss to\sss $\mathbb{A}^{\fff U}$\dnsp.\oss  \eproof

\mypar{Theorem.}{adapted-relative}
\emph{Under\sss the assumptions of\pss Theorem\qss \ref{adapted}\qss or\qss \ref{adapted-paracompact},\oss
suppose\sss that\qss $Y$\dss is\dss a subcomplex or closed subset\sss of\pss $X$\sss respectively.\oss
If\qss $t\trf(\trf y\trf)\dff,\qff y\qff \in\qff Y$\sss
is\trs a\sss strictly\sss adapted\dss trivialization\dss for\sss 
$A_{\dff y}\dff,\qff y\qff \in\qff Y$\dnsp,\oss
then\sss $t$\sss can\sss be extended\dss to a strictly\sss adapted\dss trivialization for\sss
$A_{\dff x}\dff,\qff x\qff \in\qff X$\nnsp.\oss}

\proof
The proof\trs is\dss a standard\sss modification of\dss the proofs
of\qss Theorems\qss \ref{adapted}\qss and\qss \ref{adapted-paracompact}.\oss  \eproof

\mypar{Corollary.}{adapted-unique}
\emph{Under\sss the assumptions of\pss Theorem\qss \ref{adapted},\oss
every\sss two strictly\sss adapted\dss trivializations are homotopic\sss
in\dss the class of\dss strictly\sss adapted\dss trivializations.\oss}\vspace{-0.125pt}

\proof
It\dss is\dss sufficient\sss to apply\trs Theorem\qss \ref{adapted-relative}\qss
to\sss the subset\sss $X\dff \times\dff \{\trf 0\dff,\qff 1\trf\}$
of\trs $X\dff \times\dff [\dff 0\dff,\qff 1\dff]$\nnsp.\oss  \eproof

\newpage
\mysection{Analytical\qss index\qss of\pss strictly\qss Fredholm\qss families}{analytic-index-strict}

\myuppar{The\sss topological\sss category 
$\mathcal{P}{\nsp}\hat{\mathcal{S}}$\nsp\dnsp.}
Let\sss us\sss recall\dss the definition of\dss 
$\mathcal{P}{\nsp}\hat{\mathcal{S}}$\sss
from\qss \cite{i2}.\oss
This\dss is\dss a\sss topological\sss category\sss having\sss
as objects\sss triples\sss
$(\trf V\fff,\pff H_{\dff -}\dff,\pff H_{\dff +}\trf)$\nnsp,\oss
where\sss $V$\sss is\dss a\sss finitely dimensional\sss subspace of\sss $H$\nnsp,\pss
$H_{\dff -}$\sss and\sss $H_{\dff +}$\sss are infinitely dimensional,\oss
and\vspace{0.75pt}
\[
\quad
H
\off =\off 
H_{\dff -}\dff \oplus\dff V\dff \oplus\dff H_{\dff +}
\off.
\]

\vspace{-12pt}\vspace{0.75pt}
In other\sss terms,\pss $V$\sss is\dss an object\sss of\sss $\hat{\mathcal{S}}$\sss
and\sss
$H\dff \ominus\dff V
\off =\off
H_{\dff -}\dff \oplus\dff H_{\dff +}$\sss
is\dss a polarization of\dss $H\dff \ominus\dff V$\dnsp.\oss
The\sss topology on\sss the space of\dss such\sss triples\dss is\dss defined\sss
by\sss the usual\sss topology on\sss subspaces $V$
and\sss the norm\sss topology on\sss subspaces\sss
$H_{\dff -}\dff,\off H_{\dff +}$.\oss 
This space\dss 
is\dss ordered\sss 
by\sss the relation\sss $\leq$\nnsp,\oss 
where\vspace{0.75pt}
\[
\quad
(\trf V\fff,\off H_{\dff -}\dff,\off H_{\dff +}\trf)
\off \leq\off
(\trf V\fff'\fff,\off H\fff'_{\dff -}\dff,\off H\fff'_{\dff +}\trf)
\quad
\]

\vspace{-39pt}\vspace{0.75pt}
\[
\quad
\mbox{if}\dff\quad
V\off \subset\off V\fff'
\quad
\mbox{and}\quad
H\fff'_{\dff -}\qff \subset\pff H_{\dff -}
\qff,\quad
H\fff'_{\dff +}\qff \subset\pff H_{\dff +}
\pff.
\]

\vspace{-12pt}\vspace{0.75pt}
The morphisms of\dss
$\mathcal{P}{\nsp}\hat{\mathcal{S}}$\sss
are defined\sss by\sss this order,\oss
i.e.\qss a morphism\sss $P\qff \ttoo\qff P\fff'$\sss
exists\dss if\trs and\dss only\trs if\trs
$P\qff \leq\qff P\fff'$\nnsp,\oss and\sss in\sss this case it\dss is\dss unique.\oss
Morphisms have\sss the form\vspace{0.75pt}
\[
\quad
(\trf V\fff,\off H_{\dff -}\dff,\off H_{\dff +}\dff)
\qff \ttoo\qff
(\trf U_{\dff -}\qff \oplus\qff
V\qff \oplus\qff U_{\dff +}\dff,\off 
H_{\dff -}\dff \ominus\dff U_{\dff -}\dff,\off 
H_{\dff +}\dff \ominus\dff U_{\dff +}
\trf)
\qff,
\]

\vspace{-12pt}\vspace{0.75pt}
where\sss $U_{\dff -}\dff,\pff U_{\dff +}$\sss
are finitely dimensional\sss subspaces of\dss  $H_{\dff -}\dff,\pff H_{\dff +}$\sss
There\dss is\dss a\sss forgetting\sss functor\sss
$\pi\dff \colon\dff
\mathcal{P}{\nsp}\hat{\mathcal{S}}
\qff \ttoo\qff
\hat{\mathcal{S}}$\sss
taking\sss $(\trf V\fff,\off H_{\dff -}\dff,\off H_{\dff +}\dff)$\sss
to\sss $V$\sss
and\sss taking a morphism\vspace{0.75pt}
\[
\quad
(\trf V\fff,\off H_{\dff -}\dff,\off H_{\dff +}\trf)
\off \ttoo\off
(\trf V\fff'\fff,\off H\fff'_{\dff -}\dff,\off H\fff'_{\dff +}\trf)
\quad
\]

\vspace{-12pt}\vspace{0.75pt}
to\sss the morphism\sss
$V\qff \ttoo\qff V\fff'$\sss of\dss 
$\hat{\mathcal{S}}$\sss defined\sss by\sss 
the orthogonal\sss decomposition\vspace{0.75pt}
\[
\quad
V\fff'
\off =\off
U_{\dff -}\qff \oplus\qff
V\qff \oplus\qff U_{\dff +}
\pff,
\]

\vspace{-12pt}\vspace{0.75pt}
where\sss
$U_{\dff -}
\off =\off 
H_{\dff -}\dff \ominus\qff H\fff'_{\fff -}$
and\dss
$U_{\dff +}
\off =\off 
H_{\dff +}\dff \ominus\qff H\fff'_{\fff +}$.\oss
The geometric realization\vspace{0.75pt}
\[
\quad
\num{\pi}\dff \colon\dff
\num{\mathcal{P}{\nsp}\hat{\mathcal{S}}}
\qff \ttoo\qff
\num{\hat{\mathcal{S}}}
\]

\vspace{-12pt}\vspace{0.75pt}
is\dss a\sss homotopy equivalence.\oss
See\qss \cite{i2},\oss Theorem\qss 10.3.\oss

\myuppar{The index\sss map of\dss strictly\dss Fredholm\dss families.}
As before,\oss
let\sss $\mathbb{H}$\sss be a\dss 
Hilbert\dss bundle over $X$\nnsp.\oss
Let\sss $\mathbb{A}$\sss
be a\dss strictly\dss Fredholm\dss family of\dss operators in\sss the fibers of\sss $\mathbb{H}$\nnsp.\oss
As before,\oss we\sss think about\sss $\mathbb{H}$\sss as about\sss a family\sss
$H_{\dff x}\dff,\qff x\qff \in\qff X$ of\dss Hilbert\sss spaces parameterized\sss by $X$\nnsp,\oss
and denote\sss $\mathbb{A}$\sss also as\sss 
$A_{\dff x}\dff \colon\dff H_{\dff x}\qff \ttoo\qff H_{\dff x}$,\trs
$x\qff \in\qff X$\nnsp.\oss
We will\sss assume\sss that\sss $X$\sss is\dss a\sss paracompact space.\oss
By\trs Theorem\qss \ref{adapted-paracompact}\qss there exists a\sss trivialization of\dss
the bundle $\mathbb{H}$ adapted\sss to\sss the family $\mathbb{A}$\nnsp.\oss
Let\sss us\sss fix such a\sss trivialization.\oss
Then we may consider $\mathbb{A}$ as a family of\dss operators in\sss
the fibers of\dss the\sss trivial\dss bundle\sss
$X\dff \times\dff H\qff \ttoo\qff X$\sss
and\sss hence as a family\sss
$A_{\dff x}\dff \colon\dff H\qff \ttoo\qff H$,\trs
$x\qff \in\qff X$\sss
of\dss operators\sss in\sss $H$\nnsp.

The\sss fixed\sss trivialization of\dss the bundle $\mathbb{H}$ allows\sss to
identify\sss the\sss topological\sss category\sss 
$\hat{\mathcal{S}}\dff(\trf \mathbb{H}\trf)$\sss
with\sss the product\sss
$X\dff \times\dff \hat{\mathcal{S}}$\nnsp,\oss
where $X$\sss is\dss considered as a category\sss having $X$ as\sss the space
of\dss objects and only\sss the identity\sss morphisms.\oss
Moreover,\oss the functor\vspace{1.5pt}
\[
\quad
\mathbb{A}_{\qff U,\qff \bm{\varepsilon}}\dff \colon\dff
X_{\dff U}\qff \ttoo\qff \hat{\mathcal{S}}\dff(\trf \mathbb{H}\trf)
\off =\off
X\dff \times\dff \hat{\mathcal{S}}
\]

\vspace{-12pt}\vspace{1.5pt}
constructed\sss in\dss Section\qss \ref{analytic-index-section}\qss for\sss $\mathbb{A}$\sss
as a family of\dss operators in\sss the fibers of\dss $\mathbb{H}$\sss
can\sss be identified\sss with\sss the product\sss of\dss the projection\sss functor\sss
$\pr\dff \colon\dff X_{\dff U}\qff \ttoo\qff X$\sss
with\sss the functor\vspace{1.5pt}
\[
\quad
\mathbb{A}_{\qff U,\qff \bm{\varepsilon}}\dff \colon\dff
X_{\dff U}\qff \ttoo\qff \hat{\mathcal{S}}
\]

\vspace{-12pt}\vspace{1.5pt}
constructed\sss for\sss $\mathbb{A}$ as a family of\dss operators in\sss $H$\nnsp.\oss
Up\sss to homotopy\sss the geometric realization\vspace{1.5pt}
\[
\quad
\num{\mathbb{A}_{\qff U,\qff \bm{\varepsilon}}}
\qff \colon\qff
\num{X_{\dff U}}
\qff \ttoo\qff 
\num{\hat{\mathcal{S}}}
\]

\vspace{-12pt}\vspace{1.5pt}
does not\sss depend on\sss the choice of\dss the adapted\dss trivialization,\oss
because\sss any\sss two choices are homotopic in\sss the class of\dss 
trivializations of\dss $\mathbb{H}$\nnsp.\oss
This fact\sss follows\sss from\sss the contractibility of\dss
the unitary\sss group\sss $\mathcal{U}\dff(\trf H\trf)$
and does not\sss depend on\trs Corollary\qss \ref{adapted-unique}.\oss

The strictly\dss Fredholm\dss property allows\sss to\sss lift\sss the functor\sss
$\mathbb{A}_{\qff U,\qff \bm{\varepsilon}}\dff \colon\dff
X_{\dff U}\qff \ttoo\qff \hat{\mathcal{S}}$\sss
to a functor\sss to\sss the category\sss
$\mathcal{P}{\nsp}\hat{\mathcal{S}}$\dnsp.\oss
In\sss more details,\oss suppose\sss that\sss
$U_{\dff a}\dff,\pff a\qff \in\qff \Sigma$\sss
is\dss an open covering\sss of\dss $X$\nnsp,\oss
and\sss 
$\varepsilon_{\dff a}\qff >\qff 0\dff,\pff a\qff \in\qff \Sigma$\sss
are such\sss that\sss
$(\trf A_{\dff z}\dff,\qff \varepsilon_{\dff a}\trf)$\sss
is\dss an enhanced\sss operator for 
$z\qff \in\qff U_{\dff a}$\nsp,\oss
i.e.\qss such\sss that\sss the pairs\sss
$(\trf U_{\dff a}\dff,\qff \varepsilon_{\dff a}\trf)$\sss
are adapted\sss to $\mathbb{A}$\nnsp.\oss
As in\dss Section\qss \ref{analytic-index-section},\oss
for each\sss 
$\sigma\qff \in\qff \Sigma^{\dff \fin}$\dss 
let\sss 
\vspace{1.5pt}
\[
\varepsilon_{\dff \sigma}
\off =\off
\min\nolimits_{\pff a\qff \in\qff \sigma}\qff \varepsilon_{\dff a}
\pff.
\]

\vspace{-12pt}\vspace{1.5pt}
Let\sss us\sss assign\sss to each object\sss $(\trf z\fff,\qff \sigma\trf)$\sss
of\trs $X_{\dff U}$\sss the object\sss
$(\trf V\fff,\pff H_{\dff -}\dff,\pff H_{\dff +}\trf)$\sss
of\trs $\mathcal{P}{\nsp}\hat{\mathcal{S}}$\dnsp,\oss
where\vspace{1.5pt}
\[
\quad
V
\off =\off
\image P_{\qff [\dff -\dff \varepsilon\fff,\qff \varepsilon\dff]}\trf(\trf A_{\dff z}\trf)
\quad
\mbox{and}
\]

\vspace{-36pt}
\[
\quad
H_{\dff -}
\off =\off
\image P_{\qff \leq\qff -\qff \varepsilon}\trf(\trf A_{\dff z}\trf)
\qff,\quad
H_{\dff +}
\off =\off
\image P_{\qff \geq\qff \varepsilon}\trf(\trf A_{\dff z}\trf)
\qff.
\]

\vspace{-12pt}\vspace{1.5pt}
This\sss assignment\sss preserves\sss the orders 
and\sss hence defines a functor\vspace{1.5pt}
\[
\quad
\mathcal{P}\mathbb{A}_{\qff U,\qff \bm{\varepsilon}}\dff \colon\dff
X_{\dff U}
\qff \ttoo\qff 
\mathcal{P}{\nsp}\hat{\mathcal{S}}
\pff.
\]

\vspace{-12pt}\vspace{1.5pt}
Clearly,\oss the following diagram\dss is\dss commutative.\vspace{-1.5pt}\vspace{-1pt}
\[
\quad
\begin{tikzcd}[column sep=bboom, row sep=hqboom]
&
\mathcal{P}{\nsp}\hat{\mathcal{S}}
\arrow[dd, "\dis \dff \pi"]
\\
X_{\dff U}
\arrow[ru, "\dis \mathcal{P}\mathbb{A}_{\qff U,\qff \bm{\varepsilon}}" near end]
\arrow[rd, "\dis \mathbb{A}_{\qff U,\qff \bm{\varepsilon}}"' near end]
&
\\
&
\hat{\mathcal{S}}\qff.
\end{tikzcd}
\]

\vspace{-12pt}\vspace{-1.5pt}
In other words,\oss the functor\sss
$\mathcal{P}\mathbb{A}_{\qff U,\qff \bm{\varepsilon}}$\sss
is\dss a\dss lift\dss of\dss the functor\sss
$\mathbb{A}_{\qff U,\qff \bm{\varepsilon}}$\nnsp.\oss
By\sss passing\sss to\sss the geometric realizations
we get\sss a continuous map\vspace{1.5pt}\vspace{0.125pt}
\[
\quad
\num{\mathcal{P}\mathbb{A}_{\qff U,\qff \bm{\varepsilon}}}\qff \colon\qff
\num{X_{\dff U}}
\qff \ttoo\qff 
\num{\mathcal{P}{\nsp}\hat{\mathcal{S}}}
\pff.
\]

\vspace{-12pt}\vspace{1.5pt}\vspace{0.125pt}
Corollary\qss \ref{adapted-unique}\qss implies\sss that\sss
the homotopy class of\dss this map\dss is\dss independent\sss from\sss 
the choice of\dss adapted\dss trivialization.\oss 
Since\dss $X$\sss is\dss triangulable and\sss hence paracompact,\oss there\dss is\dss a homotopy\sss
inverse $s$ of\dss $\num{\pr}$\sss and\sss get\sss a map\vspace{2pt}
\[
\quad
\num{\mathcal{P}\mathbb{A}_{\qff U,\qff \bm{\varepsilon}}}
\qff \circ\qff
s
\qff \colon\qff
X
\qff \ttoo\qff 
\num{\mathcal{P}{\nsp}\hat{\mathcal{S}}}
\pff,
\]

\vspace{-12pt}\vspace{2pt}
which we call\sss a\qss
\emph{polarized\dss index\dss map}\pss
of\dss $\mathbb{A}$\nnsp.\oss
Its\sss homotopy class carries\sss the same information as\sss the index,\oss
but\sss the map\sss
$\num{\mathcal{P}\mathbb{A}_{\qff U,\qff \bm{\varepsilon}}}
\qff \circ\qff
s$ 
keeps information
about\sss subspaces\sss
$\image P_{\qff \geq\qff \varepsilon}\trf(\trf A_{\dff z}\trf)$\nnsp.

\myuppar{The\dss Grassmannian\dss bundle of\dss $\mathbb{A}$\nnsp.}
Recall\sss that\sss 
two closed subspaces\sss $K\fff,\qff K\fff'\qff \subset\qff H$\sss
are said\sss to be\qss \emph{commensurable}\pss if\dss the intersection\sss
$K\dff \cap\dff K\fff'$\sss has finite codimension\sss
in\sss both $K$ and\sss $K\fff'$\nnsp.\oss
Given a subspace\sss $K\qff \subset\qff H$\sss of\dss infinite dimension
and codimension\sss in\sss $H$\nnsp,\oss
the\qss
\emph{restricted\dss Grassmannian}\dss
$\gr\trf(\trf K\trf)$\sss is\dss the space of\dss closed subspaces of\dss $H$\sss
commensurable with\sss $K$\sss with\sss the norm\sss topology.\oss
Up\sss to a homeomorphism\sss 
$\gr\trf(\trf K\trf)$\sss do not\sss
depend on\sss $K$\dnsp.\oss
For each $x\qff \in\qff X$\trs let\sss\vspace{2pt}
\[
\quad
\gr\trf(\trf x\trf)
\off =\off
\gr\qff\bigl(\qff
\image P_{\qff \geq\qff \varepsilon}\trf(\trf A_{\dff x}\trf)
\qff\bigr)
\pff
\]

\vspace{-12pt}\vspace{2pt}
for some\sss $\varepsilon\qff >\qff 0$\sss such\sss that\sss
$(\trf A_{\dff x}\fff,\qff \varepsilon\trf)$\sss
is\dss an enhanced operator.\oss
Clearly,\pss $\gr\trf(\trf x\trf)$ does not\sss 
depend on\sss the choice of\sss $\varepsilon$\nnsp.\oss
Since\sss the family\sss $\mathbb{A}$\dss
is\dss strictly\trs Fredholm,\oss
the family\sss 
$\gr\trf(\trf x\trf)\dff,\qff
x\qff \in\qff X$\sss
forms a\sss locally\sss trivial\sss bundle having as\sss the fiber\sss
the restricted\dss Grassmannian\sss $\gr\trf(\trf x\trf)$\nnsp.\oss
We will\sss denote\sss the\sss total\sss space of\dss this bundle by\sss
$\gr\trf(\trf \mathbb{A}\trf)$\sss 
and\sss the bundle\sss itself\dss by\vspace{2pt}
\[
\quad
\bm{\pi}\trf(\trf \mathbb{A}\trf)\dff \colon\dff
\gr\trf(\trf \mathbb{A}\trf)
\qff \ttoo\qff 
X
\pff.
\]

\vspace{-12pt}\vspace{2pt}
We will\sss call\sss sections of\dss this bundle\qss
\emph{weak\sss spectral\sss sections}\pss
for\sss the family\sss $\mathbb{A}$\nnsp.\oss

\myuppar{A universal\dss Grassmannian\sss bundle.}
Given an object\dss 
$P\off =\off (\trf V\fff,\pff H_{\dff -}\dff,\pff H_{\dff +}\trf)$\sss
of\trs $\mathcal{P}{\nsp}\hat{\mathcal{S}}$\dnsp,\dff\oss 
let\sss\vspace{1.5pt}
\[
\quad
\gr\trf(\trf P\trf)
\off =\off
\gr\trf(\trf H_{\dff +}\trf)
\pff.
\]

\vspace{-12pt}\vspace{1.5pt}
The points of\dss the geometric realization\sss $\num{\mathcal{P}{\nsp}\hat{\mathcal{S}}}$\sss 
can\sss be represented\sss by\sss weighted sums\vspace{2pt}
\begin{equation}
\label{p-weighted-sum}
\quad
t_{\trf 0}\trf P_{\dff 0}\pff +\pff
t_{\dff 1}\trf P_{\dff 1}\pff +\pff
\ldots\pff +\pff
t_{\dff n}\trf P_{\dff n}
\pff,
\end{equation}

\vspace{-12pt}\vspace{2pt}
where\sss
$t_{\dff 0}\dff,\qff t_{\dff 1}\dff,\qff \ldots\dff,\qff t_{\dff n}
\qff \geq\qff
0$\nnsp,\qss
$t_{\dff 0}\qff +\qff t_{\dff 1}\qff +\qff \ldots\qff +\qff t_{\dff n}
\off =\off
1$\nnsp,\oss
and\dss 
$P_{\dff 0}\fff,\qff
P_{\dff 1}\fff,\qff
\ldots\fff,\qff
P_{\dff n}$\sss 
are objects of\dss 
$\mathcal{P}{\nsp}\hat{\mathcal{S}}$
such\sss that\sss
$P_{\dff 0}\off <\off
P_{\dff 1}\off <\off
\ldots\off <\off
P_{\dff n}$\nsp.\oss
See\qss \cite{i2},\oss Section\qss 13.\oss
One can easily check\sss that\dss if\trs $P\fff,\qff P\fff'$\sss
are objects of\dss $\mathcal{P}{\nsp}\hat{\mathcal{S}}$ and\dss
$P\qff \leq\qff P\fff'$\nnsp,\oss
then\sss
$\gr\trf(\trf P\trf)
\off =\off
\gr\trf(\trf P\fff'\trf)$\nnsp.\oss
It\sss follows\sss that\sss the restricted\dss Grassmannian\sss
$\gr\trf(\trf P_{\dff i}\trf)$\sss 
does not\sss depends on\sss $i$\sss and\sss hence depends only on\sss
the point\sss 
$p\qff \in\qff \num{\mathcal{P}{\nsp}\hat{\mathcal{S}}}$\sss
represented\dss by\sss the weighted sum\qss (\ref{p-weighted-sum}).\oss
In\sss particular,\oss we can set\sss
$\gr\trf(\trf p\trf)
\off =\off
\gr\trf(\trf P_{\dff i}\trf)$\nnsp.\oss
Let\sss $\mathbf{G}$\sss be\sss the space of\dss pairs\sss
$(\trf p\fff,\qff K\trf)$\sss 
such\sss that\sss
$p\qff \in\qff \num{\mathcal{P}{\nsp}\hat{\mathcal{S}}}$\sss
and\dss $K\qff \in\qff \gr\trf(\trf p\trf)$\nnsp,\oss
and\sss let\vspace{0pt}
\[
\quad
\bm{\pi}\dff \colon\dff
\mathbf{G}
\qff \ttoo\qff 
\num{\mathcal{P}{\nsp}\hat{\mathcal{S}}}
\]

\vspace{-12pt}\vspace{0pt}
be\sss the projection
$(\trf p\fff,\qff K\trf)
\off \longmapsto\off
p$\nnsp.\oss
The map $\bm{\pi}$\sss
is\dss a\sss locally\sss trivial\sss bundle
having\sss restricted\dss Grassmannians\sss $\gr\trf(\trf P\trf)$\sss as fibers.\oss
See\qss \cite{i2},\oss the beginning of\trs Section\qss 13.\oss

\myuppar{The\sss polarized\sss index\sss maps and\dss
Grassmannian\dss bundles.}
As above,\oss
let\sss $X$\sss be a\sss paracompact\sss space,\pss 
$\mathbb{H}$\sss be a\dss Hilbert\dss bundle over $X$\nnsp,\oss 
and $\mathbb{A}$\sss be a strictly\dss Fredholm\dss family of\dss operators in\sss
the fibers of\dss $\mathbb{H}$\nnsp.\oss
Then we can construct\sss a\sss functor\vspace{0pt}
\[
\quad
\mathcal{P}\mathbb{A}_{\qff U,\qff \bm{\varepsilon}}\dff \colon\dff
X_{\dff U}
\qff \ttoo\qff 
\mathcal{P}{\nsp}\hat{\mathcal{S}}
\pff.
\]

\vspace{-12pt}\vspace{0pt}
By\sss choosing a homotopy\sss inverse 
$s\qff \colon\qff
X\qff \ttoo\qff \num{X_{\dff U}}$\dss 
of\dss 
$\num{\pr}$
we get\sss the polarized\sss index\sss map\vspace{1.5pt}
\[
\quad
\num{\mathcal{P}\mathbb{A}_{\qff U,\qff \bm{\varepsilon}}}
\qff \circ\qff
s
\qff \colon\qff
X
\qff \ttoo\qff 
\num{\mathcal{P}{\nsp}\hat{\mathcal{S}}}
\pff.
\]

\vspace{-12pt}\vspace{1.5pt}
By\trs Theorem\qss \ref{covering-categories}\qss
we can assume\sss that\sss
$\num{\pr}\pff \circ\pff s
\off =\off
\id_{\qff X}$\nsp.\oss

\mypar{Theorem.}{induced-bundle}
\emph{Suppose\sss that\trs $X$ is\dss paracompact.\oss
The bundle\dss
$\bm{\pi}\trf(\trf \mathbb{A}\trf)\dff \colon\dff
\gr\trf(\trf \mathbb{A}\trf)
\qff \ttoo\qff 
X$\sss
is\dss equal\sss to\sss the bundle\sss
induced\sss from\dss
$\bm{\pi}\dff \colon\dff
\mathbf{G}
\qff \ttoo\qff 
\num{\mathcal{P}{\nsp}\hat{\mathcal{S}}}$\trs
by\dss the polarized\dss index\dss map\dss
$\num{\mathcal{P}\mathbb{A}_{\qff U,\qff \bm{\varepsilon}}}
\qff \circ\qff
s$\nnsp.\oss}

\proof
Since\sss
$\num{\pr}\pff \circ\pff s
\off =\off
\id_{\qff X}$\nsp,\oss
for every\sss $x\qff \in\qff X$\sss the point\vspace{1.0pt}
\[
\quad
\num{\mathcal{P}\mathbb{A}_{\qff U,\qff \bm{\varepsilon}}}
\qff \circ\qff
s\qff(\trf x\trf)
\]

\vspace{-12pt}\vspace{1.0pt}
is\dss equal\sss to a weighted sum\qss (\ref{p-weighted-sum})\qss
such\sss that\sss for\sss each\dss
$P_{\fff i}
\off =\off
(\trf V_{\fff i}\fff,\pff H_{\dff i\qff -}\dff,\pff H_{\dff i\qff +}\trf)$\sss
we have\vspace{1.0pt}
\[
\quad
H_{\dff i\qff +}
\off =\off
\image P_{\qff \geq\qff \varepsilon\trf(\trf i\trf)}\qff(\trf A_{\dff x}\trf)
\]

\vspace{-12pt}\vspace{1.0pt}
for some\sss $\varepsilon\trf(\trf i\trf)\qff >\qff 0$\nnsp.\oss
It\dss follows\sss that\sss the fiber of\dss the induced\sss bundle over\sss $x$\sss
is\dss equal\sss to\sss $\gr\trf(\trf x\trf)$\nnsp,\oss
i.e.\qss to\sss the fiber of\dss $\bm{\pi}\trf(\trf \mathbb{A}\trf)$\sss
over\sss $x$\nnsp.\oss
Therefore\sss the induced\sss bundle\dss is\dss equal\sss to
$\bm{\pi}\trf(\trf \mathbb{A}\trf)$\nnsp.\oss  \eproof

\mypar{Theorem.}{weak-spectral-sections}
\emph{Suppose\sss that\trs $X$ is\dss paracompact.\oss
The analytical\dss index of\trs 
a strictly\trs Fredholm\dss family\sss $\mathbb{A}$\sss vanishes\dss
if\trs and\sss only\trs if\trs
there exists\sss a weak\sss spectral\sss section\sss 
for\sss $\mathbb{A}$\nnsp.\oss}\vspace{-0.125pt}

\proof
Suppose\sss that\sss the analytical\dss index\sss vanishes.\oss
Since\sss
$\num{\pi}\dff \colon\dff
\num{\mathcal{P}{\nsp}\hat{\mathcal{S}}}
\qff \ttoo\qff
\num{\hat{\mathcal{S}}}$\sss
is\dss a homotopy equivalence,\oss
the polarized\sss analytical\dss index also vanishes.\oss
Together\sss with\trs Theorem\qss \ref{induced-bundle}\qss
this implies\sss that\sss the bundle\sss
$\bm{\pi}\trf(\trf \mathbb{A}\trf)\dff \colon\dff
\gr\trf(\trf \mathbb{A}\trf)
\qff \ttoo\qff 
X$\sss
is\dss trivial.\oss
Hence\sss it\sss admits a section.\oss
Conversely,\oss if\dss there exists a weak\sss spectral\sss section\sss for $\mathbb{A}$\nnsp,\oss
then\sss the polarized\sss index\sss map can\sss be\sss lifted\sss to a map\sss
$X\qff \ttoo\qff \mathbf{G}$\nnsp.\oss
By\qss \cite{i2},\oss Theorem\qss 13.6,\oss the space $\mathbf{G}$\sss
is\dss contractible,\oss
and\sss hence\sss the polarized\sss index\sss map\dss and\sss index\sss map are
homotopic\sss to constant\sss maps.\oss   \eproof

\mysection{Discrete-spectrum\qss families}{discrete-spectrum}

\myuppar{Discrete-spectrum\dss families.}
As before,\pss
$\mathbb{H}$\sss is\dss a\dss 
Hilbert\dss bundle over $X$\nnsp,\oss
thought\sss also a family\sss
$H_{\dff x}\dff,\qff x\qff \in\qff X$ of\dss of\dss its fibers.\oss
Let\sss $\mathbb{A}$\sss
be a\sss family of\dss self-adjoint\sss operators in\sss the fibers of\sss $\mathbb{H}$\nnsp,\oss
denoted also by\sss 
$A_{\dff x}\dff \colon\dff H_{\dff x}\qff \ttoo\qff H_{\dff x}$,\trs
$x\qff \in\qff X$\nnsp.\oss
We will\sss say\sss that\sss $\mathbb{A}$\sss is\dss a\qss
\emph{discrete-spectrum\dss family}\oss if\dss for every\sss
$\lambda\qff \in\qff \rrr$\sss the family\sss $\mathbb{A}\qff -\qff \lambda$\sss
is\dss a\dss Fredholm\dss family.\oss
Of\dss course,\oss in\sss the\sss last\sss formula $\lambda$ stands for 
$\lambda$\sss times\sss the family of\dss the identity operators.\oss
As usual,\oss we assume\sss that\sss operators $A_{\dff x}$\sss
are neither essentially\sss positive,\oss
nor essentially\sss negative.\oss

If\sss $\mathbb{A}$\sss is\dss a discrete-spectrum\dss family,\oss
then\sss every\sss operator\sss $A_{\dff x}$\sss has discrete spectrum,\oss
i.e.\qss its essential\sss spectrum\dss is\dss empty.\oss
Indeed,\oss for every\sss $\lambda\qff \in\qff \rrr$\sss
the operator\sss $A_{\dff x}\qff -\qff \lambda$\sss is\dss Fredholm\sss
and\sss hence $0$ does not\sss belong\sss to its essential\sss spectrum.\oss
It\sss follows\sss that\sss $\lambda$\sss
does not\sss belong\sss to its essential\sss spectrum of\sss $A_{\dff x}$\nsp.\oss
Since\sss this\dss is\dss true for  every\sss $\lambda\qff \in\qff \rrr$\nnsp,\oss
the essential\sss spectrum of\sss $A_{\dff x}$\sss is\dss empty,\oss
i.e.\trs $A_{\dff x}$\sss indeed\sss has discrete spectrum.\oss
But\dss being a discrete-spectrum\dss family\dss
is\dss a much stronger\sss property\dss than\sss being a\dss
Fredholm\dss family of\dss operators with discrete spectrum.\oss
Since operators $A_{\dff x}$\sss have discrete spectrum,\oss
they cannot\sss be bounded.\oss
By\sss this reason\sss in\sss this section we work\sss in\sss the
framework of\dss closed densely defined operators.\oss

For\sss the rest\sss of\dss this section\sss $\mathbb{A}$\sss is\dss assumed\sss
to be a discrete-spectrum\sss family.\oss
For each\sss $x\qff \in\qff X$\sss
the spectrum\sss $\sigma\dff(\trf A_{\dff x}\trf)$\sss
consists of\dss a double infinite 
sequence of\dss eigenvalues\vspace{1.5pt}\vspace{-0.125pt}
\[
\quad
\ldots
\off <\off
\lambda_{\dff -\dff 1}\dff(\dff x\trf)
\off <\off
\lambda_{\dff 0}\dff(\dff x\trf)
\off <\off
\lambda_{\dff 1}\dff(\dff x\trf)
\off <\off
\ldots
\]

\vspace{-12pt}\vspace{1.5pt}
and\sss the eigenspaces\sss $H_{\dff n}\dff(\dff x\trf)$\sss corresponding\sss to each\sss 
$\lambda_{\dff n}\dff(\dff x\trf)$\nnsp,\qss $n\qff \in\qff \zzz$\nnsp,\oss
are finitely dimensional.\oss
Moreover,\oss the spectrum\sss 
$\sigma\dff(\trf A_{\dff x}\trf)$\sss 
has no accumulation\sss points in\sss $\rrr$ and\sss hence\vspace{1.0pt}
\[
\quad
\lim\nolimits_{\dff n\dff \to\dff -\dff \infty}\qff \lambda_{\dff n}\dff(\dff x\trf)
\off =\off 
-\qff \infty
\off\quad
\mbox{and}\off\quad
\lim\nolimits_{\dff n\dff \to\dff \infty}\qff \lambda_{\dff n}\dff(\dff x\trf)
\off =\off 
\infty
\pff.
\]

\vspace{-12pt}\vspace{1.0pt}
The eigenspaces\sss $H_{\dff n}\dff(\dff x\trf)$\sss define a decomposition
of\dss $H_{\dff x}$\sss into\sss an orthogonal\sss direct\sss sum\vspace{1.0pt}
\begin{equation}
\label{spectral-decomposition}
\quad
H_{\dff x}
\off =\off
\bigoplus\nolimits_{\dff n\qff \in\qff \zzz}\qff H_{\dff n}\dff(\dff x\trf)
\end{equation}

\vspace{-12pt}\vspace{1.0pt}
Strictly speaking,\oss the numberings of\dss eigenvalues\sss 
$\lambda_{\dff n}\dff(\dff x\trf)$\sss
and\sss of\dss subspaces\sss $H_{\dff n}\dff(\dff x\trf)$\sss
are well\sss defined only up\sss to
shifts\sss $n\off \longmapsto\off n\qff +\qff a$\sss of\dss subscripts,\oss
where\sss $a$\sss is\dss an\sss integer.\oss
But,\oss in\sss fact,\oss only\sss the order\sss $n\qff <\qff m$\sss on\sss $\zzz$\sss
matters,\oss and\dss the order\dss is\dss well\sss defined.\oss

The decomposition\qss (\ref{spectral-decomposition})\qss depends continuously on\sss $x$\sss
only\sss in\sss a relatively\sss weak\sss sense.\oss
Let\sss $x\qff \in\qff X$\sss and\sss $a\fff,\qff b\qff \in\qff \rrr$\nnsp.\oss
Suppose\sss that\sss $a\qff <\qff b$\sss and\sss that\sss
$a\fff,\qff b\qff \not\in\qff \sigma\dff(\trf A_{\dff x}\trf)$\nnsp.\oss
Then\sss there exists a neighborhood\sss $U_{\dff x}\qff \subset\qff X$\sss of\sss $x$\sss
such\sss that\sss
$a\fff,\qff b\qff \not\in\qff \sigma\dff(\trf A_{\dff y}\trf)$\sss
for every\sss $y\qff \in\qff U_{\dff x}$\nsp.\oss
Let\vspace{0.75pt}
\[
\quad
H_{\dff y}\trf (\trf a\fff,\qff b\trf)
\off =\off
\bigoplus\qff H_{\dff n}\dff(\trf y\trf)
\pff,
\]

\vspace{-12pt}\vspace{0.75pt}
where\sss the direct\sss sum\dss is\dss
taken over $n$ such\sss that\sss 
$a\qff <\qff \lambda_{\dff n}\dff(\trf y\trf)\qff <\qff b$\nnsp.\oss
Clearly,\pss
$H_{\dff y}\trf (\trf a\fff,\qff b\trf)$\sss
does not\sss depend on\sss the choice of\dss numbering of\dss the eigenvalues.\oss
Since\sss the operators\sss $A_{\dff y}\qff -\qff \lambda$\sss are\dss
Fredholm\dss for all\sss $\lambda\qff \in\qff [\dff a\fff,\qff b\trf]$\nnsp,\oss
this sum\dss is\dss finite and\dss the spaces\sss
$H_{\dff y}\trf (\trf a\fff,\qff b\trf)$\sss are finitely dimensional.\oss
The decomposition\qss (\ref{spectral-decomposition})\qss
continuously depends on\sss the parameter\sss in\sss the sense\sss that\sss the subspace\sss
$H_{\dff y}\trf (\trf a\fff,\qff b\trf)$\sss
continuously depends on\sss $y\qff \in\qff U_{\dff x}$\nsp.\oss
This follows\sss from\sss the assumption\sss that\sss the families\sss
$\mathbb{A}\qff -\qff \lambda$\sss are\dss Fredholm.\oss
Since\sss spaces\sss $H_{\dff y}\dff (\trf a\fff,\qff b\trf)$\sss
are finitely dimensional,\oss there\dss is\dss only one natural\sss notion
of\dss continuous dependence of\dss these subspaces on $y$\nnsp.\oss

The definition of\dss subspaces\sss
$H_{\dff y}\trf (\trf a\fff,\qff b\trf)$\sss
makes sense also for\sss $a\off =\off -\qff \infty$\sss
and\sss $b\off =\off \infty$\nnsp.\oss
Clearly,\pss
$H_{\dff y}\trf (\trf -\qff \infty\fff,\qff \infty\trf)
\off =\off
H_{\dff y}$,\oss
and\sss the fibers\sss $H_{\dff y}$\sss should\sss be\sss thought\sss
as continuously depending on $y$\nnsp.\oss
For\sss general\sss bundles\sss $\mathbb{H}$\sss
and\sss $a\fff,\qff b\qff \in\qff \rrr$\sss 
there\dss is\dss even\sss no suitable notion of\dss
the continuous dependence of\dss
$H_{\dff y}\trf (\trf -\qff \infty\fff,\qff b\trf)$
and\sss
$H_{\dff y}\trf (\trf a\fff,\qff \infty\trf)$\sss
on\sss $y$\nnsp.\oss
But\sss if\dss $\mathbb{A}$\sss is\dss strictly\dss Fredholm,\oss
then\sss these subspaces continuously depend on $y$ almost\sss
by\sss the definition.\oss
In\sss more details,\oss 
we can assume\sss that\sss for some\sss $\varepsilon\qff >\qff 0$\sss 
the pair\sss $(\trf U_{\dff x}\dff,\qff \varepsilon\trf)$\sss
is\dss adapted\sss to $\mathbb{A}$\nnsp.\oss
If\dss $\varepsilon\qff <\qff a$\nnsp,\oss
then\vspace{1.5pt}
\[
\quad
\image
P_{\qff \geq\qff \varepsilon}\trf(\dff A_{\dff y}\dff)
\off =\off
\image
P_{\qff [\dff \varepsilon\fff,\dff a\dff]}\trf(\dff A_{\dff y}\dff)
\qff \oplus\qff
H_{\dff y}\trf (\trf a\fff,\qff \infty\trf)
\pff.
\]

\vspace{-12pt}\vspace{1.5pt}
In an adapted\sss trivialization\sss the\sss left\sss hand side of\dss this equality\sss
continuously depends on $y$ in\sss the norm\sss topology,\oss
and\sss the summand\sss
$\image
P_{\qff [\dff \varepsilon\fff,\dff a\dff]}\trf(\dff A_{\dff y}\dff)
\off =\off
H_{\dff y}\trf (\trf \varepsilon\fff,\qff a\trf)$\sss
also continuously depends on $y$\nnsp.\oss
Hence\sss $H_{\dff y}\trf (\trf a\fff,\qff \infty\trf)$\sss
continuously depends on $y$\nnsp.\oss 
The case of\dss $a\qff <\qff \varepsilon$\sss is\dss similar,\oss
and\sss the case of\dss $a\off =\off \varepsilon$\sss is\dss trivial.\oss

\myuppar{Spectral\sss sections.}
Let\sss us\sss assume\sss that\sss $\mathbb{A}$\sss is\dss not\sss
only a dis\-crete-spectrum\sss family,\oss but\sss also\dss is\dss
a strictly\dss Fredholm\dss family.\oss
Then\sss the\dss Grassmannian\dss bundle\sss
$\bm{\pi}\trf(\trf \mathbb{A}\trf)\dff \colon\dff
\gr\trf(\trf \mathbb{A}\trf)
\qff \ttoo\qff 
X$\sss
is\dss defined.\oss
A section\sss 
$\mathbb{S}\dff \colon\dff
X
\qff \ttoo\qff
\bm{\pi}\trf(\trf \mathbb{A}\trf)$\sss
of\dss this bundle can\sss be considered as a family
of\dss subspaces\sss $S_{\fff x}\qff \subset\qff H_{\dff x}$\nsp,\qss
$x\qff \in\qff X$\sss such\sss that\sss $S_{\fff x}$\sss
is\dss commensurable with\vspace{0.0pt}
\[
\quad
H_{\qff \geq\dff k}\trf(\dff x\trf)
\off =\off
\bigoplus\nolimits_{\dff n\qff \geq\qff k}\qff H_{\dff n}\trf(\trf x\trf)
\pff,
\]

\vspace{-12pt}\vspace{0.0pt}
for some\sss $k\qff \in\qff \zzz$\sss and every\sss $x$\nnsp,\oss
or,\oss equivalently,\oss
with\sss\vspace{0pt}
\[
\quad
H_{\dff y}\trf (\trf a\fff,\qff \infty\trf)
\]

\vspace{-12pt}\vspace{0pt}
for some\sss $a\qff \not\in\qff \sigma\dff(\trf A_{\dff x}\trf)$\nnsp.\oss
Clearly,\oss these condition do not\sss depend on\sss the choice of\dss $k$\nnsp,\oss
the numberings of\dss eigenvalues,\oss or\sss of\sss $a$\nnsp.\oss
Of\dss course,\oss we are interested\sss in\sss continuous sections.\oss
Since $\mathbb{A}$\sss is\dss strictly\dss Fredholm,\oss
the continuity of\dss a family\sss
$S_{\dff x}\dff,\qff x\qff \in\qff X$\sss
is\dss well\sss defined.\oss

Following\trs Melrose\dss and\dss Piazza\qss \cite{mp}\qss we say\sss that\sss
a continuous family\sss $S_{\fff x}\dff,\qff x\qff \in\qff X$\sss
is\dss a\qss 
\emph{spectral\sss section}\pss of\sss 
$A_{\dff x}\dff,\qff x\qff \in\qff X$\dss
if\dss there exists a continuous function\dss
$r\dff \colon\dff X\qff \ttoo\qff \rrr_{\trf >\dff 0}$\sss 
such\sss that\dss\vspace{1.5pt}\vspace{0.125pt}
\begin{equation}
\label{mp-section}
\quad
P_{\qff [\trf r\dff(\dff x\trf)\fff,\pff \infty\qff)}\qff(\dff A_{\dff x}\dff)
\off \subset\off\dff
S_{\fff x}
\off \subset\off\dff
\image
P_{\qff [\qff -\qff r\dff(\dff x\trf)\fff,\pff \infty\qff)}\qff(\dff A_{\dff x}\dff)
\pff
\end{equation}

\vspace{-12pt}\vspace{1.5pt}\vspace{0.125pt}
for every\sss $x\qff \in\qff X$\nnsp.\oss
If\dss $X$\sss is\dss compact,\oss the function\sss 
$x\off \longmapsto\off r\dff(\dff x\trf)$\sss can\sss be replaced\sss by
a constant.\oss
Note\sss that\sss $r\trf(\trf x\trf)$\sss is\dss allowed\sss to be an
eigenvalue of\dss $A_{\dff x}$\nsp.\oss
Locally,\oss one can always\sss ensure\sss that\sss
$r\trf(\trf x\trf)$\sss is\dss not\sss an eigenvalue,\oss
but,\oss in\sss general,\oss one cannot\sss do\sss this globally.\oss

\myuppar{Convex\sss combinations of\dss projections.}
Our next\sss goal\dss is\dss to prove\sss that\sss every\sss
weak\sss spectral\sss section can\sss be deformed\dss into\sss
a spectral\sss section.\oss
To\sss this end\sss we need\sss the following elementary\sss tool,\oss
a subspaces version of\dss convex combinations of\dss points.\oss
Let\vspace{1.5pt}
\[
\quad
H_{\dff 0}\off \supset\off 
H_{\dff 1}\off \supset\off 
\ldots\off \supset\off 
H_{\dff n}
\]

\vspace{-12pt}\vspace{1.5pt}
be a finite sequence of\dss closed subspaces of\sss $H$\nnsp.\oss
Let\sss $p_{\dff i}\dff \colon\dff H\qff \ttoo\qff H_{\dff i}$\dss
be\sss the orthogonal\sss projection onto $H_{\dff i}$\nsp,\oss
and\dss let\sss
$q_{\dff i}\dff \colon\dff 
H
\qff \ttoo\qff 
H_{\dff i}\dff \ominus\dff H_{\dff i\dff +\dff 1}$\dss
be\sss the orthogonal\sss projection onto\sss
$H_{\dff i}\dff \ominus\dff H_{\dff i\dff +\dff 1}$\nsp.\oss
Clearly,\qss $p_{\dff i}\off =\off q_{\dff i}\qff +\qff p_{\dff i\dff +\dff 1}$\nsp.\oss
Let\sss
$t_{\trf 0}\fff,\pff t_{\trf 1}\fff,\pff \ldots\fff,\pff t_{\dff n}$\sss
be non-negative real\sss numbers such\sss that\vspace{1.5pt}\vspace{0.5pt}
\[
\quad
t_{\trf 0}\qff +\qff t_{\trf 1}\qff +\qff \ldots\qff +\qff t_{\dff n}
\off =\off
1
\qff,
\]

\vspace{-12pt}\vspace{1.5pt}\vspace{0.5pt}
and\sss let\sss 
$s_{\dff i}
\off =\off
t_{\trf 0}\qff +\qff t_{\trf 1}\qff +\qff \ldots\qff +\qff t_{\dff i}$\nsp.\oss
The convex combination of\dss 
$p_{\dff i}$\sss 
with coefficients $t_{\dff i}$\sss is\vspace{3pt}
\[
\quad
t_{\trf 0}\dff p_{\dff 0}\qff +\qff 
t_{\trf 1}\dff p_{\dff 1}\qff +\qff 
\ldots\qff +\qff 
t_{\dff n}\dff p_{\dff n}
\off =\off
s_{\trf 0}\dff q_{\dff 0}\qff +\qff 
s_{\trf 1}\dff q_{\dff 1}\qff +\qff 
\ldots\qff +\qff 
s_{\dff n\dff -\dff  1}\dff q_{\dff n\dff -\dff 1}\qff +\qff 
p_{\dff n}
\pff.
\]

\vspace{-12pt}\vspace{3pt}
If\qss 
$K\qff \subset\qff H$\sss
is\dss a subspace such\sss that\sss $p_{\dff n}$\sss is\dss injective on\sss $K$\nnsp,\oss
then\sss the above identity\sss implies\sss that\vspace{1.5pt}\vspace{0.5pt}
\[
\quad
t_{\trf 0}\dff p_{\dff 0}\qff +\qff 
t_{\trf 1}\dff p_{\dff 1}\qff +\qff 
\ldots\qff +\qff 
t_{\dff n}\dff p_{\dff n}
\]

\vspace{-12pt}\vspace{1.5pt}\vspace{0.5pt}
is\dss also injective on\sss $K$\nnsp.\oss
The\sss image\vspace{1.5pt}\vspace{1pt}
\[
\quad
L
\off =\off
\left(\trf
t_{\trf 0}\dff p_{\dff 0}\qff +\qff 
t_{\trf 1}\dff p_{\dff 1}\qff +\qff 
\ldots\qff +\qff 
t_{\dff n}\dff p_{\dff n}
\trf\right)\dff (\trf K\trf)
\pff
\]

\vspace{-12pt}\vspace{1.5pt}\vspace{1pt}
is\dss contained\sss in\sss $H_{\dff 0}$\sss and
can\sss be\sss thought\sss as\sss the\qss \emph{convex\sss combination}\pss\vspace{1.5pt}\vspace{1pt}
\[
\quad
L
\off =\off
t_{\trf 0}\dff p_{\dff 0}\dff(\trf K\trf)\qff +\qff 
t_{\trf 1}\dff p_{\dff 1}\dff(\trf K\trf)\qff +\qff 
\ldots\qff +\qff 
t_{\dff n}\dff p_{\dff n}\dff(\trf K\trf)
\]

\vspace{-12pt}\vspace{1.5pt}\vspace{1pt}
of\dss subspaces\sss $p_{\dff i}\dff(\trf K\trf)$\nnsp.\oss
Let\sss 
$p_{\trf K}\dff \colon\dff H\qff \ttoo\qff K$\sss 
be\sss the orthogonal\sss projection.\oss
The maps\sss $p_{\trf K}$\sss and\vspace{1.5pt}\vspace{1pt}
\[
\quad
\left(\trf
t_{\trf 0}\dff p_{\dff 0}\qff +\qff 
t_{\trf 1}\dff p_{\dff 1}\qff +\qff 
\ldots\qff +\qff 
t_{\dff n}\dff p_{\dff n}
\trf\right)
\dff \circ\dff
p_{\dff K}
\]

\vspace{-12pt}\vspace{1.5pt}\vspace{1pt}
are connected\sss by\sss the standard\sss linear homotopy,\oss
and\dss this\sss homotopy\sss defines a canonical\sss path of\dss subspaces
connecting\sss $K$ and\dss $L$\nnsp.\oss

\mypar{Theorem.}{two-types-sections}
\emph{Suppose\sss that\trs $X$ is\dss a paracompact\sss space
and\dss
$\mathbb{A}$\sss
is\dss a\sss discrete-spectrum\sss and\sss strictly\trs Fredholm\dss family.\oss
Then every\sss weak\sss spectral\sss section\sss $\mathbb{S}$\sss of\qss
$\mathbb{A}$\sss
is\dss homotopic\sss to
a spectral\sss section\sss of\qss $\mathbb{A}$\sss
in\sss the class of\dss sections of\trs the bundle\dss}\vspace{1.5pt}
\[
\quad
\bm{\pi}\trf(\trf \mathbb{A}\trf)\dff \colon\dff
\gr\trf(\trf \mathbb{A}\trf)
\qff \ttoo\qff 
X\pff.
\]

\vspace{-12pt}\vspace{1.5pt}
\emph{If\pss the restriction of\dss the section\sss $\mathbb{S}$\sss to a subspace\dss
$Y\qff \subset\qff X$\sss is\dss a spectral\sss section
of\dss the family\sss $A_{\dff y}\dff,\qff y\qff \in\qff Y$\dnsp,\oss
then\sss the homotopy can\sss be chosen\sss in\sss the class of\dss
sections\sss with\sss this property.\oss}

\proof
Let\sss us\sss write\sss the section\sss $\mathbb{S}$\sss as a family\sss
$S_{\dff x}\dff,\qff x\qff \in\qff X$\nnsp.\oss
Let\sss us define\sss the subspaces\vspace{1.5pt}
\[
\quad
H_{\qff <\dff k}\trf(\dff x\trf)
\off =\off
\bigoplus\nolimits_{\dff n\qff <\qff k}\qff H_{\dff n}\trf(\trf x\trf)
\]

\vspace{-12pt}\vspace{1.5pt}
in\sss terms of\dss the decomposition\qss (\ref{spectral-decomposition}).\oss
Let\sss us\sss temporarily\sss fix a point\sss $z\qff \in\qff X$\nnsp.\oss
The subspace\sss $S_{\fff z}$\sss is\dss commensurable with\sss
$H_{\qff \geq\dff 0}\trf(\trf z\trf)$\sss
and\sss hence\sss the intersection\sss
$S_{\fff z}\dff \cap\dff H_{\qff <\dff 0}\trf(\trf z\trf)$\sss
is\dss finitely dimensional.\oss
By\sss a compactness argument\sss this 
implies\sss that\vspace{1.5pt}
\[
\quad
\bigl(\trf
S_{\fff z}\dff \cap\dff H_{\qff <\dff 0}\trf(\trf z\trf)
\trf\bigr)
\qff \cap\qff
H_{\qff <\dff -\qff n}\dff(\trf z\trf)
\off =\off
0
\]

\vspace{-12pt}\vspace{1.5pt}
for sufficiently\sss large $n$\nnsp.\oss
See\sss the beginning of\dss the proof\dss of\trs 
Lemma\qss 11.1\qss in\qss \cite{i2}\qss for\sss the details.\oss
It\sss follows\sss that\sss
$S_{\dff z}\dff \cap\dff H_{\qff <\dff -\qff n}\trf(\trf z\trf)
\off =\off
0$\sss
for sufficiently\sss large $n$\nnsp.\oss
By\sss the definition,\vspace{1.5pt}
\[
\quad
H_{\qff <\dff -\qff n}\trf(\trf z\trf)
\off =\off
\image
P_{\qff \leq\qff \nu\trf(\dff z\trf)}\qff(\trf A_{\dff z}\trf)
\off =\off
\kernel
P_{\qff \geq\qff \nu\trf(\trf z\trf)}\qff(\trf A_{\dff z}\trf)
\]

\vspace{-12pt}\vspace{1.5pt}
for some\sss $\nu\trf(\trf z\trf)\qff \in\qff \rrr$\sss
different\sss from\sss the eigenvalues of\dss $A_{\dff x}$\nsp.\oss
Therefore\sss the restriction\sss\vspace{1.5pt}
\[
\quad
P_{\qff \geq\qff \nu\trf(\trf z\trf)}\qff(\trf A_{\dff z}\trf)
\qff \bigl|\pff S_{\fff z}
\]

\vspace{-12pt}\vspace{1.5pt}
is\dss injective.\oss
We claim\sss that,\oss moreover,\vspace{1.5pt}
\[
\quad
P_{\qff \geq\qff \nu\trf(\trf z\trf)}\qff(\trf A_{\dff y}\trf)
\qff \bigl|\pff S_{\fff y}
\]

\vspace{-12pt}\vspace{1.5pt}
is\dss injective for every\sss $y$\sss in\sss 
some neighborhood\sss $U$\sss of\sss $z$\nnsp.\oss
In\sss order\sss to prove\sss this,\oss let\sss us choose a pair\sss
$(\trf U\fff,\qff \varepsilon\trf)$\sss
adapted\dss to\sss $\mathbb{A}$\sss and such\sss that\sss $z\qff \in\qff U$\sss
and a\sss local\dss trivialization\sss $t_{\trf U}$ of\dss 
$\mathbb{H}$\sss defined over $U$\sss and strictly\sss adapted\sss to\sss
$(\trf U\fff,\qff \varepsilon\trf)$\sss
and\sss $\mathbb{A}$\nnsp.\oss
Then\sss the families\sss
$t_{\trf U}\trf(\trf S_{\dff y}\trf)$\nnsp,\qss $y\qff \in\qff U$\sss
and\vspace{1.5pt}
\[
\quad
t_{\trf U}\qff
\bigl(\qff
\image
P_{\qff \geq\qff \nu\trf(\trf z\trf)}\qff(\trf A_{\dff y}\trf)
\qff\bigr)\dff,\off
y\qff \in\qff U
\]

\vspace{-12pt}\vspace{1.5pt}
are norm continuous families of\dss subspaces of\sss $H$\nnsp.\oss
This reduces our claim\sss to\sss the following one.\oss
Let\sss $T_{\dff y}\dff,\qff y\qff \in\qff U$\sss be a norm continuous
family of\dss closed subspaces of\sss $H$\nnsp,\oss
and\sss let\sss $P_{\dff y}\dff,\qff y\qff \in\qff U$\sss
be a norm continuous family\sss of\dss projections.\oss
If\dss the restriction\sss $P_{\dff z}\qff |\qff T_{\dff z}$\sss
is\dss injective and\sss 
$P_{\dff z}\trf(\trf T_{\dff z}\trf)$\sss is\dss closed,\oss 
then\sss $P_{\dff y}\qff |\qff T_{\dff y}$\sss
is\dss injective for all\sss $y$ sufficiently close\sss to $z$\nnsp.\oss
In order\sss to prove\sss this,\oss
consider\sss the orthogonal\sss complement\sss $F$\sss of\dss 
$P_{\dff z}\trf(\trf T_{\dff z}\trf)$\sss
in\dss $\image P_{\dff z}$\nsp.\oss
Then\sss the map\sss
$F\dff \oplus\dff T_{\dff z}
\qff \ttoo\qff
\image P_{\dff z}$\sss
equal\sss to\sss the identity on\sss $F$\sss and\sss to\dss $P_{\dff z}$\sss
on\sss $T_{\dff z}$\sss
is\dss an\sss isomorphism.\oss
By\sss the open\sss mapping\sss theorem\sss  
similar\sss maps\sss
$F\dff \oplus\dff T_{\dff y}
\qff \ttoo\qff
\image P_{\dff y}$\sss
are isomorphisms for all\sss $y$ sufficiently close\sss to $z$\nnsp.\oss
Therefore\dss $P_{\dff y}\qff |\qff T_{\dff y}$\sss
is\dss injective for all\sss $y$ sufficiently close\sss to $z$\nnsp.\oss
This proves our claim.\oss

Since $X$\sss is\dss paracompact,\oss there exist\sss a\sss locally\sss finite open covering\sss
$G_{\dff i}\dff,\qff i\qff \in\pff I$\sss of\dss $X$\sss
refining\sss the covering\sss $U_{\dff z}\fff,\qff z\qff \in\qff X$\sss
and a partition of\dss unity\sss
$t_{\trf i}\dff,\qff i\qff \in\pff I$\sss subordinated\sss to\sss this covering.\oss
Since\sss $G_{\dff i}\dff,\qff i\qff \in\pff I$\sss 
refines\sss $U_{\dff z}\fff,\qff z\qff \in\qff X$\nnsp,\oss
for every\sss $i\qff \in\pff I$\dss there\dss is\dss a number\sss
$\nu\trf(\dff i\trf)\qff \in\qff \rrr$\sss such\sss that\sss $\nu\trf(\dff i\trf)$\sss
is\dss not\sss an eigenvalue of\dss $A_{\dff x}$\sss and\vspace{1.5pt}
\[
\quad
P_{\qff \geq\qff \nu\trf(\trf z\trf)}\qff(\trf A_{\dff x}\trf)
\qff \bigl|\pff S_{\fff x}
\]

\vspace{-12pt}\vspace{1.5pt}
is\dss injective\sss
for\sss every\sss
$x\qff \in\qff G_{\dff i}$\nsp.\oss
Let\sss us\sss consider\sss the convex combinations\vspace{2.5pt}
\[
\quad
L_{\trf x}
\off =\off
\sum\nolimits_{\qff i\qff \in\qff I\vphantom{I^i}}\qff 
t_{\trf i}\trf(\dff x\trf)\qff
\left(\trf
P_{\qff \geq\qff \nu\trf(\dff i\trf)}\qff(\trf A_{\dff x}\trf)
\trf\right)\dff
(\trf S_{\fff x}\trf)
\] 

\vspace{-12pt}\vspace{1.5pt}
interpreted\sss as explained\sss before\sss the\sss theorem.\oss
These convex combinations make sense because for every\sss $x$\sss
only a finite number of\dss coefficients\sss
$t_{\dff i}\dff(\dff x\trf)$\sss
is\dss non-zero,\oss and\sss the images\dss
$\image 
P_{\qff \geq\qff \nu\trf(\dff i\trf)}\qff(\trf A_{\dff x}\trf)$\sss
are\sss
linearly\sss ordered\sss by\sss the inclusion.\oss
It\sss follows\sss that\vspace{1.5pt}\vspace{-0.5pt}
\begin{equation}
\label{partial-mp}
\quad
L_{\trf x}
\off \subset\off
\image
P_{\qff \geq\qff \mu\trf(\trf x\trf)}\qff(\trf A_{\dff x}\trf)
\pff,
\end{equation} 

\vspace{-12pt}\vspace{1.5pt}\vspace{-0.5pt}
where\sss $\mu\trf(\dff x\trf)$\sss is\dss the minimum of\dss
the numbers\sss $\nu\trf(\dff i\trf)$\sss over\sss $i$\sss such\sss that\sss
$t_{\dff i}\trf(\dff x\trf)\off \neq\off 0$\nnsp.\oss
Moreover,\pss $S_{\fff x}$\sss is\dss connected\sss with\sss
$L_{\dff x}$\sss by a canonical\sss path of\dss subspaces.\oss
These canonical\sss paths form a homotopy\sss between\sss the sections\sss
$x\off \longmapsto\off S_{\fff x}$\sss and\dss 
$x\off \longmapsto\off L_{\trf x}$\nnsp.\oss

We would\sss like\sss the inclusions\qss (\ref{partial-mp})\qss to hold\sss for
some continuous function\sss $\mu\trf(\trf x\trf)$\nnsp.\oss 
This can\sss be done\dss if\dss we carefully choose\sss the partition of\dss unity.\oss
First,\oss given a\sss locally\sss finite open covering\sss
$G_{\dff i}\dff,\qff i\qff \in\pff I$\sss as above,\oss
there exists a closed covering\sss $F_{\dff i}\dff,\qff i\qff \in\pff I$\sss
such\sss that\sss $F_{\dff i}\qff \subset\qff G_{\dff i}$\sss for every $i$\nnsp.\oss
See,\oss for example,\oss \cite{i1},\oss Theorem\qss 7.5.\oss
Since $X$\nnsp,\oss being a paracompact\dss Hausdorff\dss space,\oss
is\dss normal,\sss there exist\sss open sets $U_{\dff i}\dff,\qff i\qff \in\pff I$\sss
such\sss that\sss $F_{\dff i}\qff \subset\qff U_{\dff i}$\sss and\dss
$\overline{U}_{\dff i}\qff \subset\qff G_{\dff i}$\sss for every $i$\nnsp,\oss
where\sss $\overline{U}_{\dff i}$\sss is\dss the closure of\sss $U_{\dff i}$\nsp.\oss
Using\sss the fact\sss that\sss $X$\sss is\dss normal\sss once more,\oss
we see\sss that\dss there exist\sss continuous functions\sss
$s_{\dff i}\dff \colon\dff X\qff \ttoo\qff \rrr$\sss
equal\sss to $1$ on\sss $\overline{U}_{\dff i}$\sss and\sss to $0$ on\sss
$X\qff \smallsetminus\qff G_{\dff i}$\nsp.\oss
We may assume\sss that\sss the partition of\dss unity\sss
$t_{\dff i}\dff,\qff i\qff \in\pff I$\sss is\dss subordinated\sss to\sss
the covering\sss $U_{\dff i}\dff,\qff i\qff \in\pff I$\nnsp.\oss
Then\sss $s_{\dff i}\dff(\dff x\trf)\off =\off 1$\sss if\dss
$t_{\dff i}\dff(\dff x\trf)\off \neq\off 0$\nnsp.\oss
Let\sss us\sss redefine\sss $\mu\dff(\dff x\trf)$ as\vspace{1.5pt}
\[
\quad
\mu\trf(\trf x\trf)
\off =\off\dff
\sum\nolimits_{\qff i\qff \in\qff I\vphantom{I^i}}\qff 
s_{\dff i}\trf(\trf x\trf)\qff \nu\trf(\dff i\trf)
\pff.
\]

\vspace{-12pt}\vspace{1.5pt}
Since\sss this\dss is\dss a\sss locally\sss finite sum and\sss $\nu\trf(\dff i\trf)$
are constant,\sss the\sss function\sss $\mu\trf(\trf x\trf)$\sss is\dss continuous.\oss
At\sss the same\sss time\sss $\mu\trf(\trf x\trf)$\sss is\dss greater or equal\sss
than\sss the sum of\dss of\dss
the numbers\sss $\nu\trf(\trf i\trf)$\sss over\sss $i$\sss such\sss that\sss
$t_{\dff i}\trf(\trf x\trf)\off \neq\off 0$\nnsp.\oss 
Hence\sss $\mu\trf(\trf x\trf)$\sss is\dss greater or equal\sss
than\sss the maximum of\dss these numbers.\oss
It\sss follows\sss that\qss (\ref{partial-mp})\qss holds\sss for\dss this $\mu$\nnsp.\oss

Let\sss
$L^{\dff \perp}_{\trf x}
\off =\off
H\dff \ominus\dff L_{\trf x}$\nsp.\oss
Clearly,\qss $L^{\dff \perp}_{\trf x}$\sss
is\dss commensurable with\sss
$P_{\dff \leq\dff 0}\trf (\dff A_{\dff x}\dff)$\sss
and\dss\vspace{1.5pt}\vspace{-0.5pt}
\begin{equation}
\label{half-proof}
\quad
L^{\dff \perp}_{\trf x}
\off \supset\off
\image
P_{\qff <\pff \mu\trf(\dff x\trf)}\qff(\dff A_{\dff x}\dff)
\pff.
\end{equation}

\vspace{-12pt}\vspace{1.5pt}\vspace{-0.5pt}
By\sss interchanging\sss the roles of\dss positive and\sss
negative numbers in\sss the above construction,\oss
we can deform\sss the map\sss
$L^{\dff \perp}$\sss
to a map\sss
$M^{\dff \perp}$\sss
such\sss that\vspace{1.5pt}\vspace{-0.5pt}
\[
\quad
M^{\dff \perp}_{\trf x}
\off \subset\off
\image
P_{\qff \leq\pff \mu^{\dff \perp}\trf(\dff x\trf)}\qff(\trf A_{\dff x}\trf)
\pff
\] 

\vspace{-12pt}\vspace{1.5pt}\vspace{-0.5pt}
for some continuos function\sss $\mu^{\dff \perp}\trf(\dff x\trf)$\nnsp.\oss
Moreover,\oss it\sss is\dss clear\sss from\sss the construction\sss that\sss
the inclusion\qss (\ref{half-proof})\qss holds\sss with\sss 
$M^{\dff \perp}_{\trf x}$\sss
in\sss the role of\sss $L^{\dff \perp}_{\trf x}$\nsp.\oss
By\sss passing from\sss the subspaces\sss $M^{\dff \perp}_{\trf x}$\sss
to\sss their orthogonal\sss complements\sss $M_{\trf x}$\sss
we get\sss a section\sss $M_{\trf x}\dff,\qff x\qff \in\qff X$\sss 
homotopic\sss to\sss the section\sss $S_{\fff x}\dff,\qff x\qff \in\qff X$\sss 
and such\sss that\dss 
$P_{\qff >\pff \mu^{\dff \perp}\trf(\dff x\trf)}\trf (\trf A_{\dff x}\trf)
\off \subset\off
M_{\trf x}
\off \subset\off
\image
P_{\qff \geq\qff \mu\trf(\dff x\trf)}\qff(\trf A_{\dff x}\trf)$\nnsp.\oss
It\sss follows\sss that\sss the family\sss $M_{\trf x}\dff,\qff x\qff \in\qff X$\sss
is\dss a spectral\sss section\sss of\dss
$A_{\dff x}\dff,\qff x\qff \in\qff X$\nnsp.\oss
This proves\dss the first\sss claim of\dss the\sss theorem.\oss
In order\sss to prove\sss the second claim,\oss
it\dss is\dss sufficient\sss to observe\sss that\sss if\dss the subspaces\sss
$S_{\fff x}$\sss satisfy\qss (\ref{mp-section})\qss for\sss $x\qff \in\qff Y$\nnsp,\oss
then\sss this\sss remains\sss the case during\sss the whole deformation\sss
from\sss $S_{\fff x}$\sss to\sss $M_{\trf x}$\nsp,\oss
except\sss that\sss one may\sss need\sss to replace\sss the function\sss $r\trf(\trf x\trf)$\sss 
by a\sss larger one.\oss  \eproof

\mypar{Corollary.}{index-and-mp-sections}
\emph{Suppose\sss that\trs $X$ is\dss paracompact\sss
and\dss $\mathbb{A}$\sss is\dss a\sss discrete-spectrum\sss 
and\sss strictly\trs Fredholm\dss family.\oss
Then\sss the analytical\dss index\sss of\qss $\mathbb{A}$ vanishes\dss if\trs and\dss only\trs if\qss
there\sss exists\sss a spectral\sss section\sss of\qss $\mathbb{A}$\nnsp.\oss}\vspace{-1.125pt}

\proof
By\trs Theorem\qss \ref{weak-spectral-sections},\oss
if\dss the analytical\dss index\sss of\dss $\mathbb{A}$ vanishes,\oss
then\sss there exists a weak\sss spectral\sss section\sss for\sss $\mathbb{A}$\nnsp.\oss
Since\sss $X$\sss is\dss paracompact,\oss
Theorem\qss \ref{two-types-sections}\qss implies\sss that\sss in\sss this case\sss there
exists a spectral\sss section\sss of\sss $\mathbb{A}$\nnsp.\oss
The other\sss implication\sss follows\sss from\dss
Theorem\qss \ref{weak-spectral-sections}.\oss  \eproof

\myuppar{Remark.}
The really\sss important\sss part\dss is\dss the\qss ``only\trs if''\qss one.\oss
It\dss is\dss independent\sss from\sss the\qss ``if''\qss part\sss of\trs
Theorem\qss \ref{weak-spectral-sections},\oss and\sss hence\dss
from\sss the\sss theorem\qss \cite{i2}\qss 
about\sss the contractibility of\dss $\mathbf{G}$\nnsp.\oss

\myuppar{Remark.}
The proof\dss of\qss Theorem\qss \ref{two-types-sections}\qss
resulted\sss from an attempt\sss to understand\dss the\sss last\sss 
paragraph\sss in\sss the proof\dss of\qss Proposition\trs 1\trs
of\pss Melrose\dss and\trs Piazza\qss \cite{mp}.\oss
Melrose\dss and\trs Piazza\qss \cite{mp}\qss
claim\sss that\sss a weak spectral\sss section can\sss be\sss
transformed\sss into a spectral\sss section\qss
\emph{``simply\dss by\sss smoothly\dss truncating\sss
the eigenfunction\sss expansion''.}\oss
Since\sss the eigenvalues change with $x$\nnsp,\oss it\dss is\dss only\sss
natural\sss to\sss truncate expansions\qss \emph{between}\pss the eigenvalues.\oss
This can\sss be always done\sss locally,\oss
but\sss if\dss we can do\sss this globally,\oss
then we already\sss know\sss that\sss there exists 
a spectral\sss section.\oss
The proof\dss of\qss Theorem\qss \ref{two-types-sections}\qss
provides a natural\sss way\sss to overcome\sss this difficulty.\oss

\myuppar{A\sss theorem of\qss Melrose--Piazza.}
The classical\sss analytical\dss index\dss is\dss defined\sss only\sss under\sss
much stronger assumptions\sss than\sss that\sss of\trs Corollary\qss \ref{index-and-mp-sections}.\oss
As we will\sss see in\dss Section\qss \ref{classical-analytic-index-sa},\oss
which\dss is\dss independent\sss of\dss the present\sss one,\oss
when\sss the classical\sss analytical\dss index\dss is\dss defined,\oss
it\sss agrees with\sss the analytical\dss index as defined\sss in\dss
Section\qss \ref{analytic-index-section}.\oss
See\dss Theorem\qss \ref{sa-index-agree}.\oss
For such families\sss the analytical\dss index\sss in\dss
Corollary\qss \ref{index-and-mp-sections}\qss can\sss be understood\sss
in\sss the classical\sss sense,\oss
and\sss this\sss turns\sss its\sss conclusion\sss 
into\sss the conclusion of\trs Proposition\trs 1\trs
of\pss Melrose\dss and\trs Piazza\qss \cite{mp}.\oss
The assumptions are still\sss weaker\sss than\sss the assumptions
of\qss Melrose\dss and\trs Piazza\qss \cite{mp},\oss
who required\sss that\sss the space\sss $X$\sss to be compact\sss
and\sss $\mathbb{A}$\sss to be a continuous family of\dss 
differential\sss operators of\dss oder $1$\nnsp.\oss
The\sss last\sss assumption\dss is\dss not\sss really\sss used\sss in\dss 
Melrose--Piazza\dss proof,\oss
but\sss the proof\dss depends on continuity\sss properties
of\dss the family\sss which are much stronger\sss 
than\sss being\sss strictly\trs Fredholm.\oss

Melrose\dss and\trs Piazza\qss \cite{mp}\qss
gloss over\sss the definition of\dss analytical\dss index and\sss the meaning\sss
of\dss trivializations of\trs Hilbert\dss bundles,\oss 
and a straightforward\sss interpretation of\dss their\sss proof\dss works only\sss
for\sss operators in a fixed\dss Hilbert\dss space.\oss
The author\sss learned about\sss this issue\sss from\dss M.\dss Prokhorova\qss \cite{p1}.\oss
The above proof\dss overcomes\sss these difficulties\sss
thanks\dss to\trs Theorem\qss \ref{adapted-paracompact}\qss and\sss 
the new definition of\dss analytical\dss index.\oss
There are also a\sss less\sss general\sss versions of\qss Melrose--Piazza\trs theorem,\oss
which are still\sss more\sss general\sss than\sss the original\sss one,\pss
and\dss which can\sss be proved\sss using only\trs Theorem\qss
\ref{adapted}\qss or\qss \ref{adapted-paracompact}\qss and\sss
the classical\sss definition of\dss analytical\dss index.\oss
See\qss \cite{i3}.\oss

\mysection{The\qss classical\qss index\qss of\qss Fredholm\qss families}{classical-analytic-index}

\myuppar{Two definitions of\dss the analytical\dss index.}
The goal\sss of\dss this section\dss is\dss to prove\sss that\sss
the definition of\dss analytical\dss index\sss of\trs Fredholm\dss families 
from\dss Section\qss \ref{analytic-index-fredholm-non-sa}\qss is\dss equivalent\sss
to\sss the classical\sss one when\sss the\sss latter\dss applies.\oss
Let\sss $B_{\dff x}\dff \colon\dff H_{\dff x}\qff \ttoo\qff H_{\dff x}$\nsp,\qss
$x\qff \in\qff X$\nnsp,\oss
be such a family of\dss operators acting in\sss the fibers of\dss
a\sss Hilbert\dss bundle\sss $H_{\dff x}\dff,\qff x\qff \in\qff X$\nnsp.\oss
The most\sss classical\sss definition of\trs
Atiyah--Singer\qss \cite{as4}\qss applies only\sss
when\sss $X$\sss is\dss compact,\oss and we will\sss discuss\sss this case first.\oss\vspace{1pt}

The classical\sss analytical\dss index of\dss the family\sss 
$B_{\dff x}\dff,\qff x\qff \in\qff X$\nnsp,\oss
is\dss an element\sss of\dss the group\sss $K\trf(\trf X\trf)$\nnsp,\oss
which can be defined either in\sss terms of\dss finitely dimensional\sss
vector bundles on\sss $X$\nnsp,\oss or\sss as\sss the group of\dss
homotopy classes of\dss maps\sss to\sss
a classifying space for\sss $K$\dnsp-theory.\oss
The space\sss $\mathcal{F}$\sss of\trs Fredholm\dss operators\dss
is\dss a well\dss known such classifying space.\oss
Since\sss there are canonical\sss homotopy equivalences\sss
$\mathcal{F}\qff \ttoo\qff \num{\mathcal{S}}\qff \ttoo\qff \num{S}$\nnsp,\oss
we can\sss interpret\sss our analytical\dss index as a map\sss to\sss $\mathcal{F}$\dnsp.\oss
This almost\sss turns\sss the problem of\dss comparing\sss two definition\sss
into a\sss triviality\sss when\sss the\dss Hilbert\dss bundle\sss 
$H_{\dff x}\dff,\qff x\qff \in\qff X$\nnsp,\oss is\dss trivial.\oss
But\sss there\dss is\dss something\sss to discuss\sss
in\sss the general\sss case.\oss

\myuppar{The classical\sss analytical\dss index of\trs Fredholm\dss families.}
As usual,\oss we will\sss denote\sss the\dss Hilbert\dss bundle\sss
$H_{\dff x}\dff,\qff x\qff \in\qff X$\nnsp,\oss
and\sss the family\sss
$B_{\dff x}\dff,\qff x\qff \in\qff X$\nnsp,\oss
also by\sss $\mathbb{H}$\sss and\sss $\mathbb{B}$\sss respectively.\oss 
Classically,\oss the analytical\dss index of\sss $\mathbb{B}$\sss
is\dss defined\sss only\sss under\sss very strong assumptions about\sss
$\mathbb{B}$\nnsp.\oss
See\qss \cite{as4}.\oss
But\sss the classical\sss definition works\sss without\sss 
any changes also for general\trs Fredholm\dss families.\oss\vspace{1pt}

Suppose\sss that\sss $X$\sss is\dss compact.\oss
Then a classical\sss argument\sss due\sss to\dss Atiyah\dss shows\sss that\sss there exists a
finitely dimensional\dss Hilbert\sss space\sss $V$\sss
and a continuous family\sss of\dss maps\sss
$q_{\dff x}\dff \colon\dff V\qff \ttoo\qff H_{\dff x}$\nsp,\qss $x\qff \in\qff X$\nnsp,\oss
such\sss that\sss the map\sss
$Q_{\dff x}\dff \colon\dff
V\dff \oplus\trf H_{\dff x}\qff \ttoo\qff H_{\dff x}$\sss
defined\sss by\sss the formula\sss\vspace{3pt}
\[
\quad
Q_{\dff x}\dff(\trf v\dff \oplus\dff u\trf)
\off =\off
q_{\dff x}\trf(\dff v\trf)\qff +\qff B_{\dff x}\trf(\dff u\trf)
\]

\vspace{-12pt}\vspace{3pt}
is\dss surjective for every\sss $x\qff \in\qff X$\nnsp.\oss
Moreover,\oss the standard construction\sss leads\sss to a family\sss of\dss maps\sss
$q_{\dff x}$\nsp,\qss $x\qff \in\qff X$\nnsp,\oss having\sss the following\sss property.\oss
For every\sss $x\qff \in\qff X$\sss
there exists a neighborhood\sss $U_{\fff x}$\sss of\sss $x$\sss
and\sss $\varepsilon\off =\off \varepsilon_{\dff x}\qff >\qff 0$\sss
such\sss that\sss for\sss $y\qff \in\qff U_{\dff x}$\sss the subspaces\dss 
$\image\dff P_{\qff [\dff 0\fff,\qff \varepsilon\dff]}\trf(\trf \num{B^*_{\dff y}}\trf)$\dss
are finitely dimensional,\oss continuously depend on $y$\nnsp,\oss
and contain\sss the image\sss $\image\dff q_{\dff x}$\nsp.\oss
Then\sss the family of\dss finitely dimensional\sss vector spaces\sss
$\kernel\dff Q_{\dff x}\dff,\qff x\qff \in\qff X$\nnsp,\oss
is\dss a finitely dimensional\sss vector bundle,\oss
which we denote by\sss $\mathbb{K}$\nnsp.\oss
Let\sss $\mathbb{V}$\sss be\sss the\sss trivial\sss bundle\sss
$V\dff \times\dff X\qff \ttoo\qff X$\nnsp.\oss
The\qss \emph{classical\sss analytical\dss index}\pss of\dss $\mathbb{B}$\dss
is\dss defined as\sss the element\sss of\dss $K\trf(\trf X\trf)$\sss
represented\sss by\sss the difference\sss
$[\dff \mathbb{K}\dff]\qff -\qff [\dff \mathbb{V}\dff]$\nnsp.\oss\vspace{1pt}

\mypar{Theorem.}{fredholm-agree}
\emph{If\qss $X$ is\dss compact,\oss 
then\dss two definitions of\dss the analytical\dss index of\pss $\mathbb{B}$\sss agree.\oss}

\proof
Let\sss
$P_{\dff x}\dff \colon\dff
V\dff \oplus\trf H_{\dff x}
\qff \ttoo\qff
V\dff \oplus\trf H_{\dff x}$\sss
be\sss the composition of\sss $Q_{\dff x}$\sss with\sss the inclusion of\sss
$H_{\dff x}$\sss into\sss $V\dff \oplus\trf H_{\dff x}$\sss
as\sss the second summand.\oss
The family\sss
$P_{\dff x}\dff,\pff x\qff \in\qff X$\nnsp,\oss
can be considered as a bundle map\sss
$P\dff \colon\dff
\mathbb{V}\dff \oplus\dff \mathbb{H}
\qff \ttoo\qff
\mathbb{V}\dff \oplus\dff \mathbb{H}$\nnsp.\oss
Clearly,\oss its cokernel\dss is\sss $\mathbb{V}$\dnsp,\oss
and,\oss by\sss the definition,\oss its kernel\dss is\sss $\mathbb{K}$\nnsp.\oss
At\sss the same\sss time\sss $P_{\dff x}\dff,\pff x\qff \in\qff X$\nnsp,\oss
is\dss a\dss Fredholm\dss family.\oss

Let\sss us define a family of\dss operators\sss
$R_{\dff x}\dff \colon\dff
V\dff \oplus\trf H_{\dff x}
\qff \ttoo\qff 
V\dff \oplus\trf H_{\dff x}$\nsp,\qss
$x\qff \in\qff X$\nnsp,\oss
by\sss the formula\sss\vspace{1.5pt}
\[
\quad
R_{\dff x}\trf(\trf v\dff \oplus\dff u\trf)
\off =\off
v\dff \oplus\trf B_{\dff x}\trf(\dff u\trf)
\qff.
\]

\vspace{-12pt}\vspace{1.5pt}
Clearly,\oss the family\sss
$R_{\dff x}\dff,\qff x\qff \in\qff X$\nnsp,\oss
is\dss also a\dss Fredholm\dss family.\oss
The index\sss maps\sss to\sss 
$\num{\mathcal{S}\dff(\trf \mathbb{V}\dff \oplus\dff \mathbb{H}\trf)}$\sss 
of\dss the families\sss
$R_{\dff x}\dff,\qff x\qff \in\qff X$\nnsp,\oss
and\sss 
$B_{\dff x}\dff,\off x\qff \in\qff X$\nnsp,\oss
can\sss be arranged\sss to be not\sss only homotopic,\oss but\sss equal.\oss
At\sss the same\sss time\sss the family\sss
$R_{\dff x}\dff,\qff x\qff \in\qff X$\nnsp,\oss
is\dss obviously\trs Fredholm\dss homotopic\sss to\sss the family\sss
$P_{\dff x}\dff,\pff x\qff \in\qff X$\nnsp,\oss
and\sss hence\sss their\sss index maps are homotopic.\oss
Therefore\sss it\dss is\dss sufficient\sss to prove\sss that\sss
the index\sss map of\dss the family\sss
$P_{\dff x}\dff,\pff x\qff \in\qff X$\nnsp,\oss
corresponds\sss to\sss
$[\dff \mathbb{K}\dff]\qff -\qff [\dff \mathbb{V}\dff]$\nnsp.\oss

Let\sss $x\qff \in\qff X$\nnsp.\oss
Since\sss 
$\kernel\dff P_{\dff x}^{\fff *}\qff P_{\dff x}
\off =\off
\kernel\dff P_{\dff x}$\sss
and\sss $P_{\dff x}$\sss is\dss Fredholm,\oss
there exists\sss $\varepsilon\qff >\qff 0$\sss
such\sss that\vspace{1.5pt}
\[
\quad
\kernel\dff Q_{\dff x}
\off =\off
\kernel\dff P_{\dff x}
\off =\off
\image\dff P_{\qff [\dff 0\fff,\qff \varepsilon\dff]}\qff
\left(\qff \bnum{P_{\dff x}}\qff\right)
\qff.
\]

\vspace{-12pt}\vspace{1.5pt}
Moreover,\oss we can assume\sss that\sss
$(\trf P_{\dff x}\dff,\qff \varepsilon\trf)$\sss
is\dss an enhanced operator.\oss
Clearly,\oss\vspace{1.5pt}
\[
\quad
\kernel\dff Q_{\dff z}
\off =\off
\kernel\dff P_{\dff z}
\off \subset\off
\image\dff P_{\qff [\dff 0\fff,\qff \varepsilon\dff]}\qff
\left(\qff \bnum{P_{\dff z}}\qff\right)
\qff.
\]

\vspace{-12pt}\vspace{1.5pt}
for every\sss $z\qff \in\qff X$\nnsp.\oss
Since\sss $P_{\dff x}\dff,\pff x\qff \in\qff X$\nnsp,\oss
is\dss a\dss Fredholm\dss family,\oss
the spaces\sss
$\image\dff P_{\qff [\dff 0\fff,\qff \varepsilon\dff]}\qff
\left(\qff \bnum{P_{\dff z}}\qff\right)$\sss
continuously depend on\sss $z$\sss for\sss $z$\sss sufficiently
close\sss to $x$\nnsp.\oss
Since\sss the spaces\sss $\kernel\dff Q_{\dff z}$\sss
also continuously\sss depend on $z$\nnsp,\oss
it\sss follows\sss that\vspace{1.5pt}
\[
\quad
\kernel\dff Q_{\dff z}
\off =\off
\kernel\dff P_{\dff z}
\off =\off
\image\dff P_{\qff [\dff 0\fff,\qff \varepsilon\dff]}\qff
\left(\qff \bnum{P_{\dff z}}\qff\right)
\qff
\]

\vspace{-12pt}\vspace{1.5pt}
for\sss $z$\sss sufficiently close\sss to $x$\nnsp.\oss
Since\sss the space\sss $X$\sss is\dss compact,\oss
for sufficiently small\sss $\varepsilon\qff >\qff 0$\sss these equalities
will\dss hold\sss for every\sss $z\qff \in\qff X$\nnsp.\oss
Similarly,\oss for sufficiently small\sss $\varepsilon\qff >\qff 0$\vspace{1.5pt}
\[
\quad
V\dff \oplus\dff 0
\off =\off
\kernel\dff Q_{\dff z}^{\fff *}
\off =\off
\kernel\dff P_{\dff z}^{\fff *}
\off =\off
\image\dff P_{\qff [\dff 0\fff,\qff \varepsilon\dff]}\qff
\left(\qff \bnum{P_{\dff z}^{\fff *}}\qff\right)
\qff
\]

\vspace{-12pt}\vspace{1.5pt}
for every\sss $z\qff \in\qff X$\nnsp.\oss
Then\sss $(\trf Q_{\dff z}\dff,\qff \varepsilon\trf)$\sss
is\dss an enhanced operator\sss for every\sss $z\qff \in\qff X$\nnsp.\oss

Now we can construct\sss an\sss index\sss map for\sss
$P_{\dff x}\dff,\pff x\qff \in\qff X$\nnsp,\oss using\sss
the covering\sss $U$\sss of\sss $X$\sss by\sss single open set\sss $X$\sss
and\sss $\varepsilon$\sss as\sss the corresponding\sss number.\oss
Then\sss $\num{X_{\dff U}}\off =\off X$\sss and\sss the homotopy\sss inverse\sss
$s\dff \colon\dff X\qff \ttoo\qff \num{X_{\dff U}}$\sss is\dss
the identity\sss map.\oss
The index\sss map\sss
$X\qff \ttoo\qff \num{\mathcal{S}\dff(\trf \mathbb{H}\trf)}$\sss
takes\sss point\sss $x\qff \in\qff X$\sss to\sss the point\sss of\sss
$\num{\mathcal{S}}$\sss represented\sss by\sss the object\vspace{1.5pt}
\[
\quad
\bigl(\qff
\kernel\dff P_{\dff x}\dff,\pff
V
\qff\bigr)
\off =\off
\bigl(\qff
\kernel\dff Q_{\dff x}\dff,\pff
V\fff \times\dff x
\qff\bigr)
\]

\vspace{-12pt}\vspace{1.5pt}
of\sss $\mathcal{S}\dff(\trf \mathbb{H}\trf)$\dnsp.\oss
Now we can either compose\sss this index\sss map\sss with\sss 
$\num{\mathcal{S}\dff(\trf \mathbb{H}\trf)}
\qff \ttoo\qff 
\num{S\trf(\trf \mathbb{H}\trf)}
\qff \ttoo\qff
\num{S}$\sss
to get\sss an\sss index\sss map\sss
$X\qff \ttoo\qff \num{S}$\nnsp,\oss
or use a\sss trivialization of\dss the bundle\sss $\mathbb{H}$\sss 
to get\sss an\sss index\sss map\sss
$X\qff \ttoo\qff \num{\mathcal{S}}$\nnsp.\oss
Either\sss way,\oss the resulting\sss index\sss map\sss
represents\sss the virtual\sss bundle\sss $\mathbb{K}\qff -\qff \mathbb{V}$\sss
at\sss least\sss on\sss the intuitive\sss level.\oss
In\sss fact,\oss it\dss is\dss easy\sss to construct\sss a map\sss 
$X\qff \ttoo\qff \mathcal{F}$\dnsp,\oss
i.e.\qss a family of\dss Fredholm operators in a fixed\dss Hilbert\dss space,\oss
with\sss this index\sss map and\sss representing\sss
$[\dff \mathbb{K}\dff]\qff -\qff [\dff \mathbb{V}\dff]
\qff \in\qff
K\trf(\trf X\trf)$\nnsp.\oss
We\sss leave\sss this\sss task\sss to\sss the reader.\oss  \eproof

\myuppar{Families parameterized\sss by\sss paracompact\sss spaces.}
When\sss $X$\sss is\dss not\sss compact,\oss the index\dss is\dss still\sss an element\sss of\sss
$K\trf(\trf X\trf)$\nnsp,\oss but\sss it\dss is\dss not\sss reasonable\sss
to define\sss $K\trf(\trf X\trf)$\sss in\sss terms of\dss finitely dimensional\sss
vector\sss bundles on\sss $X$\nnsp.\oss 
Instead,\pss $K\trf(\trf X\trf)$\sss is\dss defined\sss in\sss terms
of\dss homotopy classes\sss to a classifying space,\oss 
with\sss $\mathcal{F}$\sss being a natural\sss choice.\oss
The corresponding\sss definition of\dss index\sss was given\sss by\dss Segal\qss \cite{sfc}\qss
under\sss the assumption\sss that\sss $X$\sss is\dss paracompact.\oss

Sigal's\dss definition applies\sss to families\sss of\dss operators\sss
$D_{\dff x}\dff,\qff x\qff \in\qff X$\sss which\sss he called\qss
\emph{Fredholm\sss operators on}\qss $X$\nnsp.\oss
They are defined as families invertible modulo 
compact\sss families,\oss i.e.\qss such\sss that\sss
there exists a family of\dss operators\sss
$C_{\dff x}\dff \colon\dff 
H_{\dff x}\qff \ttoo\qff H_{\dff x}\dff,\qff 
x\qff \in\qff X$\sss
such\sss that\sss operators\vspace{3pt}
\[
\quad
\id\qff -\pff
D_{\dff x}\dff \circ\dff C_{\dff x}\dff,\off x\qff \in\qff X\dff,
\quad\fff
\mbox{and}\quad\dff
\id\qff -\pff
C_{\dff x}\dff \circ\dff D_{\dff x}\dff,\off x\qff \in\qff X\dff,
\]

\vspace{-12pt}\vspace{3pt}
are compact\sss and,\oss moreover,\oss
define a\qss \emph{compact\dss operator}\qss $\mathbb{H}\qff \ttoo\qff \mathbb{H}$\nnsp.\oss
We will\sss not\sss need\sss the definition of\dss compact\sss operator,\oss
but\sss note\sss that\sss the family of\dss zero operators\dss is\dss a compact\sss operator.\oss
Hence a family of\dss invertible operators\dss is\dss a\dss Fredholm\dss operator\sss in\dss
Segal's\dss sense.\oss
But,\oss in\sss general,\oss such a family\dss is\dss not\sss a\dss Fredholm\dss family\sss in our sense.\oss

Segal\dss considers families\sss of\dss operators\sss
$D_{\dff x}\dff \colon\dff 
H_{\dff x}\qff \ttoo\qff H_{\dff x}\dff,\qff 
x\qff \in\qff X$\sss
as families of\dss chain\sss complexes\sss
$0\qff \ttoo\qff H_{\dff x}\qff \ttoo\qff H_{\dff x}\qff \ttoo\qff 0\dff,\qff x\qff \in\qff X$\dss
(in\sss fact,\oss his\sss theory applies\sss to general\sss chain complexes).\oss
Segal\qss \cite{sfc}\qss
calls\sss a family of\dss operators\sss
$\varphi_{\dff x}\dff \colon\dff
H\qff \ttoo\qff H_{\dff x}\dff,\qff x\qff \in\qff X$\qss
\emph{finite}\oss if\dss locally,\oss 
over open subsets\sss $U\qff \subset\qff X$\sss covering\sss $X$\nnsp,\oss
it\sss can be factored\sss as\sss the composition of\dss the projection\sss
to a finitely dimensional\sss
subbundles of\dss $U\dff \times\dff H\qff \ttoo\qff U$\sss and a continuous map\sss
from\sss this subbundle\sss to\sss $\mathbb{H}$\nnsp.\oss
Suppose\sss that\sss
$D_{\dff x}\dff,\qff x\qff \in\qff X$\sss 
is\dss a\dss Fredholm\dss operator\sss on\sss $X$\nnsp.\oss
Then\sss there exists a family of\trs Fredholm\dss operators\sss
$E_{\dff x}\dff \colon\dff 
H\qff \ttoo\qff H\dff,\qff x\qff \in\qff X$\sss 
continuous as a map\sss 
$\mathbb{E}\dff \colon\dff
X\qff \ttoo\qff\mathcal{F}$\sss 
to $\mathcal{F}$\sss with\sss the norm\sss topology,\oss
and\sss finite\sss families of\dss operators\sss
$\theta_{\dff x}\dff,\qff \phi_{\dff x}\dff \colon\dff
H\qff \ttoo\qff H_{\dff x}\dff,\qff x\qff \in\qff X$\sss
such\sss that\vspace{-1.5pt}
\begin{equation}
\label{segal}
\quad
\begin{tikzcd}[column sep=sboom, row sep=sboom]
0
\arrow[r]
&
H
\arrow[d, "\dis \dff \theta_{\dff x}"]
\arrow[r, "\dis E_{\dff x}"]
&
H
\arrow[d, "\dis \dff \phi_{\dff x}"]
\arrow[r]
&
0
\\
0
\arrow[r]
&
H_{\dff x}
\arrow[r, "\dis D_{\dff x}"]
&
H_{\dff x}
\arrow[r]
&
0
\end{tikzcd}
\end{equation}

\vspace{-12pt}\vspace{4.5pt}
is\dss a chain\sss homotopy equivalence.\oss
The homotopy class of\dss the family\sss
$E_{\dff x}\dff,\qff x\qff \in\qff X$\sss
considered as a map\sss 
$\mathbb{E}\dff \colon\dff X\qff \ttoo\qff \mathcal{F}$\sss
depends only on\sss $D_{\dff x}\dff,\qff x\qff \in\qff X$\nnsp.\oss
Therefore one can define\sss the analytical\dss index of\trs Fredholm\dss operator\sss
$D_{\dff x}\dff,\qff x\qff \in\qff X$\sss on\sss $X$\sss
as\sss the homotopy class of\sss
$\mathbb{E}\dff \colon\dff
X\qff \ttoo\qff \mathcal{F}$\dnsp.\oss

\mypar{Theorem.}{fredholm-agree-paracompact}
\emph{If\qss $X$ is\dss paracompact\sss
and\dss
$D_{\dff x}\dff,\qff 
x\qff \in\qff X$\sss 
is\dss a\dss Fredholm\dss family and\dss a\dss Fredholm\dss operator\sss in\dss
Segal's\dss sense,\oss
then\dss two definitions of\dss the analytical\dss index agree for\sss
$D_{\dff x}\dff,\qff 
x\qff \in\qff X$\nnsp.}

\proof
It\dss is\dss sufficient\sss to prove\sss that,\oss
given a chain\sss homotopy equivalence of\dss the form\qss (\ref{segal}),\oss
the analytical\dss index of\sss $D_{\dff x}\dff,\qff x\qff \in\qff X$\sss
in\sss the sense of\trs Section\qss \ref{analytic-index-fredholm-non-sa}\qss
is\dss equal\sss to\sss that\sss of\sss $E_{\dff x}\dff,\qff x\qff \in\qff X$\nnsp.\oss
Since\qss (\ref{segal})\qss is\dss a chain\sss homotopy equivalence,\oss
the maps\sss $\theta_{\dff x}\dff,\qff \phi_{\dff x}$\sss induce isomorphisms
of\dss the homology groups of\dss the row complexes.\oss
In\sss particular,\pss
$\theta_{\dff x}$\sss induces an\sss isomorphism of\dss kernels\sss
$\kernel\dff E_{\dff x}\qff \ttoo\qff \kernel\dff D_{\dff x}$\sss 
for every\sss $x$\nnsp.\oss
Let\sss us\sss temporarily\sss fix\sss an arbitrary\sss $x\qff \in\qff X$\nnsp.\oss
Since\sss $E_{\dff x}$\sss and\sss $D_{\dff x}$\sss are\dss
Fredholm\dss operators,\oss
there exist\sss $\varepsilon\fff,\qff \varepsilon'\qff >\qff 0$\sss 
such\sss that\sss\vspace{1.5pt}\vspace{1.5pt}
\[
\quad
\kernel\dff E_{\dff x}
\off =\off
\image\dff P_{\qff [\dff 0\fff,\qff \varepsilon\dff]}\trf(\trf \num{E_{\dff x}}\trf)
\quad
\mbox{and}\quad
\kernel\dff D_{\dff x}
\off =\off
\image\dff P_{\qff [\dff 0\fff,\qff \varepsilon'\dff]}\trf(\trf \num{D_{\dff x}}\trf)
\pff.
\]

\vspace{-12pt}\vspace{1.5pt}\vspace{1.5pt}
Moreover,\oss one can assume\sss that\sss 
$(\trf E_{\dff z}\dff,\qff \varepsilon\trf)$\sss
and\sss
$(\trf D_{\dff z}\dff,\qff \varepsilon'\trf)$\sss
are enhanced\dss Fredholm\dss operators for\sss $z\qff \in\qff U_{\dff x}$\sss
for a neighborhood\sss $U_{\dff x}$\sss of\sss $x$\nnsp.\oss

We claim\sss that\sss the parameters\sss $\varepsilon\fff,\qff \varepsilon'\qff >\qff 0$\sss
and\sss the neighborhood\sss $U_{\dff x}$\sss can\sss be chosen\sss in
such a way\sss that\sss the composition\dss 
$P_{\qff [\dff 0\fff,\qff \varepsilon'\dff]}\trf(\trf \num{D_{\dff z}}\trf)
\qff \circ\qff
\theta_{\dff z}$\dss 
induces an\sss isomorphism\vspace{4.5pt}
\begin{equation}
\label{induced-ker}
\quad
\image\dff P_{\qff [\dff 0\fff,\qff \varepsilon\dff]}\trf(\trf \num{E_{\dff z}}\trf)
\off \ttoo\off
\image\dff P_{\qff [\dff 0\fff,\qff \varepsilon'\dff]}\trf(\trf \num{D_{\dff z}}\trf)
\pff
\end{equation}

\vspace{-12pt}\vspace{4.5pt}
for every\sss $z\qff \in\qff U_{\dff x}$.\oss
Let\sss us\sss consider\sss the family of\dss restrictions\vspace{1.5pt}\vspace{1.5pt}
\[
\quad
\theta\fff'_{\fff z}\dff \colon\dff
\image\dff P_{\qff [\dff 0\fff,\qff \varepsilon\dff]}\trf(\trf \num{d_{\dff z}}\trf)
\qff \ttoo\qff
H_{\dff z}
\qff,\quad
z\qff \in\qff U_{\dff x}
\]

\vspace{-12pt}\vspace{1.5pt}\vspace{1.5pt}
induced\sss by\sss $\theta_{\dff z}$.\oss
Since\sss the family\sss
$\image\dff P_{\qff [\dff 0\fff,\qff \varepsilon\dff]}\trf(\trf \num{E_{\dff z}}\trf)\dff,\qff
z\qff \in\qff U_{\fff x}$\sss
is\dss a finitely dimensional\sss vector bundle,\oss
the family 
$\theta\fff'_{\fff z}\dff,\qff z\qff \in\qff U_{\dff x}$\sss
is\dss continuous.\oss

Since\sss $\theta\fff'_{\fff x}$\sss is\dss injective,\oss
after\sss replacing\sss $U_{\dff x}$\sss
by a smaller neighborhood\sss
we can assume\sss that\sss $\theta\fff'_{\fff z}$\sss
is\dss injective for every\sss $z\qff \in\qff U_{\dff x}$\nsp.\oss
Moreover,\oss we can assume\sss that\vspace{1.5pt}\vspace{1.5pt}
\[
\quad
c\qff \norm{v}
\off \leq\off
\norm{\theta\fff'_{\fff z}\trf(\trf v\trf)}
\]

\vspace{-12pt}\vspace{1.5pt}\vspace{1.5pt}
for some\sss $c\qff >\qff 0$\sss
and every\sss $z\qff \in\qff U_{\dff x}$\nsp,\qss
$v\qff \in\qff
\image\dff P_{\qff [\dff 0\fff,\qff \varepsilon\dff]}\trf(\trf \num{d_{\dff z}}\trf)$\nnsp.\oss
Similarly,\oss we can assume\sss that\vspace{4.5pt}
\[
\quad
\norm{\phi_{\fff z}\trf (\qff E_{\dff z}\trf(\trf v\trf)\qff)}
\off \leq\off
c'\qff \norm{E_{\dff z}\trf(\trf v\trf)}
\]

\vspace{-12pt}\vspace{4.5pt}
for some\sss $c'\qff >\qff 0$\sss
and every\sss $z\qff \in\qff U_{\dff x}$\nsp,\qss
$v\qff \in\qff
\image\dff P_{\qff [\dff 0\fff,\qff \varepsilon\dff]}\trf(\trf \num{E_{\dff z}}\trf)$\nnsp.\oss
Then\vspace{4.5pt}
\[
\quad
\norm{D_{\dff z}\trf \left(\qff \theta\fff'_{\fff z}\trf(\trf v\trf)\trf\right)}
\off =\off
\norm{D_{\dff z}\dff \circ\dff \theta\fff'_{\fff z}\trf(\trf v\trf)}
\off =\off
\norm{\phi_{\dff z}\dff \circ\dff E_{\dff z}\trf(\trf v\trf)}
\off \leq\off
c'\qff \norm{E_{\dff z}\trf(\trf v\trf)}
\]

\vspace{-36pt}\vspace{3pt}
\[
\quad
\hspace{17.2em}
\leq\off
c'\qff \varepsilon\qff \norm{v}
\off \leq\off
c'\qff \varepsilon\qff c^{\dff -\dff 1}\qff
\norm{\theta\fff'_{\fff z}\trf(\trf v\trf)}
\pff.
\]

\vspace{-12pt}\vspace{4.5pt}
If\dss 
$c'\qff \varepsilon\qff c^{\dff -\dff 1}\qff \leq\qff \varepsilon'$\nsp,\oss
this implies\sss that\sss the projection\sss
$P_{\qff [\dff 0\fff,\qff \varepsilon'\dff]}\trf(\trf \num{D_{\dff z}}\trf)$\sss
is\dss injective on\sss
the image of\trs $\theta\fff'_{\fff z}$\sss
and\sss hence\sss the maps\sss
$P_{\qff [\dff 0\fff,\qff \varepsilon'\dff]}\trf(\trf \num{D_{\dff z}}\trf)
\qff \circ\fff\qff
\theta\fff'_{\fff z}$\sss
with\sss $z\qff \in\qff U_{\dff x}$\sss
are injective.\oss
Equivalently,\oss the maps\qss (\ref{induced-ker})\qss are injective for\sss
$z\qff \in\qff U_{\dff x}$\nsp.\oss
The sources and\sss the\sss targets of\dss these maps form\sss vector bundles,\oss
and\qss (\ref{induced-ker})\qss is\dss an\sss isomorphism\sss for\sss $z\off =\off x$\sss
and\dss is\dss injective\sss for\sss every\sss $z\qff \in\qff U_{\dff x}$\nsp.\oss
It\sss follows\sss that\qss (\ref{induced-ker})\qss is\dss an\sss isomorphism\sss for\sss
every\sss $z\qff \in\qff U_{\dff x}$\nsp.\oss
This proves our claim.\oss

Similarly,\oss the map\sss
$\phi_{\dff x}$\sss induces an\sss isomorphism of\dss cokernels\sss
$\coker\dff E_{\dff x}\qff \ttoo\qff \coker\dff D_{\dff x}$\nsp,\oss
or,\oss equivalently,\oss an\sss isomorphism of\dss kernels of\dss
the adjoint\sss operators\sss
$\kernel\dff E_{\dff x}^{\fff *}\qff \ttoo\qff \kernel\dff D_{\dff x}^*$\nsp.\oss
As above,\oss there exist\sss
$\delta\fff,\qff \delta'\qff >\qff 0$\sss 
such\sss that\sss\vspace{1.5pt}
\[
\quad
\kernel\dff E_{\dff x}^{\fff *}
\off =\off
\image\dff P_{\qff [\dff 0\fff,\qff \delta\dff]}\trf(\trf \num{E_{\dff x}^{\fff *}}\trf)
\quad
\mbox{and}\quad
\kernel\dff D_{\dff x}^*
\off =\off
\image\dff P_{\qff [\dff 0\fff,\qff \delta'\dff]}\trf(\trf \num{D_{\dff x}^*}\trf)
\pff,
\]

\vspace{-12pt}\vspace{1.5pt}
and,\oss moreover,\pss 
$(\trf E_{\dff x}^{\fff *}\dff,\qff \delta\trf)$\sss
and\sss
$(\trf D_{\dff x}^*\dff,\qff \delta'\trf)$\sss
are enhanced\dss Fredholm\dss operators.\oss
Then\sss
$(\trf E_{\dff z}^{\fff *}\dff,\qff \delta\trf)$\sss
and\sss
$(\trf D_{\dff z}^*\dff,\qff \delta'\trf)$\sss
are enhanced\dss Fredholm\dss operators for\sss $z\qff \in\qff U^{\fff *}_{\dff x}$\sss
for a neighborhood\sss $U^{\fff *}_{\dff x}$\sss of\sss $x$\nnsp.\oss
As above,\oss the parameters\sss 
$\delta\fff,\qff \delta'\qff >\qff 0$\sss
and\sss the neighborhood\sss $U^{\fff *}_{\dff x}$\sss can\sss be chosen\sss in
such a way\sss that\sss for every\sss $z\qff \in\qff U^{\fff *}_{\dff x}$\sss
the composition\sss 
$P_{\qff [\dff 0\fff,\qff \delta'\dff]}\trf(\trf \num{D_{\dff z}^*}\trf)
\qff \circ\qff
\phi_{\dff z}$\sss 
induces an\sss isomorphism\vspace{1.5pt}
\begin{equation}
\label{induced-coker}
\quad
\image\dff P_{\qff [\dff 0\fff,\qff \delta\dff]}\trf(\trf \num{E^{\fff *}_{\dff z}}\trf)
\off \ttoo\off
\image\dff P_{\qff [\dff 0\fff,\qff \delta'\dff]}\trf(\trf \num{D_{\dff z}^*}\trf)
\pff.
\end{equation}

\vspace{-12pt}\vspace{1.5pt}
In\sss the definition of\dss index\sss the same parameter\sss $\varepsilon\qff >\qff 0$\sss
is\dss used\sss both\sss for\sss the operator and\dss its adjoint,\oss
although\sss this\dss is\dss largely a matter of\dss convenience.\oss
By\sss this reason\sss we would\sss like\sss to have\sss
$\varepsilon\off =\off \delta$\sss and\sss
$\varepsilon'\off =\off \delta'$\nnsp.\oss
This can\sss be achieved\sss by choosing\sss first\sss a sufficiently small\sss positive
numbers as\sss $\varepsilon'\off =\off \delta'$\nnsp.\oss
Then a sufficiently small\sss positive number will\sss work as both\sss
$\varepsilon$\sss and\sss $\delta$\nnsp.\oss
We should also replace both\sss $U_{\dff x}$\sss and\sss $U^{\fff *}_{\dff x}$\sss
by\sss their\sss intersection\sss 
$U_{\dff x}\dff \cap\dff U^{\fff *}_{\dff x}$\nsp.\oss

Now we can use\sss $U_{\dff x}$\sss and\sss 
$\varepsilon\off =\off \varepsilon_{\dff x}$\sss
and\sss
$\varepsilon'\off =\off \varepsilon'_{\dff x}$\sss
for\sss $x\qff \in\qff X$\sss
to construct\sss continuous functors\vspace{1.5pt}\vspace{-0.375pt}
\[
\quad
\mathbb{E}_{\qff U,\qff \bm{\varepsilon}}\qff \colon\qff
X_{\dff U}\qff \ttoo\qff \mathcal{S}
\quad
\mbox{and}\quad\dff
\mathbb{D}_{\qff U,\qff \bm{\varepsilon'}}\qff \colon\qff
X_{\dff U}\qff \ttoo\qff \mathcal{S}\dff(\trf \mathbb{H}\trf)
\qff
\]

\vspace{-12pt}\vspace{1.5pt}\vspace{-0.375pt}
and\sss then\sss the corresponding\sss geometric realizations\oss\vspace{1.5pt}
\[
\quad
\num{\mathbb{E}_{\qff U,\qff \bm{\varepsilon}}}\qff \colon\qff
\num{X_{\dff U}}\qff \ttoo\qff \num{\mathcal{S}}
\quad
\mbox{and}\quad\dff
\num{\mathbb{D}_{\qff U,\qff \bm{\varepsilon'}}}\qff \colon\qff
\num{X_{\dff U}}\qff \ttoo\qff \num{\mathcal{S}\dff(\trf \mathbb{H}\trf)}
\pff.
\]

\vspace{-12pt}\vspace{1.5pt}
But\sss the functors and\sss maps in\sss these pairs have different\sss targets.\oss
Let\sss $\mathbb{T}$\sss be\sss the\sss trivial\sss bundle\sss
$\pr\dff \colon\dff X\dff \times\dff H\qff \ttoo\qff H$\nnsp.\oss
By using a\sss trivialization of\dss the bundle\sss $\mathbb{H}$\sss we can\sss
identify\sss it\sss with\sss $\mathbb{T}$\dnsp.\oss
Then\sss we can\sss reinterpret\sss
$\mathbb{E}_{\qff U,\qff \bm{\varepsilon}}$\sss
and\sss
$\mathbb{D}_{\qff U,\qff \bm{\varepsilon'}}$\sss
as continuous functors\vspace{1.5pt}
\[
\quad
\mathbb{E}_{\qff U,\qff \bm{\varepsilon}}\qff \colon\qff
X_{\dff U}\qff \ttoo\qff \mathcal{S}\dff(\trf \mathbb{T}\trf)
\quad
\mbox{and}\quad\dff
\mathbb{D}_{\qff U,\qff \bm{\varepsilon'}}\qff \colon\qff
X_{\dff U}\qff \ttoo\qff \mathcal{S}\dff(\trf \mathbb{T}\trf)
\pff.
\]

\vspace{-12pt}\vspace{1.5pt}
One can\sss hope\sss that\sss the isomorphisms\qss (\ref{induced-ker})\qss
and\qss (\ref{induced-coker})\qss define\sss
isomorphisms of\dss these functors,\oss but\sss 
these isomorphisms are not\sss morphisms of\sss 
$\mathcal{S}\dff(\trf \mathbb{T}\trf)$\nnsp.\oss
The solution of\dss this problem\dss is\dss to\sss lift\sss 
these functors\sss to\sss $S/\fff \mathbb{T}$\nnsp.\oss
If\dss we choose sufficiently small\sss neighborhoods\sss $U_{\dff x}$\nsp,\oss
the bundles\vspace{1.5pt}
\[
\quad
\image\dff P_{\qff [\dff 0\fff,\qff \varepsilon\dff]}\trf(\trf \num{E_{\dff z}}\trf)
\dff,\off
z\qff \in\qff U_{\dff x}
\quad
\mbox{and}\quad\dff
\image\dff P_{\qff [\dff 0\fff,\qff \varepsilon\dff]}\trf(\trf \num{E^{\fff *}_{\dff z}}\trf)
\dff,\off
z\qff \in\qff U_{\dff x}
\]

\vspace{-12pt}\vspace{1.5pt}
will\sss be\sss trivial.\oss
By choosing\sss trivializations of\dss these bundles we can\sss lift\sss
$\mathbb{E}_{\qff U,\qff \bm{\varepsilon}}$\sss to a functor\sss
$\widetilde{\mathbb{E}}_{\qff U,\qff \bm{\varepsilon}}
\dff \colon\dff
X_{\dff U}\qff \ttoo\qff S/\fff \mathbb{T}$\nnsp.\oss
Let\sss use\sss these\sss trivializations and\sss 
the isomorphisms\qss (\ref{induced-ker})\qss and\qss (\ref{induced-coker})\qss
to\sss trivialize\sss the corresponding\sss bundles for\sss $D$\nnsp.\oss
Then\sss the resulting\sss functor\dss
$\widetilde{\mathbb{D}}_{\qff U,\qff \bm{\varepsilon}}
\dff \colon\dff
X_{\dff U}\qff \ttoo\qff S/\fff \mathbb{T}$\sss
will\sss be not\sss only\sss isomorphic,\oss but\sss equal\sss to\dss
$\widetilde{\mathbb{E}}_{\qff U,\qff \bm{\varepsilon}}$\nnsp.\oss
By composing\sss it\sss with\sss
$S/\fff \mathbb{T}\qff \ttoo\qff S$\sss
and\sss then\sss taking\sss the geometric realizations,\oss
we will\sss get\sss equal\sss index\sss maps.\oss  \eproof

\newpage
\mysection{The\qss classical\qss index\qss of\qss
self-adjoint\qss Fredholm\qss families}{classical-analytic-index-sa}

\myuppar{The\dss Atiyah--Singer\dss definition of\dss the analytical\dss index.}
In\sss the papers of\qss Atiyah,\qss Patodi,\oss and\trs Singer\qss 
\cite{as5},\qss \cite{as},\qss \cite{aps}\qss
the analytical\dss index of\dss families of\dss 
self-adjoint\dss Fredholm\dss operators\dss 
is\dss defined only\sss for\sss the 
families of\dss bounded operators 
in a fixed\dss Hilbert\dss space.\oss
The approach of\trs these papers 
naturally extends\sss to operators in\dss Hilbert\dss bundles.\oss
Let\sss $A_{\dff x}\dff \colon\dff H_{\dff x}\qff \ttoo\qff H_{\dff x}$\nsp,\qss
$x\qff \in\qff X$\nnsp,\oss denoted also by\sss $\mathbb{A}$\nnsp.\oss 
be a\dss Fredholm\dss family of\dss self-adjoint\sss operators in\sss the fibers of\dss
a\sss Hilbert\dss bundle\sss $H_{\dff x}\dff,\qff x\qff \in\qff X$\nnsp,\oss
denoted also by\sss $\mathbb{H}$\nnsp.\oss
Following\qss \cite{as},\qss \cite{aps},\oss let\sss us consider\sss
the family\sss $\mathbb{B}$\sss of\trs Fredholm\dss operators 
defined\sss by\vspace{1.5pt}
\begin{equation}
\label{as-loop}
\quad
B_{\dff x\fff,\qff t}
\off =\off
\id_{\trf H}\qff \cos t
\pff +\pff
i\trf A_{\dff x}\qff \sin t
\qff,
\quad\off
x\qff \in\qff X\dff,\off\qff
t\qff \in\qff [\trf 0\fff,\qff \pi\trf]
\pff,
\end{equation}

\vspace{-12pt}\vspace{1.5pt}
where\sss
$\id_{\trf H}$\sss is\dss the identity operator\sss
$H\qff \ttoo\qff H$\sss and\dss
$i\off =\off \sqrt{\dff -\qff 1}$\nnsp.\oss
By\trs Theorem\qss \ref{fr-to-fr}\qss below,\oss
this\dss is\dss a\dss Fredholm\dss family.\oss
By a\qss \emph{homotopy}\qss of\dss such a\sss family,\oss
or,\oss more generally,\oss a family\sss parameterized\sss by\sss
$X\dff \times\dff [\trf 0\fff,\qff \pi\trf]$\nnsp,\oss
we will\sss understand a homotopy fixed on\sss
$X\dff \times\dff 0$\sss and\sss $X\dff \times\dff \pi$\nnsp.\oss

By\sss the\qss \emph{classical\sss analytical\dss index}\pss of\dss $\mathbb{A}$\sss
we will\sss understand\sss the analytical\dss index 
of\dss $\mathbb{B}$\sss either\sss in\sss the sense of\trs 
Atiyah\dss and\dss Singer\qss \cite{as4},\oss or\sss in\sss the sense of\trs
Section\qss \ref{analytic-index-fredholm-non-sa}.\oss
The\sss latter\dss is\dss the homotopy class of\dss the index\sss map\dss 
$X\dff \times\dff [\trf 0\fff,\qff \pi\trf]\qff \ttoo\qff \num{\mathcal{S}}$\sss
with respect\sss to homotopies fixed on\sss
$X\dff \times\dff 0$\sss and\sss $X\dff \times\dff \pi$\nnsp.\oss
Classically\sss the operators\sss $A_{\dff x}$\sss
are bounded,\pss $H_{\dff x}\off =\off H$\sss for every $x$\nnsp,\oss and\sss
the\sss target\sss of\dss index\sss maps\dss is\dss $\mathcal{F}$\nnsp.\oss
Then\sss one can\sss take\sss the map\sss
$(\trf x\fff,\qff t\trf)
\off \longmapsto\off
B_{\dff x\fff,\qff t}$\sss
as\sss the index\sss map.

\mypar{Lemma.}{as-is-fredholm-lemma}
\emph{Let\sss $(\trf A\fff,\qff \varepsilon\trf)$\sss be an enhanced\sss self-adjoint\dss
Fredholm\dss operator\sss in\sss $H$\nnsp,\oss
and\sss let\dss}\vspace{3pt}
\[
\quad
B_{\dff t}
\off =\off\dff
\id_{\trf H}\qff \cos t
\pff +\pff
i\trf A\qff \sin t
\pff,
\]

\vspace{-12pt}\vspace{3pt}
\emph{where\dss $t\qff \in\qff [\trf 0\fff,\qff \pi\trf]$\nnsp.\oss
Let\qss $\delta\qff >\qff 0$\nnsp.\oss
If\qss either\dss $\sin t\off =\off 0$ and\qss $\delta\qff <\qff 1$\nnsp,\oss
or\qss $\delta^{\dff 2}\qff <\qff \cos^{\dff 2} t$\nnsp,\oss
then}\vspace{3pt}
\[
\quad
\image\dff P_{\qff [\dff 0\fff,\qff \delta\dff]}\qff
\left(\qff \num{B_{\dff t}}\qff\right)
\off =\off\dff
0
\pff.
\]

\vspace{-12pt}\vspace{3pt}
\emph{If\qss $\sin t\off \neq\off 0$\sss and\qss
$\delta^{\dff 2}
\off =\off
\cos^{\dff 2} t\qff +\qff \varepsilon^{\dff 2}\qff \sin^{\dff 2} t$\nnsp,\oss
then\qss
$\delta\qff \not\in\qff \sigma\trf(\trf \num{B_{\dff t}}\trf)$\sss and}\vspace{3pt}
\[
\quad
\image\dff P_{\qff [\dff 0\fff,\qff \delta\dff]}\qff
\left(\qff \num{B_{\dff t}}\qff\right)
\off =\off\dff
\image\dff P_{\qff [\dff 0\fff,\qff \varepsilon\dff]}\qff
\left(\qff \num{A}\qff\right)
\pff.
\]

\vspace{-12pt}\vspace{3pt}
\proof
The case of\sss $\sin t\off =\off 0$\sss is\dss trivial.\oss
Suppose\sss that\sss $\sin t\off \neq\qff 0$\nnsp.\oss
Since\sss $A^*\off =\off A$\nnsp,\oss\vspace{3pt}
\[
\quad
\num{B_{\dff t}}^{\dff 2}
\off =\off\dff
B_{\dff t}\qff B^*_{\dff t}
\off =\off\dff
B^*_{\dff t}\qff B_{\dff t}
\off =\off\dff
\id_{\trf H}\qff \cos^{\dff 2} t
\pff +\pff
A^{\dff 2}\qff \sin^{\dff 2} t
\pff.
\]

\vspace{-12pt}\vspace{3pt}
If\qss
$\lambda^2\qff \in\qff \sigma\trf(\trf B_{\dff t}\qff B^*_{\dff t}\trf)$\nnsp,\oss
then\sss the operator\vspace{3pt}
\[
\quad
\id_{\trf H}\qff\lambda^2\qff -\qff B_{\dff t}\qff B^*_{\dff t}
\off\dff =\off\qff
\id_{\trf H}\qff (\trf \lambda^2\qff -\qff \cos^{\dff 2} t\trf)
\pff -\pff
A^{\dff 2}\qff \sin^{\dff 2} t
\pff
\]

\vspace{-12pt}\vspace{3pt}
is\dss not\sss invertible.\oss 
Since\dss $\sin t\off \neq\off 0$\nnsp,\oss
this implies\sss that\vspace{1.5pt}\vspace{0.6pt}
\[
\quad
\frac{\lambda^2\qff -\qff \cos^{\dff 2} t}{\sin^{\dff 2} t}
\off \in\off
\sigma\trf(\trf A^2\trf)
\pff.
\]

\vspace{-12pt}\vspace{1.5pt}\vspace{0.6pt}
Since\sss $A^2\off =\off A^{\fff *}\fff A$\sss is\dss a positive operator,\oss
this implies\sss that\sss
$\lambda^2\pff \geq\pff \cos^{\dff 2} t$\nnsp.\oss
Therefore\sss\vspace{3pt}
\[
\quad
\image\dff P_{\qff [\dff 0\fff,\qff \delta\dff]}\qff
\left(\qff \num{B_{\dff t}}\qff\right)
\off =\off\dff
\image\dff P_{\qff [\dff 0\fff,\qff \delta^{\dff 2}\dff]}\qff
\left(\qff \num{B_{\dff t}}^{\dff 2}\qff\right)
\off =\off\dff
0
\pff
\]

\vspace{-12pt}\vspace{3pt}
if\qss $\delta^{\dff 2}\qff <\qff \cos^{\dff 2} t$\nnsp.\qff\oss
If\qss 
$\lambda^2
\pff \leq\pff 
\cos^{\dff 2} t\qff +\qff \varepsilon^{\dff 2}\qff \sin^{\dff 2} t$\nnsp,\oss
then\oss\vspace{1.5pt}\vspace{0.6pt}
\[
\quad
\frac{\lambda^2\qff -\qff \cos^{\dff 2} t}{\sin^{\dff 2} t}
\off \leq\off
\varepsilon^{\dff 2}
\qff.
\]

\vspace{-12pt}\vspace{1.5pt}\vspace{0.6pt}
Since\sss $(\trf A\fff,\qff \varepsilon\trf)$\sss 
is\dss an enhanced\sss self-adjoint\sss operator,\oss
in\sss this case\vspace{1.5pt}\vspace{0.6pt}
\[
\quad
\frac{\lambda^2\qff -\qff \cos^{\dff 2} t}{\sin^{\dff 2} t}
\off =\off
\mu^{\dff 2}
\dff
\]

\vspace{-12pt}\vspace{1.5pt}\vspace{0.6pt}
for some eigenvalue\sss 
$\mu\qff <\qff \varepsilon$\sss 
of\sss $A$\sss of\dss finite multiplicity.\oss 
It\sss follows\sss that\sss
$\delta\qff \not\in\qff \sigma\trf(\trf \num{B_{\dff t}}\trf)$\sss and\vspace{3pt}
\[
\quad
\image\dff P_{\qff [\dff 0\fff,\qff \delta\dff]}\qff
\left(\qff \num{B_{\dff t}}\qff\right)
\off =\off\dff
\image\dff P_{\qff [\dff 0\fff,\qff \varepsilon\dff]}\qff
\left(\qff \num{A}\qff\right)
\pff
\]

\vspace{-12pt}\vspace{3pt}
if\qss
$\delta^{\dff 2}
\off =\off
\cos^{\dff 2} t\qff +\qff \varepsilon^{\dff 2}\qff \sin^{\dff 2} t$\nnsp.\oss
This completes\sss the proof.\oss  \eproof

\mypar{Theorem.}{fr-to-fr}
\emph{If\qss $\mathbb{A}$ is\dss a self-adjoint\trs Fredholm\dss family,\oss
then\qss $\mathbb{B}$ is\dss a\trs Fredholm\dss family.\oss}

\proof
Let\sss $x\qff \in\qff X$\nnsp.\oss
Then\sss there exists\sss $\alpha\off =\off \alpha_{\dff x}\qff >\qff 0$\sss such\sss that\sss 
$(\trf A_{\dff x}\fff,\qff \alpha\trf)$\sss 
is\dss an enhanced\sss self-adjoint\sss operator.\oss
Moreover,\oss there exists\sss 
$\beta\off =\off \beta_{\dff x}\qff \in\qff (\trf 0\fff,\qff \alpha\trf)$\sss
such\sss that\sss
$(\trf A_{\dff x}\fff,\qff \varepsilon\trf)$\sss 
is\dss an enhanced\sss self-adjoint\sss operator\sss for every\sss
$\varepsilon\qff \in\qff [\trf \beta\fff,\qff \alpha\trf]$\nnsp.\oss
Since\sss $\mathbb{A}$\sss is\dss a self-adjoint\trs Fredholm\dss family,\oss
there exists a neighborhood\sss $U\off =\off U_{\dff x}$\sss of\sss $x$\sss such\sss that\sss
$(\trf A_{\dff z}\fff,\qff \varepsilon\trf)$\sss 
is\dss an enhanced\sss self-adjoint\sss operator\sss for every\sss
$\varepsilon\qff \in\qff [\trf \beta\fff,\qff \alpha\trf]$\sss
and\sss $z\qff \in\qff U$\nnsp.\oss
Since\sss $\beta\qff <\qff \alpha$\nnsp,\oss
there exists\sss $\delta\off =\off \delta_{\dff x}\qff >\qff 0$\sss such\sss that\vspace{3pt}
\[
\quad
\cos^{\dff 2} t\qff +\qff \beta^{\dff 2}\qff \sin^{\dff 2} t
\off <\off
\delta^{\dff 2}
\off <\off
\cos^{\dff 2} t\qff +\qff \alpha^{\dff 2}\qff \sin^{\dff 2} t
\]

\vspace{-12pt}\vspace{3pt}
for all\sss $t$\sss sufficiently close\sss to\sss $\pi/2$\nnsp,\oss
say,\oss for\sss $t\qff \in\qff (\trf u\fff,\qff \pi\qff -\qff u\trf)$\sss
for some\sss $u\off =\off u_{\dff x}\qff \in\qff (\trf 0\fff,\qff \pi/2\trf)$\nnsp.\oss
For every such\sss $t$\sss the number\sss $\delta^{\dff 2}$\sss
can\sss be written\sss in\sss the form\vspace{3pt}
\[
\quad
\delta^{\dff 2}
\off =\off
\cos^{\dff 2} t\qff +\qff \varepsilon\trf(\dff t\trf)^{\dff 2}\qff \sin^{\dff 2} t
\]

\vspace{-12pt}\vspace{3pt}
for some\sss
$\varepsilon\trf(\dff t\trf)\qff \in\qff (\trf \beta\fff,\qff \alpha\trf)$\nnsp.\oss
Lemma\qss \ref{as-is-fredholm-lemma}\qss implies\sss that\sss
$\delta\qff \not\in\qff \sigma\trf(\trf \num{B_{\dff z\fff,\qff t}}\trf)$\sss and\vspace{3pt}
\[
\quad
\image\dff P_{\qff [\dff 0\fff,\qff \delta\dff]}\qff
\left(\qff \num{B_{\dff z\fff,\qff t}}\qff\right)
\off =\off\dff
\image\dff P_{\qff [\dff 0\fff,\qff \varepsilon\trf(\dff t\trf)\dff]}\qff
\left(\qff \num{A_{\dff z}}\qff\right)
\pff
\]

\vspace{-12pt}\vspace{3pt}
for every\sss $z\qff \in\qff U$\sss and\sss
$t\qff \in\qff (\trf u\fff,\qff \pi\qff -\qff u\trf)$\nnsp.\oss
Also,\oss since\sss
$(\trf A_{\dff z}\fff,\qff \varepsilon\trf)$\sss 
is\dss an enhanced\sss self-adjoint\sss operator\sss for every\sss
$\varepsilon\qff \in\qff [\trf \beta\fff,\qff \alpha\trf]$\nnsp,\oss
the spectrum\sss $\sigma\trf(\trf A_{\dff z}\trf)$\sss
is\dss disjoint\sss from\sss $[\trf \beta\fff,\qff \alpha\trf]$\sss
and\sss hence\vspace{3pt}
\[
\quad
\image\dff P_{\qff [\dff 0\fff,\qff \varepsilon\trf(\dff t\trf)\dff]}\qff
\left(\qff \num{A_{\dff z}}\qff\right)
\off =\off\dff
\image\dff P_{\qff [\dff 0\fff,\qff \alpha\dff]}\qff
\left(\qff \num{A_{\dff z}}\qff\right)
\pff.
\]

\vspace{-12pt}\vspace{3pt}
The\sss last\sss two equalities\sss imply\sss that\vspace{3pt}
\[
\quad
\image\dff P_{\qff [\dff 0\fff,\qff \delta\dff]}\qff
\left(\qff \num{B_{\dff z\fff,\qff t}}\qff\right)
\off =\off\dff
\image\dff P_{\qff [\dff 0\fff,\qff \alpha\dff]}\qff
\left(\qff \num{A_{\dff z}}\qff\right)
\pff
\]

\vspace{-12pt}\vspace{3pt}
for every\sss $z\qff \in\qff U$\sss and\sss
$t\qff \in\qff (\trf u\fff,\qff \pi\qff -\qff u\trf)$\nnsp.\oss
In\sss turn,\oss this implies\sss that\sss the pair\sss
$(\trf B_{\dff z\fff,\qff t}\dff,\qff \delta\trf)$\sss
is\dss an enhanced operator\sss for every\sss
for every\sss $z\qff \in\qff U$\sss and\dss
$t\qff \in\qff (\trf u\fff,\qff \pi\qff -\qff u\trf)$\nnsp,\oss
and\dss hence\sss the family\sss\vspace{3pt}
\[
\quad
B_{\dff z\fff,\qff t}\dff,\pff 
z\qff \in\qff U\dff,\off
t\qff \in\qff (\trf u\fff,\qff \pi\qff -\qff u\trf)
\]

\vspace{-12pt}\vspace{3pt}
is\dss a\dss Fredholm\dss family.\oss
Now,\oss let\sss us choose some\sss 
$w\off =\off w_{\dff x}\qff \in\qff (\trf u\fff,\qff \pi/2\trf)$\nnsp.\oss
Lemma\qss \ref{as-is-fredholm-lemma}\qss implies\sss that\sss 
$(\trf B_{\dff z\fff,\qff t}\dff,\qff \cos w\trf)$\sss
is\dss an enhanced operator\dss if\dss $z\qff \in\qff U$\sss
and\sss  $\cos^{\dff 2} w\pff <\pff \cos^{\dff 2} t$\nnsp,\oss
i.e.\qss if\qss
$t\qff <\qff w$\sss or\sss $t\qff >\qff \pi\qff -\qff w$\nnsp.\oss
It\sss follows\sss that\sss the family\sss
$B_{\dff z\fff,\qff t}\dff,\pff 
z\qff \in\qff U\dff,\off
t\qff \in\qff [\trf 0\fff,\qff w\trf)\dff \cup\dff (\trf \pi\qff -\qff w\fff,\qff \pi\trf]$\sss
is\dss a\dss Fredholm\dss family.\oss
It\sss follows\sss that\sss
$B_{\dff z\fff,\qff t}\dff,\pff 
z\qff \in\qff U\dff,\off
t\qff \in\qff [\trf 0\fff,\qff \pi\trf]$\sss
is\dss a\dss Fredholm\dss family.\oss
Since\sss $x\qff \in\qff U$\sss and\sss $x\qff \in\qff X$\sss is\dss arbitrary,\oss  
we see\sss that\sss 
$B_{\dff x\fff,\qff t}\dff,\pff x\qff \in\qff X$\sss 
is\dss a\dss Fredholm\dss family.\oss  \eproof\vspace{0.5pt}

\mypar{Theorem.}{s-loops-shat}
\emph{There\dss is\dss a canonical\dss homotopy\sss equivalence
$\num{\hat{\mathcal{S}}}\qff \ttoo\qff \Omega\trf \num{\mathcal{S}}$\nnsp,\oss where\sss
$\Omega\trf \num{\mathcal{S}}$\sss is\dss 
the\sss loop space of\sss $\num{\mathcal{S}}$\nnsp.\oss}\vspace{0.5pt}

\proof
This\sss immediately\sss follows from\qss \cite{i2},\oss Theorems\qss 9.8\qss and\qss 15.4.\oss 
This also\sss follows\sss from\trs Theorem\qss 15.6\qss there,\oss
and\dss is\dss implicit\sss in\sss its\sss proof.\oss  \eproof\vspace{0.5pt}

\myuppar{Comparing\sss two definitions of\dss the analytical\dss index.}
Since\sss
$B_{\dff x\fff,\qff 0}\off =\off \id_{\trf H}$\sss
and\sss
$B_{\dff x\fff,\qff \pi}\off =\off -\qff \id_{\trf H}$\sss
for every\sss $x\qff \in\qff X$\nnsp,\oss
one can consider\sss $\mathbb{B}$\sss as a family\sss
parameterized\sss by\sss the space\sss $\Sigma\dff X$\sss obtained\sss
by collapsing in\sss $X\dff \times\dff [\trf 0\fff,\qff \pi\trf]$\sss
the subspace\sss $X\dff \times\dff 0$\sss into a point,\oss
and\sss the subspace\sss
$X\dff \times\dff \pi$\sss into another\sss point.\oss
This\sss turns\sss the homotopies\sss fixed on\sss
$X\dff \times\dff 0$\sss and\sss $X\dff \times\dff \pi$\sss
into\sss the homotopies of\dss maps with\sss the domain\sss $\Sigma\dff X$\nnsp.\oss
In view of\dss this\sss the analytical\dss index of\dss the family\sss $\mathbb{B}$\sss
can\sss be defined as\sss the homotopy class of\dss an\sss index\sss map\sss
$\Sigma\dff X\qff \ttoo\qff \num{\mathcal{S}}$\nnsp.\oss
By\trs Theorem\qss \ref{fredholm-agree}\qss this definition\dss is\dss
equivalent\sss to\sss the classical\sss one.\oss\vspace{0.5pt} 

So,\oss the classical\sss analytical\dss index of\dss the family $\mathbb{A}$
is\dss the homotopy class of\dss the index map\sss
$\Sigma\dff X\qff \ttoo\qff \num{\mathcal{S}}$\sss
of\dss the family\sss $\mathbb{B}$\nnsp.\oss
The analytical\dss index of $\mathbb{A}$ in\sss the sense of\trs
Section\trs \ref{analytic-index-section}\trs is\dss the homotopy class
of\dss an\sss index\sss map
$X\dff \ttoo\trf \num{\hat{\mathcal{S}}}$\nnsp.\oss
In\sss order\sss to compare\sss these definitions,\pss
recall\sss that\sss there\dss is\dss a canonical\sss one-to-one
correspondence between\sss maps
$\Sigma\dff X\qff \ttoo\qff \num{\mathcal{S}}$\sss
and\sss maps
$X\qff \ttoo\qff \Omega\trf \num{\mathcal{S}}$\nnsp.\oss
The same\dss is\dss true for\sss homotopies,\oss
and we can consider\sss the classical\sss analytical\dss index of $\mathbb{A}$
as\sss the homotopy class of\dss a map\sss
$X\qff \ttoo\qff \Omega\trf \num{\mathcal{S}}$\sss
deserving\sss to be called\dss the\qss \emph{classical\dss index\dss map}.\oss
Theorem\qss \ref{s-loops-shat}\qss
allows us\sss to interpret\sss the analytical\dss index\sss in\sss the sense of\trs 
Section\qss \ref{analytic-index-section}\qss as\sss the homotopy class of\dss
a map\sss $X\qff \ttoo\qff \Omega\trf \num{\mathcal{S}}$\nnsp,\oss
which also may\sss be called an\qss \emph{index\dss map}.\oss

\myuppar{Finite-polarized operators.}
We will\sss call\sss a self-adjoint\sss operator\sss
$A\dff \colon\dff K\qff \ttoo\qff K$\sss in a\sss Hilbert\sss space\sss $K$\dss
\emph{finite-polarized}\oss if\dss
$\norm{A}\off =\off 1$\nnsp,\oss
the essential\sss spectrum of\sss $A$\sss is\dss consists of\dss two points\sss
$-\qff 1\fff,\qff 1$\nnsp,\oss
and\dss the spectral\sss projection\sss
$P_{\trf (\dff -\qff 1\fff,\qff 1\trf)}\trf(\trf A\trf)$\sss
is\dss an operator of\dss finite rank.\oss
If\dss we omit\sss the\sss last\sss property,\oss 
we will\sss get\sss exactly\sss the operators from\sss the space\sss $\hat{F}_{\dff *}$
from\qss \cite{as},\oss Section\qss 2.\oss

\myuppar{Finite-polarized\sss replacements.}
We would\dss like\sss to replace\sss the family\sss $\mathbb{A}$\sss
by a family\sss $\mathbb{A}'$\sss of\dss finite-polarized operators\sss
$A'_{\dff x}\dff \colon\dff H_{\dff x}\qff \ttoo\qff H_{\dff x}$\sss
without\sss changing neither\sss the classical\sss analytical\dss index,\oss
nor\sss the analytical\dss index\sss in\sss the sense 
of\trs Section\qss \ref{analytic-index-section}.\oss
This can\sss be done by\sss a spectral\sss deformation\sss similar\sss to
spectral\sss deformations used\sss in\qss \cite{as}.\oss 

We will\sss assume\sss that\sss $X$\sss is\dss paracompact.\oss
Let\sss us choose an atlas for\sss $\mathbb{A}$\nnsp,\oss
i.e.\qss a family of\dss pairs\sss $(\trf U\fff,\qff \varepsilon_{\trf U}\dff)$\sss
adapted\sss to\sss $\mathbb{A}$\sss such\sss that\sss the sets\sss $U$\sss
form an open covering of\sss $X$\nnsp.\oss
We can assume\sss that\sss
$\varepsilon_{\trf U}\qff <\qff 1$\sss for every\sss $U$\nnsp.\oss
Since\sss $X$\sss is\dss paracompact,\oss we can further assume\sss that\sss
this covering\dss is\dss locally\sss finite and\sss
there exists a partition of\dss unity\sss subordinated\sss to\sss this covering.\oss
Let\sss 
$r_{\trf U}\dff \colon\dff X\qff \ttoo\qff \rrr_{\qff \geq\dff 0}$\sss 
be\sss the function from\dss this partition of\dss unity corresponding\sss to\sss $U$\nnsp.\oss
Let\vspace{1.5pt}
\[
\quad
r\trf(\trf x\trf)
\off =\off
\sum\nolimits_{\qff U}\qff r_{\trf U}\dff(\trf x\trf)\qff \varepsilon_{\trf U}
\pff.
\]

\vspace{-12pt}\vspace{1.5pt}
Then\sss $r\trf(\trf x\trf)\qff \in\qff (\trf 0\fff,\qff 1\trf)$\sss and\sss
$r\trf(\trf x\trf)\qff \leq\qff \max_{\trf U\qff \ni\qff x}\qff \varepsilon_{\trf U}$\nsp,\oss
where\sss the maximum\dss is\dss taken over\sss $U$\sss 
such\sss that\sss $x\qff \in\qff U$\dnsp,\oss
for every\sss $x\qff \in\qff X$\nnsp.\oss
It\sss follows\sss that\sss for every\sss $x\qff \in\qff X$\sss
the essential\sss spectrum of\sss $A_{\dff x}$\sss
is\dss disjoint\sss from\sss
$(\trf -\qff r\trf(\trf x\trf)\fff,\qff r\trf(\trf x\trf)\qff)$\nnsp.\oss
For\sss $r\qff >\qff 0$\sss let\sss 
$\chi_{\dff r}\dff \colon\dff \rrr\qff \ttoo\qff \rrr$\sss
be\sss the function\sss defined\sss by\vspace{3pt}\vspace{-0.125pt}
\[
\quad
\chi_{\dff r}\dff(\trf u\trf)
\off =\off
u
\quad
\mbox{for}\quad
0\qff \leq\qff u\qff \leq\qff r/2\qff,
\]

\vspace{-36pt}
\[
\quad
\chi_{\dff r}\dff(\trf u\trf)
\off =\off
(\trf 2/r\qff -\qff 1\trf)\trf u\qff +\qff r\qff -\qff 1
\quad
\mbox{for}\quad
r/2\qff \leq\qff u\qff \leq\qff r\qff,
\]

\vspace{-36pt}
\[
\quad
\chi_{\dff r}\dff(\trf u\trf)
\off =\off
1
\quad
\mbox{for}\quad
r\qff \leq\qff u\qff,
\]

\vspace{-36pt}
\[
\quad
\mbox{and}\quad
\chi_{\dff r}\dff(\trf u\trf)
\off =\off
-\pff \chi_{\dff r}\dff(\trf -\qff u\trf)
\quad
\mbox{for\sss every}\quad
u\qff.
\]

\vspace{-12pt}\vspace{3pt}\vspace{-0.125pt}
In\sss particular,\pss $\chi_{\trf _1}$\sss is\dss the identity\sss map\sss 
$\rrr\qff \ttoo\qff \rrr$\nnsp.\oss
For\sss $x\qff \in\qff X$\sss let\vspace{1.5pt}
\[
\quad
A'_{\dff x}
\off =\off
\chi_{\dff r\trf(\trf x\trf)}\trf\bigl(\trf A_{\dff x}\trf\bigr)
\qff \colon\qff 
H_{\dff x}\qff \ttoo\qff H_{\dff x}
\pff
\]

\vspace{-12pt}\vspace{1.5pt}
Clearly,\oss operators\sss $A'_{\dff x}$\sss are finite-polarized and\sss
$A'_{\dff x}\dff,\pff x\qff \in\qff X$\sss is\dss a\dss Fredholm\dss family.\oss
This\dss is\dss our\qss \emph{finite-polarized\sss replacement}\dss $\mathbb{A}'$\nnsp.\oss
Clearly,\oss if\sss
$0\qff <\qff \varepsilon\qff <\qff r\trf(\trf x\trf)/2$\nnsp,\oss
then\vspace{1.5pt}
\[
\quad
\image\dff P_{\qff [\dff \varepsilon\fff,\qff \infty\dff)}\qff
\left(\trf A_{\dff x}\trf\right)
\off =\off\dss
\image\dff P_{\qff [\dff \varepsilon\fff,\qff \infty\dff)}\qff
\left(\trf A'_{\dff x}\trf\right)
\pff.
\]

\vspace{-12pt}\vspace{1.5pt}
It\sss follows\sss that\dss if\dss $\mathbb{A}$\sss is\dss strictly\dss Fredholm,\oss
then\sss $\mathbb{A}'$\sss is\dss also strictly\dss Fredholm.\oss

\mypar{Theorem.}{index-replacement}
\emph{For both definitions of\dss the analytical\dss index,\oss
the analytical\dss index of\dss $\mathbb{A}$\sss is\dss equal\sss to\sss
the analytical\dss index\sss of\qss its\trs finite-polarized\sss replacements\sss
$\mathbb{A}'$\dnsp.\oss}

\proof
The open covering\sss
$U_{\dff a}\dff,\pff a\qff \in\qff \Sigma$\sss of\dss $X$\sss
and\sss the family of\dss positive numbers\sss
$\varepsilon_{\dff a}\dff,\pff a\qff \in\qff \Sigma$\sss
used\sss to construct\sss the index functor\vspace{0.75pt}
\[
\quad
\mathbb{A}_{\qff U,\qff \bm{\varepsilon}}\dff \colon\dff
X_{\dff U}\qff \ttoo\qff \hat{\mathcal{S}}\dff(\trf \mathbb{H}\trf)
\pff
\]

\vspace{-12pt}\vspace{0.75pt}
don't\sss need\sss to be equal\sss to\sss the covering\sss by\sss the sets\sss $U$\sss
and\sss the numbers\sss $\varepsilon_{\trf U}$\sss used\sss to construct\sss
the finite-polarized\sss replacement\sss $\mathbb{A}'$\nnsp.\oss
By\sss choosing\sss sufficiently small\sss open sets\sss $U_{\dff a}$\sss
and\sss numbers\sss $\varepsilon_{\dff a}$\sss we can ensure\sss that\sss
$\varepsilon_{\dff a}\qff <\qff r\trf(\trf x\trf)/\fff 2$\sss
for every\sss $x\qff \in\qff U_{\dff a}$\nsp.\oss
Then\sss the same covering\sss and\sss the same numbers\sss will\sss
work\sss for\sss both\sss $\mathbb{A}$\sss and\sss $\mathbb{A}'$\nnsp.\oss
Moreover,\oss the index\sss functors will\dss be equal,\vspace{1.5pt}
\[
\quad
\mathbb{A}_{\qff U,\qff \bm{\varepsilon}}
\off =\off
\mathbb{A}'_{\qff U,\qff \bm{\varepsilon}}
\dff \colon\dff
X_{\dff U}\qff \ttoo\qff \hat{\mathcal{S}}\dff(\trf \mathbb{H}\trf)
\pff.
\]

\vspace{-12pt}\vspace{1.5pt}
Therefore\sss the index\sss maps\sss for\sss $\mathbb{A}$\sss
and\sss $\mathbb{A}'$\sss can\sss be assumed\sss to be equal.\oss
It\sss follows\sss that\sss the analytical\dss index of\dss $\mathbb{A}$\sss
in\sss the sense of\trs Section\qss \ref{analytic-index-section}\qss
is\dss equal\sss to\sss that\sss of\dss $\mathbb{A}'$\nnsp.\oss
This proves\sss the\sss theorem for\dss the analytical\dss index\sss
in\sss the sense of\trs Section\qss \ref{analytic-index-section}.\oss

Let\sss us consider\sss now\sss the classical\sss analytical\dss index.\oss
We may\sss consider\sss
the classical\sss analytical\dss index of\sss $\mathbb{A}$\sss
as\sss the homotopy class of\dss an\sss index\sss map\sss
$\Sigma\dff X\qff \ttoo\qff \num{\mathcal{S}}$\nnsp,\oss
or,\oss equivalently,\sss of\dss the index\sss map\sss
$X\dff \times\dff [\trf 0\fff,\qff \pi\trf]\qff \ttoo\qff \num{\mathcal{S}}$\nnsp,\oss
for\dss the family\sss $\mathbb{B}$\nnsp.\oss
Similarly,\oss the classical\sss analytical\dss index of\sss
$\mathbb{A}'$\sss
is\dss the homotopy class of\dss an\sss 
index map of\dss the family\sss\vspace{1.5pt}\vspace{1pt}
\[
\quad
B'_{\dff x\fff,\qff t}
\off =\off
\id_{\trf H}\qff \cos t
\pff +\pff
i\trf A'_{\dff x}\qff \sin t
\qff,
\quad\off
x\qff \in\qff X\dff,\off\qff
t\qff \in\qff [\trf 0\fff,\qff \pi\trf]
\pff,
\]

\vspace{-12pt}\vspace{1.5pt}\vspace{1pt}
which we will\sss denote also by\sss $\mathbb{B}'$\nnsp.\oss
The proof\dss of\trs Theorem\qss \ref{fr-to-fr}\qss includes a construction
of\dss open subsets\sss $U_{\dff x}\qff \subset\qff X$\sss and\sss numbers\dss
$\delta_{\dff x}\dff,\off u_{\dff x}\dff,\off w_{\dff x}$\dss
such\sss that\sss the open subsets\vspace{1.5pt}\vspace{1pt}
\[
\quad
U_{\dff x}
\dff \times\dff
(\trf u_{\dff x}\fff,\qff \pi\qff -\qff u_{\dff x}\trf)
\quad
\mbox{and}\quad
U_{\dff x}
\dff \times\trf
\bigl(\qff
[\trf 0\fff,\qff w_{\dff x}\trf)\dff \cup\dff (\trf \pi\qff -\qff w_{\dff x}\fff,\qff \pi\trf]
\qff\bigr)
\]

\vspace{-12pt}\vspace{1.5pt}\vspace{1pt}
of\dss $X\dff \times\dff [\trf 0\fff,\qff \pi\trf]$\dss
with\sss parameters\sss $\delta_{\dff x}$\sss and\sss $\cos w_{\dff x}$\sss
respectively\sss can\sss be used\sss for constructing\sss an\sss index\sss map\sss
for\sss $\mathbb{B}$\nnsp.\oss
The starting\sss point\sss of\dss this construction\dss is\dss a choice of\dss
parameters\sss $\alpha_{\dff x}\qff >\qff 0$\nnsp.\oss
One can always choose\sss parameters\sss $\alpha_{\dff x}\qff <\qff r\trf(\trf x\trf)/\fff 2$\nnsp.\oss
Then\trs Lemma\qss \ref{as-is-fredholm-lemma}\qss implies\sss that\sss
the same open subsets $U_{\dff x}$\sss and\sss numbers\dss
$\delta_{\dff x}\dff,\off u_{\dff x}\dff,\off w_{\dff x}$\dss
work also for\sss $\mathbb{B}'$\nnsp,\oss
and\dss the resulting\sss index\sss maps\sss for\sss 
$\mathbb{B}$\sss and\sss $\mathbb{B}'$\sss
are equal.\oss
This proves\sss the\sss theorem for\dss 
the classical\sss analytical\dss index.\oss  \eproof

\mypar{Theorem.}{sa-index-agree}
\emph{Suppose\sss that\trs $X$\sss is\dss paracompact.\oss
If\qss $\mathbb{A}$\sss is\dss a strictly\trs Fredholm\dss family,\oss
then\sss the\sss two\sss index\sss maps\sss
$X\qff \ttoo\qff \Omega\trf \num{\mathcal{S}}$\sss
are homotopic.\oss}

\proof
Theorem\qss \ref{index-replacement}\qss implies\sss that\dss it\dss is\dss 
sufficient\sss to prove\sss the\sss theorem\sss for a\sss
finite-polarized\sss replacement\sss $\mathbb{A}'$\sss
in\sss the role of\sss $\mathbb{A}$\nnsp.\oss
The family\sss $\mathbb{A}'$\sss
is\dss a strictly\trs Fredholm.\oss
By\trs Theorem\qss \ref{adapted-paracompact}\qss there exists a strictly\sss adapted\sss
to\sss $\mathbb{A}'$\sss trivialization of\dss the bundle\sss $\mathbb{H}$\nnsp.\oss
As explained\sss in\dss Section\qss \ref{analytic-index-strict},\oss
this allows\sss to\sss consider\sss $\mathbb{A}'$\sss as a family
of\dss self-adjoint\trs Fredholm\dss operators in\sss the\dss Hilbert\sss space\sss $H$\nnsp,\oss
i.e.\qss as a map\sss
$\mathbb{A}'\dff \colon\dff
X\qff \ttoo\qff \hat{\mathcal{F}}$\dnsp.\oss
We claim\sss that\sss this map\dss is\dss a continuous map\sss to space\sss
$\hat{\mathcal{F}}$\sss considered\sss with\sss the norm\sss topology.\oss

We may\sss assume\sss that\sss the pairs\sss 
$(\trf U\fff,\qff \varepsilon_{\trf U}\dff)$\sss
used\sss for\sss the construction of\sss $\mathbb{A}'$\sss
are strictly adapted\sss to\sss $\mathbb{A}$\nnsp.\oss
Let\sss $z\qff \in\qff X$\nnsp,\oss
and\sss let\sss
$\varepsilon\off =\off \max_{\trf U\qff \ni\qff z}\qff \varepsilon_{\trf U}$\nsp,\oss
where\sss the maximum\dss is\dss taken over\sss $U$\sss 
such\sss that\sss $x\qff \in\qff U$\dnsp.\oss
Then\sss $\varepsilon\off =\off \varepsilon_{\trf U\dff(\trf z\trf)}$\sss for some\sss 
set\sss $U\dff(\trf z\trf)$\sss from\sss the covering,\oss
and\sss $(\trf A_{\dff x}\dff,\qff \varepsilon\trf)$\sss is\dss an enhanced
operator\sss for $x\qff \in\qff U\dff(\trf z\trf)$\nnsp.\oss
Moreover,\pss $(\trf U\dff(\trf z\trf)\fff,\qff \varepsilon\trf)$\sss
is\dss strictly adapted\sss to\sss $\mathbb{A}$\nnsp.\oss 

If\dss $z\qff \not\in\qff U$\sss for some set\sss $U$\sss form\sss the covering,\oss
then\sss $z$\sss does not\sss belong\sss to\sss the support\sss of\sss $r_{\trf U}$\nnsp.\oss
Since\sss the covering\sss by\sss the sets\sss $U$\sss is\dss locally\sss finite,\oss
in\sss some neighborhood\sss $V$\sss of\sss $z$\sss one can\sss replace
in\sss the definition of\sss $r\trf(\trf x\trf)$\sss the sum over all\sss sets\sss $U$\sss
by\sss the sum over\sss set\sss $U$\sss containing\sss $z$\nnsp.\oss
Since\sss
$\varepsilon\off =\off \max_{\trf U\qff \ni\qff z}\qff \varepsilon_{\trf U}$\nsp,\oss
this implies\sss that\sss
$\varepsilon\qff \geq\qff r\trf(\trf x\trf)$\sss
for\sss $x\qff \in\qff V$\nnsp.\oss
Replacing\sss $V$\sss by\sss $V\dff \cap\dff U\dff(\trf z\trf)$\sss if\dss necessary,\oss
we can assume\sss that\sss $V\qff \subset\qff U\dff(\trf z\trf)$\nnsp.\oss
Then\sss $(\trf V\fff,\qff \varepsilon\trf)$\sss
is\dss strictly adapted\sss to\sss $\mathbb{A}$\sss
and every\sss vector\sss $v\qff \in\qff H_{\dff x}$\sss with\sss 
$x\off =\off x\trf(\trf v\trf)\qff \in\qff V$\sss
can\sss be presented as\sss the direct\sss sum\dss\vspace{0pt}
\[
\quad
v
\off =\off
v_{\dff -}\qff \oplus\qff v_{\dff 0}\qff \oplus\qff v_{\dff +}
\qff,\quad
\mbox{where}\quad
\]

\vspace{-36pt}\vspace{-0.75pt}
\[
\quad
v_{\dff -}
\off =\off
P_{\qff (\dff -\qff \infty\fff,\qff -\qff \varepsilon\trf]}\qff
\left(\trf A_{\dff x}\trf\right)\qff(\trf v\trf)\dff,\off
\quad
v_{\dff 0}
\off =\off
P_{\qff [\dff -\qff \varepsilon\fff,\qff \varepsilon\trf]}\qff
\left(\trf A_{\dff x}\trf\right)\qff(\trf v\trf)\dff,\off
\quad
v_{\dff +}
\off =\off
P_{\qff [\dff \varepsilon\fff,\qff \infty\dff)}\qff
\left(\trf A_{\dff x}\trf\right)\qff(\trf v\trf)\dff.\off
\qff
\]

\vspace{-12pt}\vspace{1.5pt}
Since\sss $\varepsilon\qff \geq\qff r\trf(\trf x\trf)$\sss
for\sss $x\qff \in\qff V$\nnsp,\oss\vspace{1.5pt}
\[
\quad
A'_{\dff x}\trf(\trf v_{\dff +}\trf)
\off =\off
v_{\dff +}
\quad
\mbox{and}\quad\dff
A'_{\dff x}\trf(\trf v_{\dff -}\trf)
\off =\off
-\qff v_{\dff -}
\]

\vspace{-12pt}\vspace{1.5pt}
for\sss $x\qff \in\qff V$\nnsp.\oss
It\sss follows\sss that\vspace{1.5pt}
\[
\quad
A'_{\dff x}
\off =\off
-\qff 
P_{\qff (\dff -\qff \infty\fff,\qff -\qff \varepsilon\trf]}\qff
\left(\trf A_{\dff x}\trf\right)
\off +\off
A'_{\dff x}\qff \circ\qff
P_{\qff [\dff -\qff \varepsilon\fff,\qff \varepsilon\trf]}\qff
\left(\trf A_{\dff x}\trf\right)
\off +\off
P_{\qff [\dff \varepsilon\fff,\qff \infty\dff)}\qff
\left(\trf A_{\dff x}\trf\right)
\qff.
\]

\vspace{-12pt}\vspace{1.5pt}
Since\sss $(\trf V\fff,\qff \varepsilon\trf)$\sss 
is\dss adapted\sss to\sss $\mathbb{A}$\nnsp,\oss
the projection\sss 
$P_{\qff [\dff -\qff \varepsilon\fff,\qff \varepsilon\trf]}\qff
\left(\trf A_{\dff x}\trf\right)$\sss
has\sss finite rank\sss and continuously depends on\sss $x$\sss
for $x\qff \in\qff V$\nnsp.\oss
Therefore\sss 
$A'_{\dff x}\qff \circ\qff
P_{\qff [\dff -\qff \varepsilon\fff,\qff \varepsilon\trf]}\qff
\left(\trf A_{\dff x}\trf\right)$\sss also continuously depends on\sss $x$\sss
for $x\qff \in\qff V$\nnsp.\oss 
Since\sss $(\trf V\fff,\qff \varepsilon\trf)$\sss is,\oss moreover,\oss 
strictly adapted\sss to\sss $\mathbb{A}$\nnsp,\oss
the projections\vspace{1.5pt}
\[
\quad
P_{\qff (\dff -\qff \infty\fff,\qff -\qff \varepsilon\trf]}\qff
\left(\trf A_{\dff x}\trf\right)
\quad
\mbox{and}\quad\dff
P_{\qff [\dff \varepsilon\fff,\qff \infty\dff)}\qff
\left(\trf A_{\dff x}\trf\right)
\]

\vspace{-12pt}\vspace{1.5pt}
continuously depend on\sss $x\qff \in\qff V$\sss in\sss the norm\sss topology.\oss
It\sss follows\sss that\sss $A'_{\dff x}$\sss
continuously depend on\sss $x\qff \in\qff V$\sss in\sss the norm\sss topology.\oss
Since\sss the point\sss $x\qff \in\qff X$\sss was arbitrary,\oss
this proves\sss that\sss the map\sss
$\mathbb{A}'\dff \colon\dff
x\off \longmapsto\off A'_{\dff x}$\sss
is\dss continuous as a map\sss to\sss
$\hat{\mathcal{F}}$\sss
with\sss the norm\sss topology.\oss

The classical\sss analytical\dss index of\sss $\mathbb{A}'$\sss
can\sss be considered as\sss the homotopy class of\dss the composition of\sss $\mathbb{A}'$\sss
with\sss the\dss Atiyah--Singer\dss map\sss 
$\alpha\dff \colon\dff
\hat{\mathcal{F}}\qff \ttoo\qff \Omega\trf \mathcal{F}$\sss
defined\sss by\sss the formula\qss (\ref{as-loop})\qss with omitted\sss parameter $x$\nnsp,\oss
and\sss then with\sss the canonical\sss homotopy equivalence\sss
$\Omega\trf \mathcal{F}\qff \ttoo\qff \Omega\trf \num{\mathcal{S}}$\nnsp.\oss
The analytical\dss index of\sss $\mathbb{A}'$\sss in\sss the sense of\trs
Section\qss \ref{analytic-index-section}\qss
can\sss be considered as\sss the homotopy class of\dss the composition of\sss $\mathbb{A}'$\sss
with\sss the homotopy equivalences\sss
$\hat{\mathcal{F}}
\qff \ttoo\qff
\num{\hat{\mathcal{S}}} 
\qff \ttoo\qff
\Omega\trf \num{\mathcal{S}}$\sss
from\dss Section\qss \ref{analytic-index-section}\qss and\trs Theorem\qss \ref{s-loops-shat}.\oss
Using\sss the  homotopy equivalence\sss
$\Omega\trf \mathcal{F}\qff \ttoo\qff \Omega\trf \num{\mathcal{S}}$\sss 
one can\sss interpret\sss this analytical\dss index as\sss the homotopy class
of\dss the composion of\sss $\mathbb{A}'$\sss
with\sss a canonical\sss homotopy equvivalence\sss
$\bm{\alpha}\dff \colon\dff
\hat{\mathcal{F}}
\qff \ttoo\qff
\Omega\trf \mathcal{F}$\dnsp.\oss
By\qss \cite{i2},\oss Theorem\qss 16.6,\oss
the maps\sss $\alpha$\sss and\dss $\bm{\alpha}$\sss are homotopic.\oss
It\sss follows\sss that\sss the maps\sss
$\alpha\dff \circ\dff \mathbb{A}'$\sss and\dss $\bm{\alpha}\dff \circ\dff \mathbb{A}'$\sss
are homotopic,\oss as also\sss their compositions with\sss
$\Omega\trf \mathcal{F}\qff \ttoo\qff \Omega\trf \num{\mathcal{S}}$\nnsp.\oss
It\sss follows\sss that\sss the\sss two index\sss maps are homotopic.\oss  \eproof

\newpage
\myappend{Polarizations\qss and\qss trivializations\qss of\qss Hilbert\qss bundles}{polarizations}

\myuppar{Hilbert\dss bundles and\sss local\dss trivializations.}
Let\sss $X$\sss be a paracompact\sss space and\dss
let\sss $H$\sss be a separable infinite dimensional\dss Hilbert\sss space.\oss
Let\sss $\mathbb{H}$\sss be a\sss locally\sss trivial\dss Hilbert\dss 
bundle over $X$ with fibers isomorphic\sss to $H$\nnsp,\oss
thought\sss as a family\sss
$H_{\dff x}\dff,\qff x\qff \in\qff X$ 
of\dss Hilbert\sss spaces parameterized\sss by $X$\nnsp.\oss
Recall\dss that\sss a\qss 
\emph{local\dss trivialization}\pss of\dss $\mathbb{H}$\sss 
can\sss be considered as a 
family of\dss Hilbert\sss space isomorphisms\sss
$t_{\trf U}\dff(\dff x\trf)\dff \colon\dff H_{\dff x}\qff \ttoo\qff H$\nnsp,\oss 
where $x$ runs over an open subset\qss $U\qff \subset\qff Y$\nnsp.\oss
One can speak also about\sss trivializations over an arbitrary subset\sss of\sss $X$\nnsp.\oss

\myuppar{Polarizations.}
A\qss \emph{polarization}\qss of\dss a\sss Hilbert\sss space $K$\sss is\dss
a presentation of\sss $K$\sss as an orthogonal\sss direct\sss sum\sss
$K\off =\off K_{\dff -}\dff \oplus\dff K_{\dff +}$\sss of\dss two
infinitely dimensional\sss closed subspaces\sss $K_{\dff -}\dff,\qff K_{\dff +}$\nsp.\oss
Clearly,\oss a polarization\sss $K\off =\off K_{\dff -}\dff \oplus\dff K_{\dff +}$\sss
is\dss uniquely determined\sss by\sss $K_{\dff +}$\nsp,\oss
and a\sss polarization of\sss $K$ can\sss be also defined as closed\sss infinitely dimensional\sss
subspace $P\qff \subset\qff K$\sss such\sss that\sss the orthogonal\sss complement\sss
$K\dff \ominus\dff P$\sss is\dss also infinitely dimensional.\oss
The set\sss of\dss polarizations of\sss $K$\sss is\dss equipped\sss by\sss the\sss
the\sss topology defined\sss by\sss the norm\sss topology of\dss orthogonal\sss
projections\sss $K\qff \ttoo\qff P$\dnsp.\oss
By a well\dss known\sss theorem of\trs Atiyah\dss and\dss Singer\qss \cite{as}\qss
the space of\dss polarizations\dss is\dss contractible.\oss

A\qss \emph{local\dss polarization}\qss of\sss $\mathbb{H}$\sss over an open set\sss
$U\qff \subset\qff X$\sss is\dss a family of\dss polarizations\sss
$P_{\fff x}\qff \subset\qff H_{\dff x}$\nsp, $x\qff \in\qff U$\sss
such\sss that\sss for some\sss trivialization\sss $t_{\trf U}$
of\sss $\mathbb{H}$ over $U$\sss  the family\sss
$t_{\trf U}\fff(\dff x\trf)\dff(\trf P_{\fff x}\trf)$\nnsp, $x\qff \in\qff U$
of\dss polarizations of\sss $H$\sss is\dss norm\sss continuous.\oss 
Two\sss local\sss polarizations\sss
$P_{\fff x}\qff \subset\qff H_{\dff x}$\nsp, $x\qff \in\qff U$\sss and\sss
$P\fff'_x\qff \subset\qff H_{\dff x}$\nsp, $x\qff \in\qff U'$\sss
are said\sss to be\qss \emph{compatible at\sss a point}\qss
$z\qff \in\qff U\dff \cap\dff U'$\sss if\dss 
there exists a neighborhood\sss 
$V\qff \subset\qff U\dff \cap\dff U'$\sss of\sss $z$
such\sss that\sss either\sss 
$P_{\fff x}$\sss is\dss a subspace of\dss finite codimension\sss in\sss $P'_{x}$ 
for every\sss $x\qff \in\qff V$\sss and\sss the family of\dss orthogonal\sss complements\sss 
$P\fff'_x\dff \ominus\dff P_{\dff x}$\nsp, $x\qff \in\qff V$\sss
is\dss continuous,\oss or\sss this condition\sss holds with\sss the roles of\dss
$P_{\fff x}$\sss and\sss $P'_{x}$ interchanged.\oss 
Two\sss local\sss polarizations\sss
$P_{\fff x}\qff \subset\qff H_{\dff x}$\nsp, $x\qff \in\qff U$\sss and\sss
$P\fff'_x\qff \subset\qff H_{\dff x}$\nsp, $x\qff \in\qff U'$\sss
are said\sss to be\qss \emph{compatible}\pss if\dss they are compatible
at\sss every\sss point\sss $z\qff \in\qff U\dff \cap\dff U'$\dnsp.\oss

A\qss \emph{polarization atlas}\qss is\dss a collection of\dss
pair-wise compatible\sss local\sss polarizations over some open subsets covering $X$\nnsp.\oss
Two polarization atlases are said\sss to be\qss \emph{equivalent}\pss if\dss their union\dss
is\dss also a polarization atlas.\oss
A\qss \emph{polarization}\qss of\sss $\mathbb{H}$\sss is\dss an equivalence class
of\dss polarization atlases.\oss
Clearly,\oss every polarization\dss is\dss represented\sss by unique maximal\sss atlas.\oss
A\sss local\sss polarization\dss is\dss said\sss to be a\qss \emph{chart}\pss
of\dss a polarization\sss $\mathbb{P}$\sss if\dss it\dss belongs\sss 
to\sss the maximal\sss atlas of\dss $\mathbb{P}$\dnsp.\oss
If\sss $\mathbb{P}$\sss is\dss a polarization of\sss $\mathbb{H}$\sss
and\sss $Y\qff \subset\qff X$\nnsp,\oss
then\sss $\mathbb{P}$\sss defines,\oss in an obvious way,\oss
a polarization of\dss the restriction\sss $\mathbb{H}\trf |\dff Y$\sss
of\sss $\mathbb{H}$\sss to\sss $Y$\dnsp,\oss
which we will\sss denote\sss by\sss $\mathbb{P}\dff |\dff Y$\dnsp.\oss

\myuppar{Adapted\sss trivializations.}
Let\sss $\mathbb{P}$\sss be a polarization of\sss $\mathbb{H}$\nnsp,\oss
and\sss $P_{\fff x}$\nsp, $x\qff \in\qff U$\sss
be a chart\sss of\sss $\mathbb{P}$\dnsp.\oss
A\sss trivialization\sss $t_{\trf V}$ of\dss $\mathbb{H}$\sss 
defined over a subset\sss $V\qff \subset\qff X$\sss 
is\dss said\sss to be\qss \emph{adapted\dss to the chart}\qss
$P_{\fff x}$\nsp, $x\qff \in\qff U$\sss
if\dss the family of\dss polarizations\sss
$t_{\qff V}\fff(\dff x\trf)\dff(\trf P_{\fff x}\trf)$\nnsp, $x\qff \in\qff U\dff \cap\dff V$
of\sss $H$\sss is\dss norm\sss continuous.\oss
A\sss trivialization\sss $t_{\qff V}$\sss is\dss said\sss 
to\sss be\qss \emph{adapted\dss to\sss the polarization}\qss $\mathbb{P}$\sss
if\sss $t_{\qff V}$\sss is\dss adapted\sss to every chart\sss of\sss $\mathbb{P}$\dnsp.\oss
The following\sss properties of\dss adapted\sss trivializations\sss immediately\sss
follow\sss from\sss the definitions and\sss will\sss be used\sss without\sss references.\oss

\myapar{Lemma.}{basic}
\emph{Let\trs $\mathbb{P}$\sss be a polarization of\qss $\mathbb{H}$\nnsp,\oss
and\trs $P_{\fff x}$\nsp, $x\qff \in\qff U$\sss and\sss
$P'_{x}$\nsp, $x\qff \in\qff U'$\sss
be some charts of\qss $\mathbb{P}$\dnsp.\oss
If\pss $V\qff \subset\qff U$\sss and a\sss trivialization\sss 
$t_{\qff V}\fff(\dff x\trf)$\nnsp, $x\qff \in\qff V$\sss 
is\dss adapted\sss to\dss
$P_{\fff x}$\nsp, $x\qff \in\qff U$\dnsp,\oss
then\sss it\dss is\dss adapted\dss to\dss
$P'_{x}$\nsp, $x\qff \in\qff U'$\dnsp.\oss
If\qss a\sss trivialization\sss 
$t_{\qff V}\fff(\dff x\trf)$\nnsp, $x\qff \in\qff V$\sss 
is\dss adapted\sss to every chart\sss from
some atlas of\qss $\mathbb{P}$\dnsp,\oss then\sss $t_{\qff V}$\sss is\dss
adapted\dss to $\mathbb{P}$\dnsp.\oss}  \eproof

\myuppar{Existence of\dss adapted\dss trivializations.}
Recall\dss that\sss every\dss Hilbert\dss bundle over a paracompact\sss space\dss
is\dss trivial.\oss
In\sss particular,\oss there exists a\sss trivialization of\sss $\mathbb{H}$\nnsp,\oss
i.e.\qss a\sss local\sss trivialization defined over\sss the whole space $X$\nnsp.\oss
We will\sss strengthen\sss this result\sss by proving\sss that\sss
for every polarization\sss $\mathbb{P}$\sss of\sss $\mathbb{H}$\sss
there exists a\sss trivialization of\sss $\mathbb{H}$ adapted\dss to $\mathbb{P}$\dnsp.\oss
We will\sss consider\sss first\sss the simpler case of\dss triangulable spaces $X$\nnsp.\oss

\myapar{Theorem.}{p-adapted}
\emph{Let\pss $\mathbb{P}$\sss be a polarization of\dss $\mathbb{H}$\nnsp.\oss
If\qss $X$ is\dss a\sss triangulable space,\oss
then\sss there exists a\sss trivialization of\trs 
$\mathbb{H}$\dss 
adapted\dss to\sss $\mathbb{P}$\dnsp.\oss}

\proof
Let\sss us\sss choose an atlas of\dss charts\sss
$P_{\fff x}^{\fff U}\qff \subset\qff H_{\dff x}$\nsp, $x\qff \in\qff U$\sss
defining\sss $\mathbb{P}$\dnsp.\oss
Replacing\sss a\sss triangulation of\sss $X$\sss
by a subdivision we can assume\sss that\sss every\sss simplex $\sigma$ of\dss
a fixed\dss triangulation of\sss $X$\sss is\dss contained\sss in\sss
one of\dss the sets\sss $U$\dnsp.\oss 
Let\sss $\ske_{\trf n}\dff X$\sss
be\sss the $n${\nnsp}th skeleton of\dss this\sss triangulation.\oss

Let\sss us\sss begin\sss by\sss choosing isomorphisms\sss
$t\trf(\dff v\trf)\dff \colon\dff H_{\dff v}\qff \ttoo\qff H$\sss
for\sss the vertices $v$ of\dss the\sss triangulation.\oss
Such\sss isomorphisms define a\sss trivialization of\dss the restriction\sss
$\mathbb{H}\trf |\fff \ske_{\trf 0}\dff X$\sss of\sss $\mathbb{H}$\nnsp.\oss
This\sss trivialization\dss is\dss obviously adapted\sss to\sss
$\mathbb{P}\trf |\fff \ske_{\trf 0}\dff X$\nnsp.\oss

Suppose\sss that\sss we already\sss found a\sss trivialization $t$
of\sss $\mathbb{H}\trf |\fff \ske_{\trf n}\dff X$ 
adapted\sss to\sss $\mathbb{P}\trf |\fff \ske_{\trf n}\dff X$\nnsp.\oss
Let\sss $\sigma$\sss be an $(\dff n\qff +\qff 1\dff)$\dnsp-simplex of\dss
our\sss triangulation,\oss and\sss let\sss $\partial\dff \sigma$ be its boundary.\oss
By\sss our assumptions\sss there exists a\sss chart\sss
$P_{\fff x}^{\fff U}$\nsp, $x\qff \in\qff U$\sss
such\sss that\sss $\sigma\qff \subset\qff U$\dnsp.\oss
By\sss the definition\sss there exists a\sss local\sss trivialization $t_{\trf U}$ of\sss $\mathbb{H}$\sss
over $U$ adapted\sss to\sss 
$P_{\fff x}^{\fff U}$\nsp, $x\qff \in\qff U$\nnsp.\oss
Then\sss the families\vspace{1.5pt}
\[
\quad
t\trf\left(\trf P_{\fff x}^{\fff U}\trf\right)\dff,\off x\qff \in\qff \partial\dff \sigma
\quad
\mbox{and}\quad
t_{\trf U}\trf\left(\trf P_{\fff x}^{\fff U}\trf\right)\dff,\off x\qff \in\qff \partial\dff \sigma
\qff
\]

\vspace{-12pt}\vspace{1.5pt}
of\dss polarizations of\sss $H$\sss are norm continuous.\oss
The contractibility\sss of\dss the space of\dss polarizations of\sss $H$\sss implies\sss
that\sss the second\sss family\dss is\dss homotopic\sss to\sss the first\sss one.\oss
The second\sss family\dss is\dss defined also for\sss
$x\qff \in\qff U$\nnsp,\oss
and we can extend\sss this homotopy\sss to $U$\nnsp.\oss
The\sss last\sss family\sss in\sss the extended\sss homotopy\dss is\dss
a new\sss local\dss trivialization\sss $t\fff'_{\dff U}$\sss
adapted\sss to\sss $P_{\fff x}^{\fff U}$\nsp, $x\qff \in\qff U$\sss
and such\sss that\sss\vspace{0pt}
\[
\quad
t\trf(\trf x\trf)\qff\left(\trf P_{\fff x}^{\fff U}\trf\right)
\off =\off
t\fff'_{\dff U}\trf(\trf x\trf)\trf\left(\trf P_{\fff x}^{\fff U}\trf\right)
\]

\vspace{-12pt}\vspace{0pt}
for\sss every\sss $x\qff \in\qff \partial\dff \sigma$\nsp.\oss
The family\sss $P_{\fff x}^{\fff U}$\nsp, $x\qff \in\qff \partial\dff \sigma$\sss
is\dss a\dss Hilbert\dss subbundle of\dss $\mathbb{H}\trf |\trf \partial\dff \sigma$\nnsp.\oss
It\dss is\dss trivial\sss as a\sss Hilbert\dss bundle,\oss
and\sss its orthogonal\sss complement\sss
$H_{\dff x}\qff \ominus\qff P_{\fff x}^{\fff U}$\nsp, $x\qff \in\qff \partial\dff \sigma$\sss
is\dss also\sss trivial.\oss
By using some\sss trivializations of\dss these bundles we can consider\sss the maps\vspace{2pt}
\[
\quad
t\fff'_{\dff U}\trf(\trf x\trf)^{\dff -\dff 1}\dff \circ\dff t\trf(\trf x\trf)\dff,\off 
x\qff \in\qff \partial\dff \sigma
\]

\vspace{-12pt}\vspace{2pt}
as isomorphisms of\dss a\dss Hilbert\sss space $K$ respecting a
polarization\sss $K\off =\off K_{\dff -}\dff \oplus\dff K_{\dff +}$\nsp.\oss
The contractibility of\dss the unitary groups\sss of\dss Hilbert\sss spaces
in\sss the compact-open\sss topology\qss
(see\dss Atiyah\dss and\dss Segal\qss \cite{ase})\qss 
implies\sss that\sss this
family\dss is\dss homotopic\sss by families of\dss maps respecting\sss the polarization\sss 
to\sss the family of\dss identity maps.\oss
By\sss taking\sss the compositions with\sss
$t\fff'_{\dff U}\dff(\trf x\trf)\dff,\pff 
x\qff \in\qff \partial\dff \sigma$\sss
we get\sss a homotopy\sss from\sss the\sss trivialization\sss
$t\fff'_{\dff U}\trf(\trf x\trf)\dff,\off 
x\qff \in\qff \partial\dff \sigma$\sss
to\sss the\sss trivialization\sss
$t\trf(\trf x\trf)\dff,\off 
x\qff \in\qff \partial\dff \sigma$\sss
by\dss trivializations adapted\dss to\sss
$P_{\fff x}^{\fff U}$\nsp, $x\qff \in\qff U$\dnsp.\oss
By extending\sss this homotopy\sss to\sss $x\qff \in\qff \sigma$\sss
we get\sss a\sss trivialization\sss 
$t\fff''_{\dff \sigma}\trf(\trf x\trf)\fff,\qff x\qff \in\qff \sigma$\sss
adapted\sss to\sss $P_{\fff x}^{\fff U}$\nsp, $x\qff \in\qff U$ 
and such\sss that\sss \vspace{0pt}
\[
\quad
t\fff''_{\dff \sigma}\trf(\trf x\trf)
\off =\off
t\trf(\trf x\trf)
\]

\vspace{-12pt}\vspace{0pt}
for every\sss $x\qff \in\qff \partial\dff \sigma$.\oss
Let\sss us extend $t$\dss from\sss $\ske_{\trf n}\dff X$\sss
to\sss $\ske_{\trf n}\dff X\dff \cup\dff \sigma$\sss by\sss $t\fff''_{\dff U}$\sss
and denote\sss the resulting\sss trivialization again\sss by $t$\nnsp.\oss
Then\sss the restriction of\sss $t$\sss to\sss $\sigma$\sss
is\dss  adapted\sss to\sss
$P_{\fff x}^{\fff U}$\nsp, $x\qff \in\qff U$\dnsp.\oss
It\sss follows\sss that\sss $t$\sss is\dss adapted\sss to\sss
$\mathbb{P}\trf |\fff \ske_{\trf n}\dff X\dff \cup\dff \sigma$\nnsp.\oss
By doing\sss this for all $(\dff n\qff +\qff 1\dff)$\dnsp-simplices simultaneously,\oss
we can extend\sss $t$\sss to a\sss trivialization of\dss
$\mathbb{H}\trf |\fff \ske_{\trf n\dff +\dff 1}\dff X$\sss adapted\dss to\sss
$\mathbb{P}\trf |\fff \ske_{\trf n\dff +\dff 1}\dff X$\nnsp.\oss
By continuing\sss in\sss this way\sss we will\sss get\sss
an\sss adapted\sss trivialization of\sss $H$\nnsp.\oss  \eproof

\myapar{Theorem.}{p-adapted-paracompact}
\emph{Let\pss $\mathbb{P}$\sss be a polarization of\dss $\mathbb{H}$\nnsp.\oss
If\qss $X$ is\dss a\sss paracompact\sss space,\oss
then\sss there exists a\sss trivialization of\trs 
$\mathbb{H}$\dss 
adapted\dss to\sss $\mathbb{P}$\dnsp.\oss}

\proof
Let\sss us\sss choose an atlas of\dss charts\sss
$P_{\fff x}^{\fff a}$\nsp, $x\qff \in\qff U_{\dff a}$\sss
defining\sss $\mathbb{P}$\dnsp,\oss
where $a$ runs over some set\sss $\Sigma$\nnsp,\oss
and\dss let\sss us choose\sss for every\sss $a\qff \in\qff \Sigma$\sss
a\sss local\dss trivializations\sss $t_{\trf a}$\sss of\sss $\mathbb{H}$ over\sss $U_{\dff a}$\sss
such\sss that\sss the family\sss 
$t_{\trf a}\fff(\dff x\trf)\dff(\trf P_{\fff x}^{\fff a}\trf)$\nnsp, $x\qff \in\qff U_{\dff a}$\sss
is\dss norm continuous.\oss

Let\sss $\Sigma\ffin$\sss and\sss $U_{\dff \sigma}\dff,\qff \varepsilon_{\dff \sigma}$\sss 
for\sss $\sigma\qff \in\qff \Sigma\ffin$\sss
be defined as in\dss Section\qss \ref{analytic-index-section},\oss 
and\sss let\sss $X_{\dff U}$\sss
be\sss the\sss topological\sss category\sss from\dss Section\qss \ref{analytic-index-section}.\oss
The space of\dss objects of\sss $X_{\dff U}$\sss is\dss the disjoint\sss union\sss
of\dss subspaces\sss $U_{\dff \sigma}$\nsp,\oss
and\dss the category\sss $X_{\dff U}$\sss is\dss defined\sss by an order\sss $\leq$\sss
on\sss $\ob\trf X_{\dff U}$\nsp.\oss
Clearly,\oss the set\sss pairs\sss $(\trf u\fff,\qff u\trf)$\nnsp,\oss
where\sss $u\qff \in\qff \ob\trf X_{\dff U}$\nsp,\oss
is\dss the union of\dss several\sss components of\dss the set\sss of\dss pairs\sss
$(\trf u\fff,\qff v\trf)$\sss of\dss comparable objects of\sss $X_{\dff U}$\nsp,\oss
i.e.\qss such\sss that\sss either\sss $u\qff \leq\qff v$\sss or\sss $v\qff \leq\qff u$\nnsp.\oss
In\sss the\sss terminology of\qss \cite{i2}\qss this means\sss that\sss
$\ob\trf X_{\dff U}$\sss is\dss a partially ordered space with\qss 
\emph{free equalities},\oss and\sss hence\sss $X_{\dff U}$\sss
can\sss be\sss treated as a\sss topological\sss simplicial\sss complex.\oss
In\sss particular,\pss $\num{X_{\dff U}}$\sss 
is\dss equal\sss to\sss the\qss ``naive''\qss
geometric realization\sss $\bbnum{X_{\dff U}}$\nnsp.\oss
We refer\sss to\qss \cite{i2},\oss Section\qss 5,\oss for\sss the
definition of\dss the geometric realization\sss 
$\bbnum{\bullet}$ and\dss its basic properties.\oss

Let\sss $\mathbb{H}^{\dff U}$ be\sss the bundle induced\sss from\sss $\mathbb{H}$\sss by\sss 
$\bbnum{\pr}\dff \colon\dff
\bbnum{X_{\dff U}}\qff \ttoo\qff X$\nnsp.\oss
The polarization\sss $\mathbb{P}$\sss induces a polarization\sss $\mathbb{P}^{\fff U}$\sss of\dss
of\sss $\mathbb{H}^{\dff U}$\dnsp.\oss
If\sss $s$\sss is\dss a homotopy\sss inverse of\sss $\num{\pr}$\sss as above,\oss
then\sss $\mathbb{H}$\sss and\sss $\mathbb{P}$\sss are equal\sss to\sss
the bundle and\sss the family of\dss operators induced\sss by\sss $s$\sss 
from\sss $\mathbb{H}^{\dff U}$\sss and\sss $\mathbb{P}^{\dff U}$ respectively.\oss
Therefore\sss it\dss is\dss sufficient\sss to prove\sss that\sss
there exists a\sss trivialization of\sss $\mathbb{H}^{\dff U}$\sss
adapted\sss to\sss $\mathbb{P}^{\dff U}$\dnsp.\oss
As in\sss the proof\dss of\trs Theorem\qss \ref{p-adapted},\oss
we will\sss prove\sss this using an\sss induction by\sss skeletons.\oss
If\dss $\sigma\qff \in\qff \Sigma\ffin$\sss and\sss $a\qff \in\qff \sigma$\nnsp,\oss
then\sss the restriction $t_{\dff \sigma}$\sss of\sss $t_{\dff a}$\sss to\sss 
$U_{\dff \sigma}\qff \subset\qff U_{\dff a}$\sss is\dss adapted\sss to $\mathbb{P}$\nnsp.\oss
Hence\sss the\sss trivializations\sss $t_{\trf \sigma}$\sss
define a\sss trivialization of\dss the restriction of\dss $\mathbb{H}^{\dff U}$\sss
to\sss the $0${\nsp}th skeleton\sss $\bbnum{\nsp\ssk_{\trf 0}\dff X_{\dff U}}$\sss
adapted\sss to\sss the restriction of\sss $\mathbb{P}^{\dff U}$\dnsp.\oss

Suppose\sss that\sss we constructed a\sss trivialization $t$ of\dss the restriction\dss
$\mathbb{H}^{\dff U}\dff |\qff \bbnum{\nsp\ske_{\trf n}\dff X_{\dff U}}$\sss
adapt\-ed\sss to\sss
$P^{\dff U}_{\dff u}\dff,\qff u\qff \in\qff \bbnum{\nsp\ske_{\trf n}\dff X_{\dff U}}$
for some $n\qff \geq\qff 0$\nnsp.\oss
An $(\dff n\qff +\qff 1\dff)$\dnsp-simplex of\dss $X_{\dff U}$\sss
is\dss defined\sss by a point\sss $x\qff \in\qff X$\sss and a strictly
decreasing sequence\vspace{1.75pt}
\[
\quad
\sigma\trf(\dff n\qff +\qff 1\dff)
\off \supsetneq\off
\sigma\trf(\dff n\trf)
\off \supsetneq\off
\ldots
\off \supsetneq\off
\sigma\trf(\trf 0\trf)
\]

\vspace{-12pt}\vspace{1.75pt}
of\dss elements of\sss $\Sigma\ffin$\dnsp.\oss
Therefore\sss the space of\dss $(\dff n\qff +\qff 1\dff)$\dnsp-simplices of\dss $X_{\dff U}$\sss
can\sss be identified\sss with\sss the disjoint\sss union\sss
of\dss the sets\sss $U_{\dff \sigma\trf(\dff n\qff +\qff 1\dff)}$\sss
over\sss all\sss such sequences.\oss
In\sss particular,\oss every component\sss of\dss the space of\sss
$(\dff n\qff +\qff 1\dff)$\dnsp-simplices
can\sss be identified\sss with\sss $U_{\dff \sigma}$\sss 
for some\sss $\sigma\qff \in\qff \Sigma\ffin$\nnsp.\oss
Let\sss us consider\sss such a component\sss $U_{\dff \sigma}$\sss
and\sss the canonical\sss map\sss\vspace{1.75pt}
\begin{equation*}
\quad
e_{\dff \sigma}\dff \colon\dff
\Delta^{n\dff +\dff 1}\dff \times\qff U_{\dff \sigma}
\qff \ttoo\qff
\bbnum{X_{\dff U}}
\off =\off
\num{X_{\dff U}}
\pff.
\end{equation*}

\vspace{-12pt}\vspace{1.75pt}
The bundle $\mathbb{H}^{\dff U}$ and\sss the polarization\sss $\mathbb{P}^{\dff U}$
induce a bundle $\mathbb{H}^{\dff \sigma}$ 
and a\sss polarization\sss $\mathbb{P}^{\dff \sigma}$ over\sss
$\Delta^{n\dff +\dff 1}\dff \times\qff U_{\dff \sigma}$\nsp,\oss
and\sss the already constructed\sss trivialization\sss $t$\sss
defines a\sss trivialization $t'$ of\sss $\mathbb{H}^{\dff \sigma}$ over\trs
$\partial\trf \Delta^{n\dff +\dff 1}\dff \times\qff U_{\dff \sigma}$\dss
adapted\sss $\mathbb{P}^{\dff \sigma}$\dnsp.\oss
Since\sss the map\sss $\num{\pr}$\sss collapses
each simplex\sss to a point,\oss the fibers of\sss
$\mathbb{H}^{\dff U}$ over a simplex can\sss be\sss treated as equal,\oss
as also\sss the polarizations\sss $P^{\dff U}_{\dff u}$\nsp.\oss
Hence\sss the\sss trivialization\sss $t_{\trf \sigma}$\sss
defines a\sss trivialization of\sss $\mathbb{H}^{\dff \sigma}$
adapted\sss to\sss $\mathbb{P}^{\dff \sigma}$\dnsp.\oss

Let\sss us\sss use\sss this\sss trivialization\sss to identify\sss
$\mathbb{H}^{\dff \sigma}$\sss with\sss 
the\sss trivial\sss bundle over\sss
$\Delta^{n\dff +\dff 1}\dff \times\qff U_{\dff \sigma}$\nsp.\oss
Now we can use\sss the contractibility of\dss 
the space of\dss polarizations 
and\dss the contractibility\sss of\sss 
the unitary\sss groups\sss in\sss the compact-open\sss topology\sss
as in\sss the proof\dss of\trs Theorem\qss \ref{adapted}\qss
to extend $t'$\sss to a\sss trivialization $t''$ 
of\sss $\mathbb{H}^{\dff \sigma}$ over\sss
$\Delta^{n\dff +\dff 1}\dff \times\qff U_{\dff \sigma}$\dss
adapted\sss to\sss $\mathbb{P}^{\dff \sigma}$\dnsp.\oss
Since we are dealing\sss with a\sss topological\sss simplicial\sss
complex,\oss the map\sss $e_{\dff \sigma}$\sss is\dss injective
and we can extend $t$\sss to\sss the image of\dss the map\sss $e_{\dff \sigma}$\sss 
by\sss using $t''$\dnsp.\oss
By\sss doing\sss this simultaneously\sss for every component\sss
of\dss the space of\sss $(\dff n\qff +\qff 1\dff)$\dnsp-simplices
we can extend\sss $t$\sss to\sss next\sss skeleton\sss
$\bbnum{\nsp\ske_{\trf n\dff +\dff 1}\dff X_{\dff U}}$\nnsp.\oss
The\sss topology of\sss 
$\num{X_{\dff U}}\off =\off \bbnum{X_{\dff U}}$\sss is,\oss
by\sss the definition,\sss the direct\sss limit\sss of\dss the\sss
topologies of\dss the skeletons.\oss
Therefore,\oss by continuing\sss in\sss this way\sss
we will\sss get\sss a\sss trivialization of\dss $\mathbb{H}^{\dff U}$
adapted\sss to\sss $\mathbb{P}^{\dff U}$\dnsp.\oss  \eproof

\myapar{Theorem.}{p-adapted-relative}
\emph{Under\sss the assumptions of\pss Theorem\qss \ref{p-adapted}\qss or\qss \ref{p-adapted-paracompact},\oss
suppose\sss that\qss $Y$\dss is\dss a subcomplex or closed subset\sss of\pss $X$\sss respectively.\oss
If\qss $t\trf(\trf y\trf)\dff,\qff y\qff \in\qff Y$\sss
is\trs a\sss trivialization of\qss $\mathbb{H}\trf |\trf Y$\sss adapted\dss to\sss 
$\mathbb{P}\trf |\trf Y$\dnsp,\oss
then\sss $t$\sss can\sss be extended\dss to an\sss adapted\dss trivialization of\qss
$\mathbb{P}$\nnsp.\oss}

\proof
The proof\trs is\dss a standard\sss modification of\dss the proofs
of\qss Theorems\qss \ref{p-adapted}\qss and\qss \ref{p-adapted-paracompact}.\oss  \eproof

\myapar{Corollary.}{p-adapted-unique}
\emph{Under\sss the assumptions of\pss Theorems\qss \ref{p-adapted}\qss or\qss \ref{p-adapted-paracompact},\oss
every\sss two\dss trivializations adapted\dss to\sss $\mathbb{P}$\sss are homotopic\sss
in\dss the class of\dss trivializations adapted\dss to\sss $\mathbb{P}$\dnsp.\oss}

\proof
It\dss is\dss sufficient\sss to apply\trs Theorem\qss \ref{p-adapted-relative}\qss
to\sss the subset\sss $X\dff \times\dff \{\trf 0\dff,\qff 1\trf\}$
of\trs $X\dff \times\dff [\dff 0\dff,\qff 1\dff]$\nnsp.\oss  \eproof

\myuppar{Remark.}
M.\dss Prokhorova\qss \cite{p2},\oss influenced\sss by\sss the first\sss version of\dss
the present\sss paper,\oss considered\sss families of\trs local\dss polarizations 
compatible in a weaker sense,\oss 
related\sss to\sss ours as compact\sss operators
are related\sss to operators of\dss finite rank.\oss
Namely,\oss the orthogonal\dss projections onto\sss 
$P_{\fff x}$ and\sss $P\fff'_x$\sss
are required\dss to differ by compact\sss operators
continuously depending on $x$\nnsp.


\vspace{\bigskipamount}

\begin{flushright}

November\qss 28\fff,\oss 2021.\oss
Revised\qss January\qss 2,\oss 2023

https\halfff:/\!/\!nikolaivivanov.com

E-mail\halfff:\oss nikolai.v.ivanov{\fff}@{\dff}icloud.com\vspace{12pt}

Department\sss of\trs Mathematics,\oss Michigan\sss State\sss University

\end{flushright}

\end{document}